\documentclass[11pt]{amsart}
\usepackage{amsmath,amsfonts,amssymb,mathrsfs,graphicx,subfigure}
\usepackage{color,psfrag,graphicx}
\usepackage{algorithm,algpseudocode}

\def\half{\frac{1}{2}}
\def\veps{\varepsilon}

\def\bm{\boldsymbol}

\def\wt{\widetilde}

\def\bp{\bm{p}}
\def\bq{\bm{q}}
\def\bg{\bm{g}}
\def\bu{\bm{u}}
\def\by{\bm{y}}

\def\bmf{\bm{f}}
\def\bx{\bm{x}}

\def\R{\mathbb R}
\def\T{\mathcal T}
\def\V{\mathcal V}
\def\argmin{\operatorname{argmin}}

\def\eg{{\it e.g.}}

\begin {document}

\title[An atomistic/continuum coupling method using enriched bases]
{An atomistic/continuum coupling method using enriched bases
\vrule height
15pt width 0pt}

\author[J. Chen]{Jingrun Chen}
\address{Mathematics Department, South Hall 6705, University of California, Santa Barbara, CA93106, USA} \email{cjr@math.ucsb.edu}
\thanks{}

\author[C. J. Garc\'{\i}a-Cervera]{Carlos J. Garc\'{\i}a-Cervera}
\address{Mathematics Department, South Hall 6707, University of California, Santa Barbara, CA93106, USA} \email{cgarcia@math.ucsb.edu}
\thanks{}

\author[X. Li]{Xiantao Li}
\address{Department of Mathematics, the Pennsylvania State University, University Park, Pennsylvania, 16802, U.S.A.} \email{xli@math.psu.edu}
\thanks{}

\date{\today}
\begin {abstract}
A common observation from an atomistic to continuum coupling method is that the error is often generated and concentrated near the interface, where the two models are combined.
  In this paper, a new method is proposed to suppress the error at the interface,
and as a consequence, the overall accuracy is improved.
The method is motivated by formulating the molecular mechanics model as a two-stage minimization problem.
In particular, it is demonstrated that the error at the interface can be considerably reduced
when new basis functions are introduced in a Galerkin projection formalism.  The  improvement of the
accuracy is illustrated by two examples. Further, comparison with some quasicontinuum-type methods is provided.

\end {abstract}

\maketitle

\section{Introduction}
\label{sec:introduction}
Multiscale modeling has been becoming an extremely important tool in many areas of applied sciences
 \cite{E:book}. By combining models from different spatial and temporal scales, one can describe and simulate complex physical processes with balanced accuracy and efficiency.

 For the modeling of mechanical properties of crystalline materials, one notable success is the  quasicontinuum (QC) method \cite{TadmorOrtizPhillips:1996}. In addition to its applications to
 numerous nano-mechanical systems, it also serves as a great example, where a conventional continuum (elasticity) model can be supplemented with an atomic-level description to appropriately incorporate
the needed microscopic structures. Since the first emergence of the QC method, there have been a lot
of recent developments aimed at understanding and improving its modeling accuracy \cite{MillerTadmor:2009, LuskinOrtner:2013}.  In particular, the Cauchy-Born approximation  \cite{BornHuang:1954, Ericksen:1984} in the QC
method, which evaluates the elastic energy or the strain-stress relation based on the underlying lattice structures and atomic interactions, has been analyzed in \cite{Blancetal:2002, EMing:2007}.

Another important issue in such multiscale methods is the consistency of the specific strategy for combining the two models.
This is challenging, especially because the two models are of different nature: The elasticity
model is continuous, and is posed in terms of strain and stress (or energy), while the atomistic model
is expressed using displacements and interatomic forces. The inconsistency is manifested, for instance,
by the appearance of ghost forces \cite{Shenoyetal:1999}. Explicit error estimates for the influence of the ghost forces can be found in \cite{MingYang:2009, ChenMing:2012, CuiMing:2013}
for models in one, two, and three dimensions. Meanwhile, a number of methods have been proposed to remove
ghost forces. These include the force correction by ``deadload" \cite{Shenoyetal:1999}, the
quasi-nonlocal QC method \cite{Shimokawa:2004}, and the geometrically consistent scheme \cite{ELuYang:2006}.
More recent efforts include
the energy-based atomistic-to-continuum coupling methods \cite{Shapeev:2011, Shapeev:2012, OrtnerZhang:2012, OrtnerZhang},
and the energy-based blended QC method \cite{VanLuskin:2011}, \textit{etc}.
Nevertheless, in high dimensions when the interface is of general geometry, removing ghost forces
in an energy-based method remains an open challenge.
Another approach is to couple the two models based on forces, \eg, \cite{Shenoyetal:1999, MillerTadmor:2002,
MillerTadmor:2009,Lietal:2013,LI12a,citeulike:2898408}.

From the perspective of numerical analysis, how the aforementioned methods approximate a solution of the atomistic model has been
of great interest over a decade; see \cite{LuskinOrtner:2013} for a review and references therein.
To address the issue of convergence, one needs to prove consistency (free of ghost forces) and some
stability conditions. Thorough studies have been conducted in the one-dimensional case.
Under some technical conditions, it has been proven that the quasi-nonlocal QC method converges linearly in $W^{1,\infty}$ norm \cite{MingYang:2009}. Its convergence rate is $\veps^{1+1/p}$ if the $W^{1,p}$ norm is used  with $1\leq p\leq \infty$ \cite{DobsonLuskin:2009b}. It has also been proven that the geometrically consistent QC
method converges linearly in $W^{1,\infty}$ norm \cite{MingYang:2009}. Convergence of the force-based QC method can be found in \cite{Ming:2008,DobsonLuskinOrtner:2010a} with a second-order accuracy \cite{DobsonLuskinOrtner:2010a}.
There are results in high dimensions as well, \eg, \cite{LuMing:2013,LuMing:2014} for a force-based method, and \cite{OrtnerShapeev}
for an energy-based atomistic/continuum method.
Also of practical importance is the adaptive approach based on a systematic error control  \cite{OdPrRoBa06,prudhomme2006error,prudhomme2009adaptive}.
In addition, convergence of the QC method has also been investigated using the notion of $\Gamma-$convergence \cite{schaffner2013analytical}.


Numerous numerical results indicate that the error, especially the deformation gradient, exhibits peaks around the interface between the atomistic and continuum regions, regardless of whether ghost forces exist in a coupling method. In the presence of ghost forces, the width of the peak can be estimated  \cite{MingYang:2009, ChenMing:2012, CuiMing:2013}. For example, in one space dimension, it is of the order $\veps|\ln\veps|$.  The width for a two-dimensional model is $\mathcal {O}(\veps^{\frac12}).$
Motivated by these observations, the current work aims to improve the accuracy at the interface. Similarly to the rep-atoms in the original QC framework, we select a number of coarse-grained (CG) variables to begin with. To motivate the method, we formulate the atomistic model as a two-stage minimization problem, from which we derive an effective model at the coarse scale.  We further draw the connection to the standard Galerkin projection. A simple calculation shows that in principle the {\it exact} effective model for a given set of CG variables can be derived; but compared to the Galerkin method, it contains an additional term in the stiffness matrix. The large error near the interface can thus be attributed to having neglected this term. In light of the serious issues with ghost
forces, we finally formulate the method as a force-based method, in which the ghost forces do not appear.

The purpose of this paper is to introduce proper approximations of the aforementioned missing term and demonstrate how the
approximation reduces the error near the interface. In particular, we propose to use the Krylov
subspace approximation and the Lanczos algorithm. With the connection to the Galerkin method, this can
be implemented by projecting the atomistic model to an extended subspace. The additional bases, referred to as {\it enriched bases}, are defined and supported near the interface, which do not require much extra computational cost. More importantly, unlike the basis functions in the standard Galerkin method, they take into account the atomic interactions and improve the
modeling accuracy.

Currently there is no analytical result that quantifies the numerical error in this approach.
 The current method does not show resemblance to any existing approach, and the analysis of the accuracy is
 still open.   But we have done extensive
numerical experiments to show how the enriched bases affect the accuracy. More specifically, we start
with a one-dimensional chain model, a test problem considered in many other works, and conduct extensive
numerical tests, including tests on the rate of convergence. In addition, we consider a one-dimensional model of fracture, which in spite of its simplicity, captures some essential aspects of crack propagation. Throughout the numerical study, our emphasis is focused on the error measured in
various norms, some of which have not been considered in previous theoretical analysis.
Issues on how to apply the current approach to two- or three-dimensional problems will be discussed in Section 2.4 and applications to more important micro-mechanical systems will be considered in subsequent works.

{


The remaining part of the paper is organized as follows. In Section \ref{sec:methodology}, we derive an effective
coarse-grained model for molecular mechanics, and introduce Galerkin methods with connection
to the QC method. Many numerical studies of Galerkin methods and other
multiscale methods are presented in Section \ref{sec:results}. Conclusions are drawn in Section
\ref{sec:conclusion}.

\section{Methodology}
\label{sec:methodology}


Consider a system with $N$ atoms interacting through the potential
function $U$, and let $\bmf_i^{\text{ex}}$ be the external force on the $i$-th
atom. At zero temperature, the total energy of the system can be
written as
\begin{equation}\label{eq: etot}
V(\by)=U(\by_1, \dots, \by_N)-\sum_{i=1}^N
\bmf_i^{\text{ex}}\cdot\by_i,
\end{equation}
where $\by_i$ is the position of the $i$-th atom in the deformed
state. An equilibrium atomic configuration is determined from the following
minimization problem:
\begin {equation*}
\{\by_1,\dots,\by_N\}=\argmin_{}V(\by)
\end {equation*}
with $\by$ subject to certain boundary condition. This is known as the molecular mechanics model.

Meanwhile, the displacement of the $i$-th atom is defined as
\[
\bu_i=\by_i-\bx_i,
\]
where $\bx_i$ is the reference position of the $i$-th atom, typically defined in an undeformed configuration. We rewrite the minimization problem in terms of the displacement as follows
\begin {equation}\label{eqn:mm}
\bu =\argmin_{}V(\bu),
\end {equation}
where $\bu= \{\bu_1,\dots,\bu_N\}$ and $V(\bu) = U(\bu_1,\dots,\bu_N) - \bu\cdot\bmf^{\text{ex}}$.
We drop the $-\bx\cdot\bmf^{\text{ex}}$ term which will not affect the solution of \eqref{eqn:mm}.
In what follows, a solution of \eqref{eqn:mm} will be considered as the {\it exact} solution for comparison.

\subsection{The motivation for the coarse-grained model}

Although $N$ is very large in real applications, only a small number of degrees of freedom
 is required to adequately describe the displacement field for systems with localized defects,
such as point defects, dislocations, and cracks. This has been the primary motivation for all
the existing atomistic to continuum coupling methods.  Conceptually, there are two
distinct regions, one of which contains local defects, where all the atoms are retained. The other is the surrounding bulk where the displacement is smooth. The first region will be referred to as {\it atomistic region}, and the latter will be called the {\it continuum region}, where a reduction of the model is possible.  For this purpose, let us first introduce general CG variables, denoted by $\bq$, without having to explicitly specify how the domains are divided.
Furthermore, let the number of CG variables be $n$ with $n\ll N$.
In practice, all the atoms near local defects are chosen, and much fewer atoms are needed in the surrounding region where the displacement is smooth. Here, we define the CG variables via  a linear operator $\Phi\in \R^{n\times N}$ which connects the fine scale variables $\bu$ to
the coarse variables $\bq$ according to the equation
\begin{equation}\label{eqn:cgv}
\bq= \Phi \bu.
\end{equation}
This is similar to the restriction operator in multigrid methods \cite{B77, McCormick87}.
We assume that $\Phi$ has full rank so that the CG variables are independently defined.  Meanwhile, $\Phi^T$ may be regarded as an interpolation operator such that $\bu \approx \Phi^T \bq$. This suggests that one can choose standard finite element basis functions to construct $\Phi.$ More details will be given later.

Problem \eqref{eqn:mm} can be formulated into two successive minimization problems \cite{EEngquist02}. To explain this, we first fix the CG variables $\bq$, and solve
\begin{equation}\label{eqn:min1}
 \min_{\bu: \Phi\bu=\bq} V(\bu)
\end{equation}
using \eqref{eqn:cgv} as constraints. Conceptually, to remove the constraints, we may let $X_0=\text{Range}(\Phi^T),$ and write $\bm u= \Phi^T (\Phi \Phi^T)^{-1}\bm q + \bm \xi$, with the free variable $\bm \xi \in X_0^\perp.$
Let the resulting minimal energy be $W_{\text{eff}}(\bm q)$.
In the second problem, we solve
\begin{equation}\label{eqn:min2}
 \min_{\bm q} W_{\text{eff}}({\bm q}).
\end{equation}
Obtaining a suitable approximation of the effective energy $W_{\text{eff}}(\bm q)$ is a general and fundamental problem in molecular modeling \cite{Leach01}.
The common goal is to derive an effective problem for $\bq$, from which efficient methods can be used to determine the CG variables. To motivate the proposed procedure,  we first consider a quadratic potential
\begin{equation}\label{eqn:harmonic}
V(\bu) = \half \bu^T A\bu - \bu^T\bmf^\text{ex}.
\end{equation}
Under this approximation, the exact solution to \eqref{eqn:min1} reads as
\begin{equation}\label{eqn:interpolation}
\bu= R \bq + Q_{\bu} A^{-1} \bmf^\text{ex}
\end{equation}
with
\begin{eqnarray*}
R&= &A^{-1}\Phi^T (\Phi A^{-1}\Phi^T)^{-1},\\
P_{\bu}&= &R\Phi,\\
Q_{\bu}&= &I - P_{\bu}.
\end{eqnarray*}
As a consequence, the effective model for $\bq$ is given by
\begin{equation*}
W_{\text{eff}}(\bq)= \half \bq^T R^T A R \bq - \bq^T R^T \bmf^\text{ex}.
\end{equation*}
If the force field is ``smooth" in the sense that it can be approximated by $\bmf^\text{ex}= \Phi^T \bg$ for some
$\bg \in \mathbb{R}^n$, then
the Euler-Lagrange equation of \eqref{eqn:min2} is given by
\begin{equation}\label{eqn:effective1}
(\Phi A^{-1}\Phi^T)^{-1} \bq= \bg.
\end{equation}
We remark that the reconstruction condition $\bmf^\text{ex}= \Phi^T \bg$ can be satisfied if a problem has localized
defects and the CG variables are chosen properly.

The above formulation can be further simplified by introducing $\Psi\in \R^{(N-n)\times N}$, which
maps the atomic information to the remaining degrees of freedom after coarse-graining.
$\Psi$ is chosen as an orthogonal matrix. Namely,
\begin{equation}
\Psi\Psi^T=I^{N-n}, \quad \Psi\Phi^T=0.
\end{equation}
Next we let $P= \Phi^T M^{-1}\Phi$ be the orthogonal projection to  $Y:=\text{Range}(\Phi^T)$
with $M=\Phi\Phi^T\in \R^{n\times n}$. Its complementary projection is $Q = \Psi^T \Psi.$
Using blockwise matrix inversion formulas, we can rewrite the matrix in \eqref{eqn:effective1} as follows
\begin{equation*}
(\Phi A^{-1}\Phi^T)^{-1} = M^{-1} \Phi A\Phi^T M^{-1} - M^{-1} \Phi A\Psi^T (\Psi A\Psi^T)^{-1} \Psi A\Phi^T M^{-1}.
\end{equation*}

It is now clear that the {\it exact} effective equation can be rewritten as
\begin{equation}\label{eqn:effective2}
(A_0 - A_1) \bp= \Phi \bmf^\text{ex},
\end{equation}
where
\begin{eqnarray}
\label{eqn:A0}
A_0&=&\Phi A\Phi^T,\\
\label{eqn:A1}
A_1&=&\Phi A\Psi^T (\Psi A\Psi^T)^{-1} \Psi A\Phi^T,\\
\label{eqn:p}
\bp &=& M^{-1}\bq.
\end{eqnarray}
These calculations, although quite simple, suggest an important approach to derive
accurate coupling methods. We further comment that
\begin{enumerate}
\item  The formulas above are exact for linearized models. For interested readers, a more detailed derivation can be found in the appendix of \cite{Li:2010}.
\item In equations \eqref{eqn:effective2}-\eqref{eqn:p}, $A_0$ and $M$ can be computed provided that
$A$ and $\Phi$ are given. However, it is difficult to compute $A_1$ since the calculation of the inverse
of $\Psi A\Psi^T$ is almost as difficult as solving the problem \eqref{eqn:harmonic}. As a result, we need to approximate $A_1$ using an alternative approach, instead of using \eqref{eqn:A1} directly. In the next section, we will provide a hierarchy of approximations. In particular, $A_0$ can be obtained using the standard Galerkin method, and $A_1$ corresponds to a Galerkin projection to  successively extended subspaces.
\item These calculations are limited to linearized models. For nonlinear models,
we will retain the Galerkin formalism, by projecting the full, nonlinear, atomistic model to the
extended subspaces.
\item In general, the external force $\bmf^{\text{ex}}$ may not be in the range of $\Phi^T$;
An example will be given later (see \eqref{eq:fnonlocal}). As a consequence, the displacement field cannot be well
approximated by basis functions in $\Phi$; see Fig. \ref{fig:Displacementnonlocal} for the
performance of piecewise linear basis functions. In this case, high-order basis functions
are more appropriate to reduce the approximation error and the enriched bases method can still
reduce the error around the interface.

\item It is also worth mentioning that the
 concepts of two-stage minimization and enriching procedure has been proposed in \cite{AbdulleLinShapeev2012,AbdulleLinShapeevII}. The enriching procedure is considered globally to achieve
a better approximation due to the internal degrees of freedom. However, the enriching procedure
in Section 2.3 is only considered locally to reduce the interfacial error with a compromise between accuracy and efficiency.
\end{enumerate}


\subsection{The standard Galerkin method} If $A_1$ is neglected in \eqref{eqn:effective2}, we have $A_0\bp = \Phi\bmf$
with $A_0=\Phi A\Phi^T$. Therefore, it is equivalent to the standard Galerkin method, which has been a useful tool for approximating solutions of PDEs \cite{LaTh03}. To explain this more precisely, we let $X_0= \text{Range}(\Phi^T).$ We seek $\bp_0 \in X_0$, such that,
\begin{equation}\label{eqn:Galerkin}
( A \bp_0, \varphi ) = (\bmf^\text{ex}, \varphi), \quad \forall \; \varphi \in X_0.
\end{equation}
Here the parenthesis stands for the standard inner product. Namely, $(\bmf^\text{ex}, \varphi)= \sum_i \bmf^\text{ex}(x_i) \varphi(x_i).$

With direct calculations, one can verify that this formulation leads to the stiffness matrix $A_0$, and
the effective equation \eqref{eqn:effective2} (without $A_1$). This type of weak forms  have been used extensively in the adaptive QC methods \cite{OdPrRoBa06,prudhomme2006error,prudhomme2009adaptive}.


As a simple illustration, we consider a test problem, where the atomistic system consists of a chain of atoms with up to the second nearest neighbor interaction. The setup of the problem is as follows: We include displacements of all the atoms on the right as the CG variables, while on the left, we choose one atom out of every few atoms. The operator $\Phi$ is defined based on the standard piecewise linear nodal bases.
The details will be explained in the next section. From  Fig. \ref{fig:DisplacementErrorStandard2}, we observe that the error is concentrated near
the interface. More examples, \eg, Figs.  \ref{fig:DisplacementErrorStandardInhomogeneous} and \ref{fig:DisplacementErrorCrack1},
further support this observation. Therefore, when the error is measured under a uniform norm, \eg, $W^{1,\infty}$, the error at the interface would dominate.

We attribute this phenomenon to the absence of $A_1$ in the effective model \eqref{eqn:Galerkin}.
In Fig. \ref{fig:DisplacementErrorStandard3}, we plot the $l^2$ norm of each row of $A_1$ as a function of the
atomic index. It is clear that  the main contribution of $A_1$ is localized around the interface.
In the following section, we describe how to incorporate $A_1$ into the model, and through numerous examples, we demonstrate the reduction of the error.

%
%

\subsection{Enriched bases method}
We now try to improve the approximation by introducing an extended approximation. In the case
when the potential is quadratic, this amounts to computing the matrix $A_1$.
Since $A_1$ is difficult to compute directly,  we will approximate it using a Krylov subspace method.
In particular, the expression \eqref{eqn:A1} involves the solution of a linear system with multiple right-hand sides, for which the Krylov subspace method is especially useful \cite{el2004block}.
More specifically, let $V= \Psi A\Phi^T$, and we define
\begin{equation}\label{eqn:Krylov1}
K_\ell\big(\Psi A\Psi^T,V\big)= \text{span}\Big\{V, \Psi A\Psi^T V, \cdots, (\Psi A\Psi^T)^\ell V \Big\}.
\end{equation}
For the general theory regarding Krylov subspaces, we refer the reader to \cite{Saad:2011, saad1980rates}. Often of practical interest are the orthogonal bases for $K_\ell\big(\Psi A\Psi^T,V\big)$, which will be denoted by $V_\ell$, and we have
\begin{equation*}
(\Psi A\Psi^T) V_\ell \approx V_\ell T,
\end{equation*}
where $T$ is a block tri-diagonal matrix with dimension $\le (\ell+1)n$. The error is a matrix with low rank.  Thanks to this approximation, we have
\begin{equation*}
(\Psi A\Psi^T)^{-1} V_\ell \approx V_\ell T^{-1}.
\end{equation*}
As a result,
\begin{equation*}
A_1 \approx E_1^T T^{-1} E_1,
\end{equation*}
where $E_1^T=[I, 0, \cdots, 0].$

To compute the basis vectors, a Block Lanczos algorithm \cite{Saad:2011} can be used. The detailed algorithm is shown below.
\begin {algorithm}
\label{alg:BlockLanczos}
\caption{Block Lanczos algorithm} \label{alg:BlockLanczos}
\begin{algorithmic}
\State Set $V_0=0$ and $Z_0=V$.
\ForAll{$j=1, 2, \cdots, \ell$,}
\State Rank revealing QR factorization of the $n\times p_{j-1}$ matrix $Z_{j-1}$: $Z_{j-1} = Q_j R_{j-1}$. $R_{j-1}$ may be a permuted upper triangular matrix;
\State Let $p_j=\text{rank}(Z_{j-1}),$ $V_j$ be the first $p_j$ columns of $Q_j$,
and $B_{j-1}$ be the first $p_j$ rows of $R_{j-1}$;
\State $Z_j \longleftarrow A V_j - V_{j-1}B_{j-1}^T$;
\State $A_j \longleftarrow V_j^T Z_j$;
\State $Z_j \longleftarrow Z_j - V_j A_j$.
\EndFor
\end{algorithmic}
\end{algorithm}

From a practical viewpoint, the main difficulty associated with implementing {\bf Algorithm} \ref{alg:BlockLanczos} is due to the fact that the
matrix $\Psi$ is difficult to access. Furthermore, the algorithm seems to be limited to linear problems. Here {\it three} steps are taken to alleviate these difficulties.

\smallskip

First, we define a Krylov subspace which is equivalent to \eqref{eqn:Krylov1}; but it can be implemented much more efficiently. In particular, we consider the Krylov space
\begin{equation}\label{eqn:Krylov2}
K_\ell\big(QA,W\big)= \text{span}\Big\{W, QA W, \cdots, (QA)^\ell W \Big\},
\end{equation}
with $W=Q A \Phi^T $. It is easy to check that {\bf Algorithm} \ref{alg:BlockLanczos}, when applied to $K_\ell\big(QA,W\big)$, yields the same result as the method applied to $\Psi^T K_\ell\big(\Psi A\Psi^T,V\big)$. The advantage of \eqref{eqn:Krylov2} is that $Q$ can be obtained as the complementary projection of $P\in \R^{n\times n}$, which can be computed from $\Phi.$

\smallskip

Secondly, we recall the previous observation that  the error introduced in the standard Galerkin method tends to center around the interface between the continuum region and the atomistic region. Therefore we only choose {\it a  subset of }basis functions in $\Phi$, which correspond to the basis functions around the interface. To make it more specific, we let $\wt{\Phi}$ be the matrix that contains $m$ additional basis functions with $m \ll n$, and we define $\wt{W}= QA\wt{\Phi}^T$, and accordingly, $K_\ell\big(QA,\wt{W}\big)$. Since the atomic interaction is assumed to be of finite range, these additional basis functions still have compact supports, which makes the Lanczos algorithm feasible in the current setting.

Once the Krylov subspace is generated, we define an extended space by
\begin{equation*}
X = X_0 \oplus K_\ell\big(QA,\wt{W}\big).
\end{equation*}
The additional basis functions depend on the force constant matrix $A$. Unlike the
basis functions in $\Phi$, they incorporate the atomic interactions. Motivated by this
observation, these basis functions will be referred to as {\it enriched bases.}

\smallskip

At the third step, we return to the fully nonlinear atomistic model \eqref{eqn:mm}.
We will project the full (and nonlinear)
 model to the subspace $X$ using a Galerkin projection: We seek $\bu\in X$, such that
\begin{equation}\label{eqn:extendedGalerkin}
( -\nabla V(\bm u), \varphi ) = (\bmf^\text{ex}, \varphi),\quad   \forall\;\varphi \in X.
\end{equation}
This weak form was also used in the adaptive QC methods \cite{OdPrRoBa06,prudhomme2006error,prudhomme2009adaptive}.
Clearly, this Galerkin formulation is based on the force balance, rather than a minimization of the energy, which
might suffer from the issues of the ghost forces.

The above three steps constitute the main components of the enriched bases method.
In the next section, we discuss how to approximate the inner products in the weak forms.

\subsection{Quadrature approximations and connections to the QC method}
To reduce the computational cost in \eqref{eqn:Galerkin} and \eqref{eqn:extendedGalerkin}, we
introduce quadrature approximations in the continuum region where the displacement is smooth. Notice that the enriched bases are only confined to a neighborhood of the interface, so the standard Galerkin method and the enriched bases methods have the same set of basis functions in the interior of the continuum region, and our discussion on the quadrature approximation applies to both methods.

 There are three terms that need to be approximated using the quadrature rule in
\eqref{eqn:Galerkin}: $\Phi\bmf^\text{ex}$, $\Phi\Phi^T$, and $\Phi A\Phi^T$. As an example, we consider
the case where the continuum region has been equipped with elements: Intervals (1D), triangles (2D), and
tetrahedrons (3D), with basis $\varphi_j$ being a standard piecewise linear function which forms
the $j$-th row of $\Phi$.

Let $i$ and $j$ be two vertices of an interval $\T$ in 1D (or a triangle
in 2D and a tetrahedron in 3D). The $(i,j)$ component of $M=\Phi\Phi^T$ can be computed as
\begin{equation*}
M_{i,j} = \sum_{k\in\T} \varphi_i(\bx_k)\varphi_j(\bx_k) \approx \frac{1}{\V_0}\int_{\T} \varphi_i(\bx)\varphi_j(\bx)\text{d}\bx,
\end{equation*}
where $\V_0$ is the volume of the unit cell. This integral can be further approximated by standard quadrature
formulas.

We can  approximate the components of $F=\Phi\bmf^\text{ex}$ and $K=\Phi A\Phi^T$ by
\begin{eqnarray*}
F_i = \sum_{k\in\Omega_i}\varphi_i(\bx_k)\bmf(\bx_k) & \approx & \int_{\Omega_i} P(\bx)\nabla \varphi_i(\bx)\text{d}\bx, \\
K_{i,j} = \sum_{k\in\Omega_i}\sum_{m\in\Omega_j}\varphi_i(\bx_k)a_{k,m}\varphi_j(\bx_m) & \approx & \int_{\Omega_i} C(\bx):\big(\nabla \varphi_i(\bx)\otimes\nabla\varphi_j(\bx)\big)\text{d}\bx,
\end{eqnarray*}
where $\Omega_i$ and $\Omega_j$ are supports of $\varphi_i$ and $\varphi_j$, respectively.
Here $P$ and $C$ can be identified as the first Piola-Kirchhoff stress tensor and the Lagrangian tangent stiffness tensor,
respectively. The detailed construction of $P$ and $C$ from the atomistic model can be found in \cite{TadmorOrtizPhillips:1996}. Both integrals can be further approximated by standard quadrature formulas.
Interested readers may find earlier discussions on quadrature approximations
in \cite{KnapOrtiz:2001, Lin07, LuskinOrtner:2009, Gunzburger:2010} for static problems
and \cite{Li:2014,WaLi03} for dynamics problems.

In the standard Galerkin method, if we use piecewise linear functions as basis functions and the mid-point rule
for the quadrature approximation for all elements in \eqref{eqn:Galerkin}, we get the local QC method \cite{TadmorOrtizPhillips:1996, EMing:2005}.


Once the quadrature approximation is used in Galerkin methods, three regions arise: The atomistic region where the force is treated ``as is", one continuum region with the quadrature approximation, which is equivalent to a finite element approximation of an elasticity model, and one region in between with no quadrature approximation. This region at the interface will be referred to as the {\it interbedded region}, which is considered as an addition to the original QC method to specifically reduce the coupling error. This is illustrated in Fig. \ref{fig:qa2d}.
 Kinematically, the number of degrees of freedom in the interbedded region is reduced by introducing rep-atoms (coarse grained elements). However, energies or forces over these elements are still calculated based on the summation of site energies or site forces per atom. 
 Therefore, atoms in the interbedded region act as {\it quasiatoms},
which play a similar transitional role to the quasi-nonlocal atoms \cite{Shimokawa:2004}.
\begin{figure}[htp]

\centering
\vspace{-0em}
\rotatebox{90}{\includegraphics[width=3in]{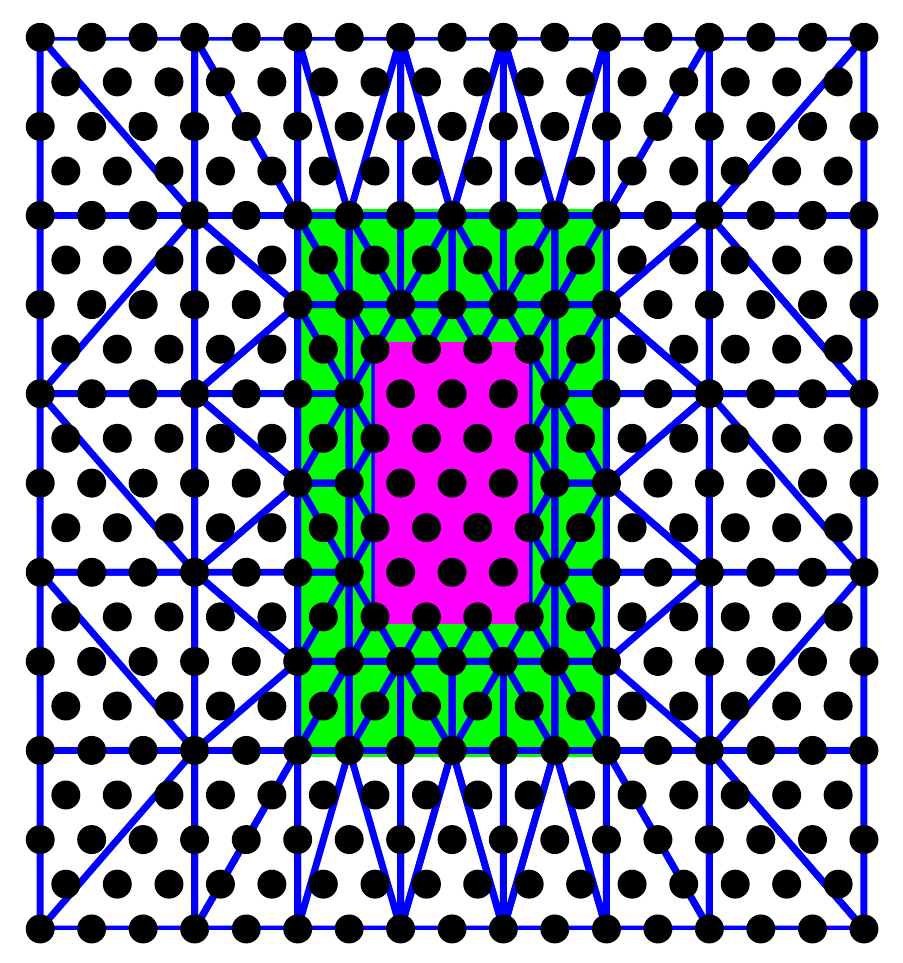}}
\vspace{-0em}
\caption{\small A schematic picture of the interbedded region and quasiatoms. The
interbedded region with quasiatoms separates the atomistic region from the continuum region,
serving as a smooth model transition.}\label {fig:qa2d}
\end{figure}

Clearly, formulas presented here and the Galerkin projection are valid in high dimensions as well.
In Fig.  \ref{fig:qa2d} we have illustrated how this can be implemented in 2D. More specifically, the enriched bases will be
introduced in the elements next to the interface (green), and in those elements, the forces are evaluated exactly. On the other
hand, in elements further away from the interface (white), the quadrature formulas shown above are used to compute
the average forces. This part is equivalent to the local QC method.  In fact, this procedure for approximating the forces has been implemented in \cite{Li:2014,WaLi03}, although in different settings.

The original QC method \cite{MingYang:2009, ChenMing:2012, CuiMing:2013} and
other multiscale methods without ghost forces \cite{LuskinOrtner:2013} share the similarity that the error
in a multiscale method is dominated around the interface. In our approach, on the one hand, the error in the
interbedded region decays quickly away from the interface between the interbedded region and the atomistic region. At the interface between the continuum region and the interbedded region, the quadrature approximation will make
an error, which is often much smaller compared with the error at the interface between the interbedded region and
the atomistic region. On the other hand, there is no energetic reduction for quasiatoms in the interbedded
region. The enriched bases can be naturally added in the framework of a Galerkin method, which can reduce
the error between the interbedded region and the atomistic region significantly.
Moreover, the quadrature approximation and the enrichment are combined in a concurrent way. One helps the
other, and vice versa. This leads to a much better approximation to the atomistic model without increasing
the computational complexity. All these will be verified by simulation results in Section \ref{sec:results}.

\section{Numerical Results}\label{sec:results}

\subsection{A one-dimensional chain model}

Consider a one-dimensional chain model over the unit interval $[0, 1]$. For comparison with some existing  multiscale methods, only up to the second-neighbor interaction is included.
The total energy is written as
\begin{equation}
 V= \sum_j \Big[ U \left (\frac{y_{j+2}-y_j}{\veps}\right ) + U\left (\frac{y_{j+1}-y_j}{\veps}\right ) - f_j^\text{ex} y_j\Big ],
\end{equation}
where $U$ is the interacting potential and $\veps$ is the lattice spacing constant.
Therefore, the equilibrium equation of the
$j$-th atom reads as
\begin{equation*}
-U'\left (\frac{y_{j+2}-y_j}{\veps} \right ) - U'\left (\frac{y_{j+1}-y_j}{\veps}\right ) + U'\left (\frac{y_j-y_{j-1}}{\veps}\right ) + U'\left (\frac{y_j-y_{j-2}}{\veps}\right ) = \veps f_j^\text{ex},\,
\end{equation*}
for $1\le j \le N.$
This simple model has been considered as a test problem in many coupling methods.

As a starting point, the linearized version of the above equation will be tested
to illustrate the idea of Galerkin methods.  Under this assumption, we have
\begin{equation*}
(2K_0+2K_1)y_j - K_0(y_{j+1}+y_{j-1}) - K_1(y_{j+2}+y_{j-2}) = \veps^2f_j^{\text{ex}},
\end{equation*}
where $K_0$ and $K_1$ are spring constants for the first and the second neighbor interactions, respectively.
Both the standard Galerkin method and the enriched bases method have been tested for $K_0$ and $K_1$ over a number of choices with similar performances, including parameters from linearized Lennard-Jones potential and linearized Morse potential. In what follows, we choose $K_0 =4$ and $K_1 = 1.4$.

\subsubsection{A Direct Comparison of the Error}
In the first set of tests, we consider a system with $1024$ atoms, where the left $512$ atoms are coarse-grained and the remaining atoms are kept. This divides the system into two sub-regions. In the coarse-grained (continuum) region, we choose one CG variable out of every eight atoms, and we use standard piecewise linear basis functions. Namely, for the nodal basis $\varphi_i(x)$, we define $q_i=\sum_j \varphi_i(x_j) u(x_j).$  In the atomistic region, every atom is chosen as a CG variable. This implies that if the $i$th atom is in the atomistic region, we include the trivial basis function $\varphi_i(x_j)=\delta_{i,j}$ in the mapping $\Phi.$ To illustrate this better, we show in Fig. \ref{fig:Basis1d} the basis functions near the interface for cases with a uniform mesh and a nonuniform mesh. The nonuniform mesh is constructed by a gradual increase of mesh size from the atomistic region to the
continuum region. Also shown in the Figure are the  extended basis functions around the interface.
\begin{figure}[htbp]

\vspace{-4em}
\centering
\subfigcapskip -3em
\subfigure[Uniform]{\label{fig:Basis1duniform}
\includegraphics[width=1.6in]{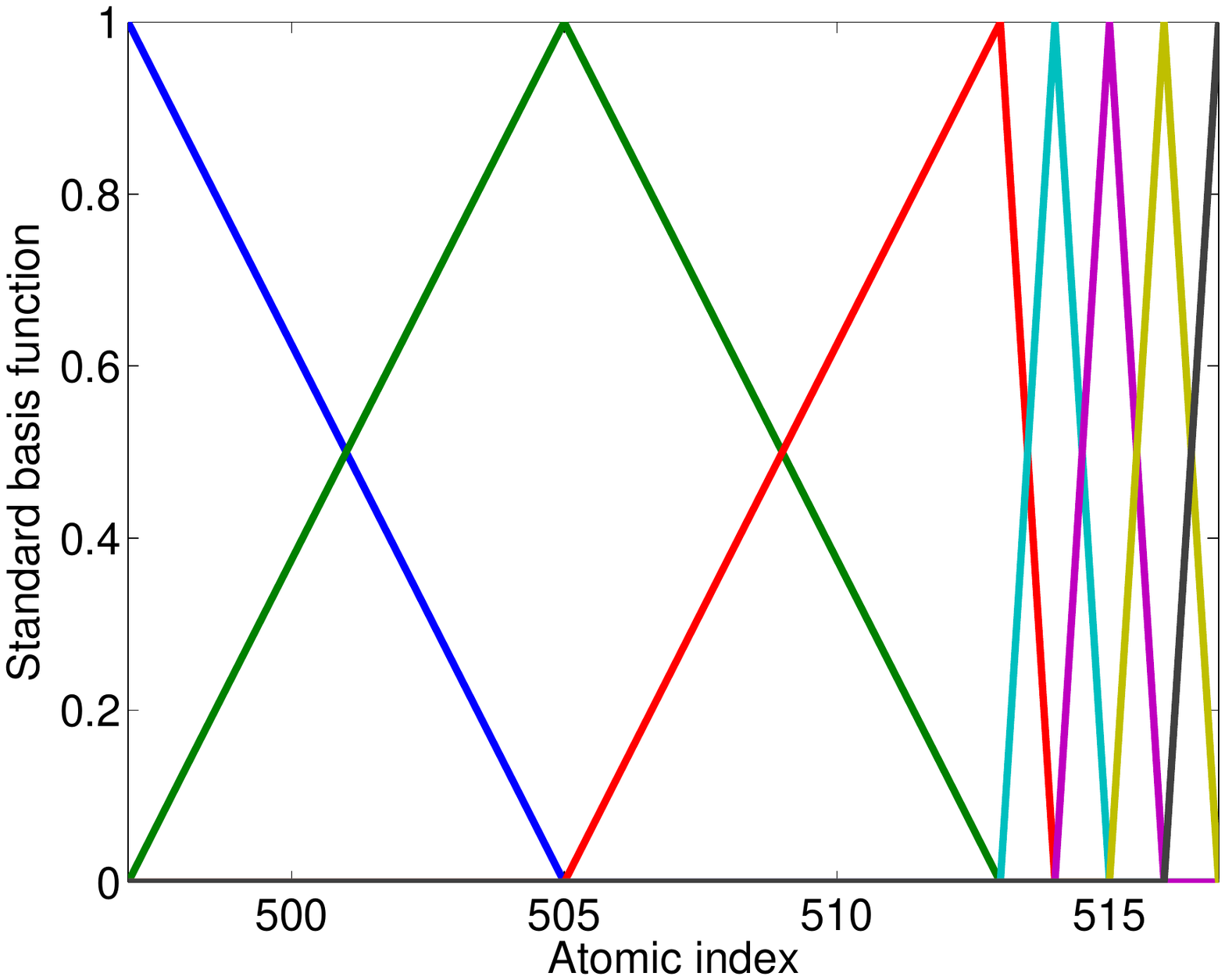}}%
\subfigure[Nonuniform]{\label{fig:Basis1dnonuniform}
\includegraphics[width=1.6in]{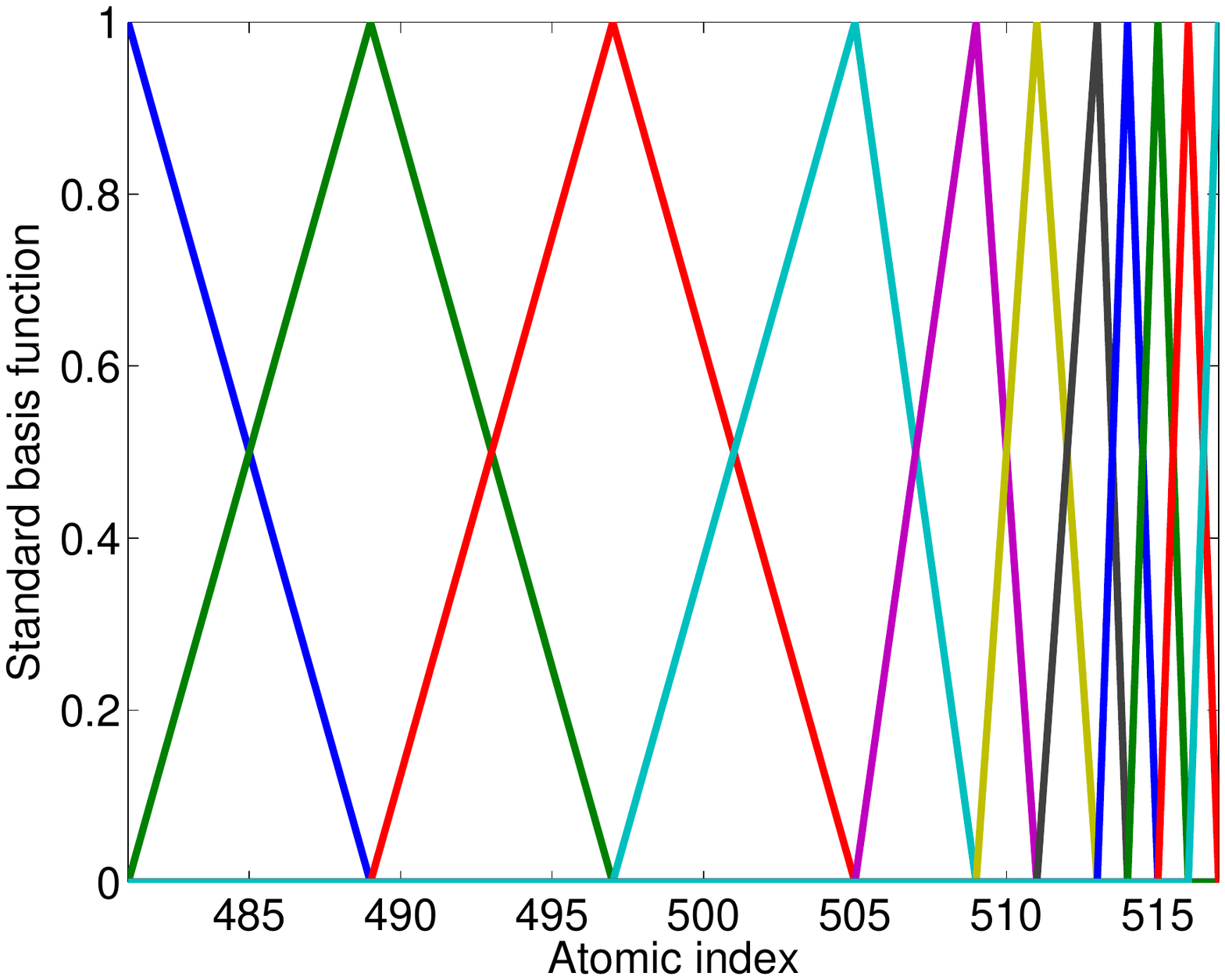}}%
\subfigure[Enriched bases]{\label{fig:EnrichedBasis}
\includegraphics[width=1.6in]{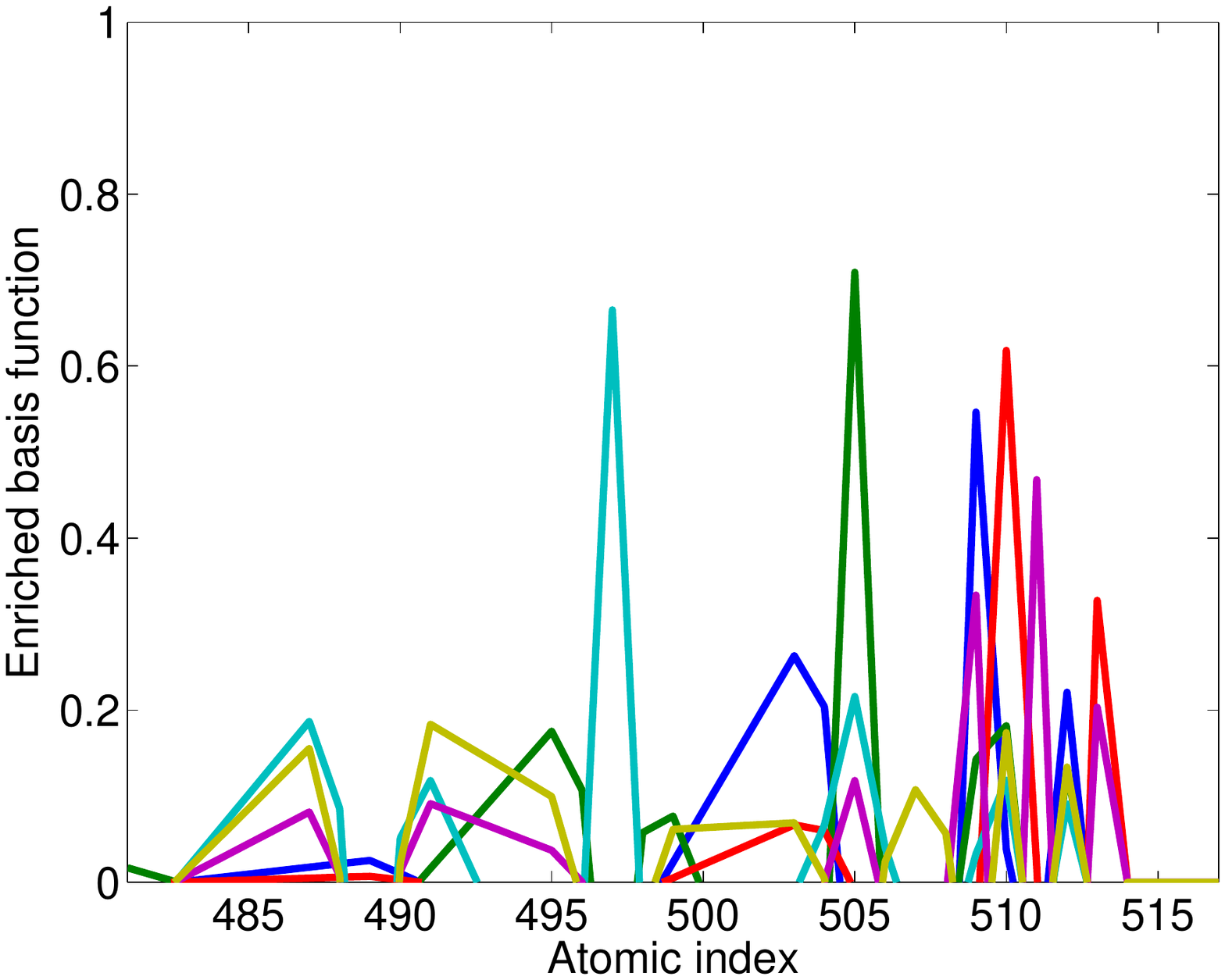}}
\vspace{-3em}
\caption{\small Basis functions around the interface in Galerkin methods with a uniform mesh, with a nonuniform mesh, and the enriched bases method with $\ell=1$, respectively. (a): Uniform mesh; (b): Nonuniform mesh;
(c): Enriched basis functions in the enriched bases method with $\ell=1$. The horizontal axis indicates the labels of the atoms near the interface. }\label{fig:Basis1d}
\end{figure}

To mimic the effect of local defects, a point force is applied on the $514$-th atom, which is close to the interface between the
atomistic region and the continuum region. Dirichlet boundary conditions are used for both the
left and the right boundaries. Namely, $u(0)=u(1)=0,$ or $u_1=u_N=0$. Fig. \ref{fig:DisplacementDelta} shows the displacement field and
the displacement gradient of the atomistic model. This solution will serve as the exact solution for the comparison purpose.

\begin{figure}[htbp]

\vspace{-6em}
\centering
\subfigcapskip -5em
\subfigure[Displacement]{\label {fig:Displacement1dDelta}
\includegraphics[width=2.5in]{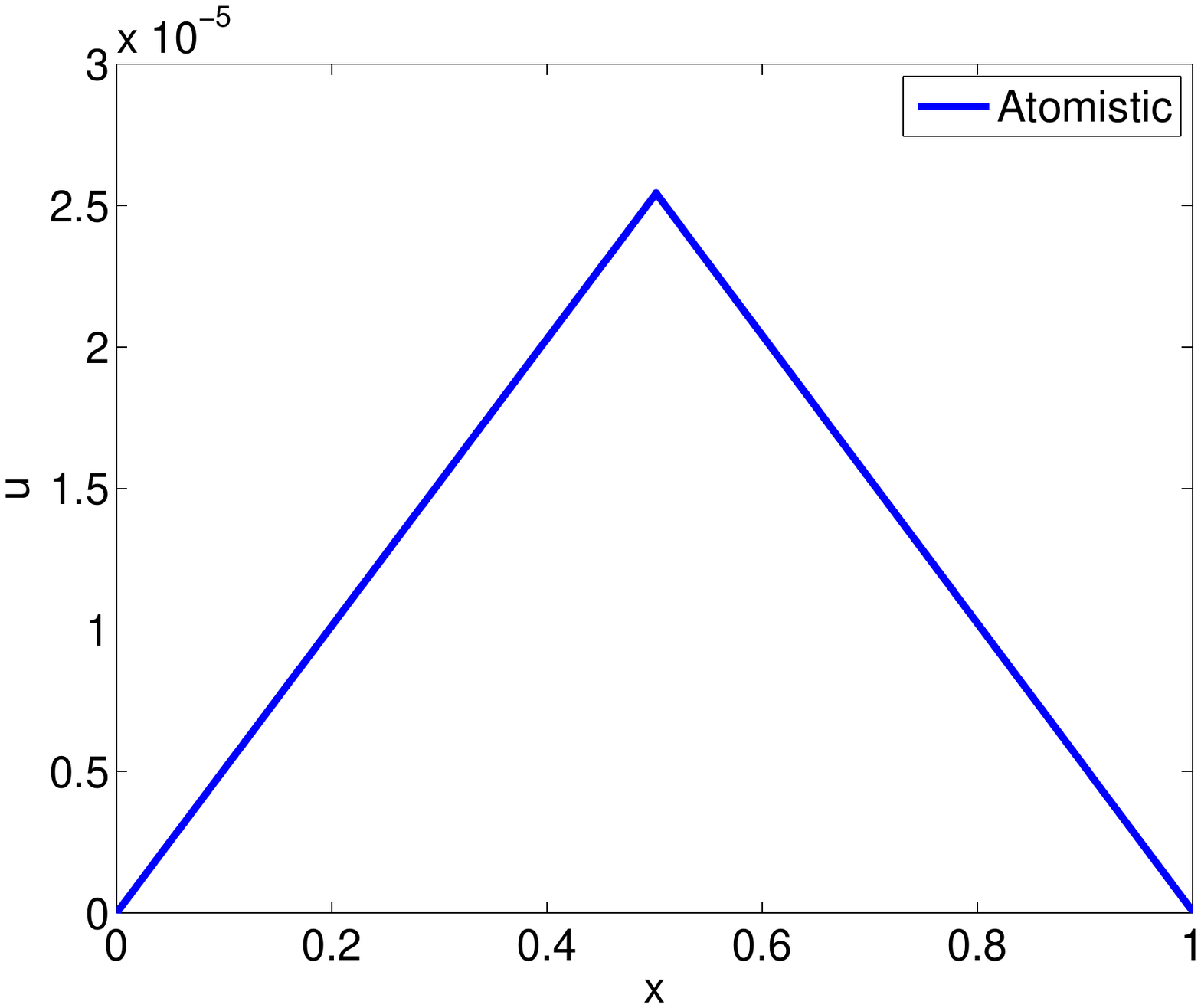}}%
\subfigure[Displacement gradient]{\label {fig:DisplacementGradient1dDelta}
\includegraphics[width=2.5in]{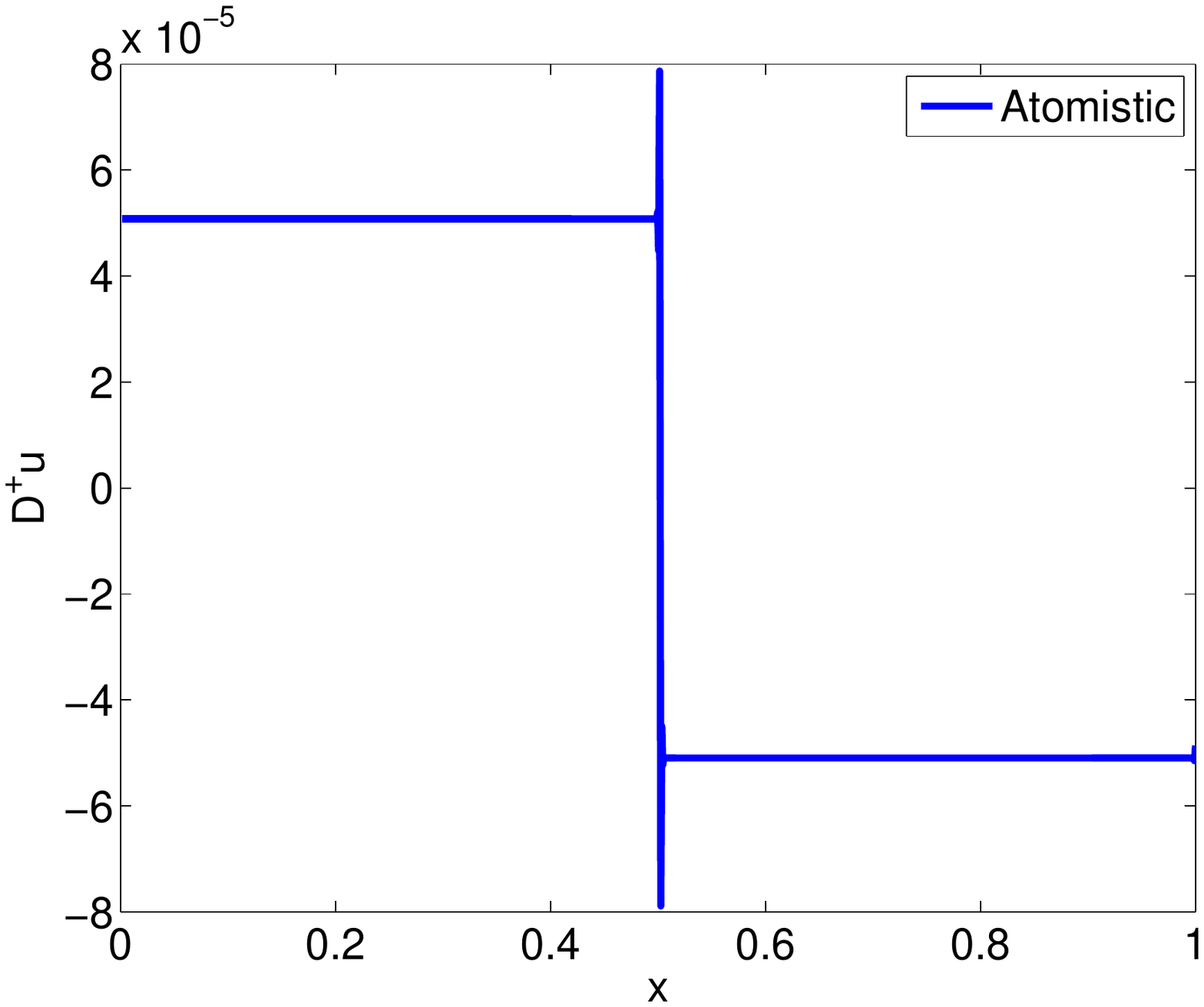}}%
\vspace{-4em}
\caption{\small Displacement and displacement gradient of a system with $1024$ atoms given by the atomistic model
with a point force applied on the $514-$th atom.
(a): Displacement; (b): Displacement gradient. This external force satisfies the requirement that $\bmf^{\text{ex}}\in \text{Range}\big(\Phi^T\big)$ since the point is inside the atomistic region.}\label{fig:DisplacementDelta}
\end{figure}

We first compare the standard Galerkin method to the atomistic model, for which the displacement error is shown in Fig. \ref{fig:DisplacementErrorStandard1}. The error clearly develops a peak at the interface ($x=\frac12$).
In addition, large error is also observed at the boundaries. This is due to the fact that the interaction is among first and second nearest neighbors, in which case the boundary condition is quite subtle.
A boundary layer may appear when one simply keeps the first and the last atom fixed.
To resolve this issue, we use an extrapolation technique. More precisely, we use displacements of
interior atoms to extrapolate displacements of ghost atoms (additional atoms introduced outside the boundary), and then substitute them into equilibrium equations of atoms near the boundary. For the present model with up to the second neighbor interaction, we first compute the displacement gradient of atom $1$ using one sided interpolation
\[
\frac{-3y_1+4y_2-y_3}{2\veps}
\]
The displacement of the ghost atom $0$ is then determined from
\begin{align*}
\frac{y_2-y_0}{2\veps} & = \frac{-3y_1+4y_2-y_3}{2\veps}.
\end{align*}
We solve for $y_0$ in the above equation and substitute it into the equilibrium equation of atom $2$. The error associated with the solution obtained using this boundary condition is shown in Figs. \ref{fig:DisplacementErrorStandard2}, from which one can clearly see the elimination of the boundary effect. Similar treatment will be used in subsequent tests, \eg, in Fig. \ref{fig:DisplacementErrorExtended2}.
\begin{figure}[htbp]

\vspace{-4em}
\centering
\subfigcapskip -3em
\subfigure[Without extrapolation]{\label{fig:DisplacementErrorStandard1}
\includegraphics[width=1.6in]{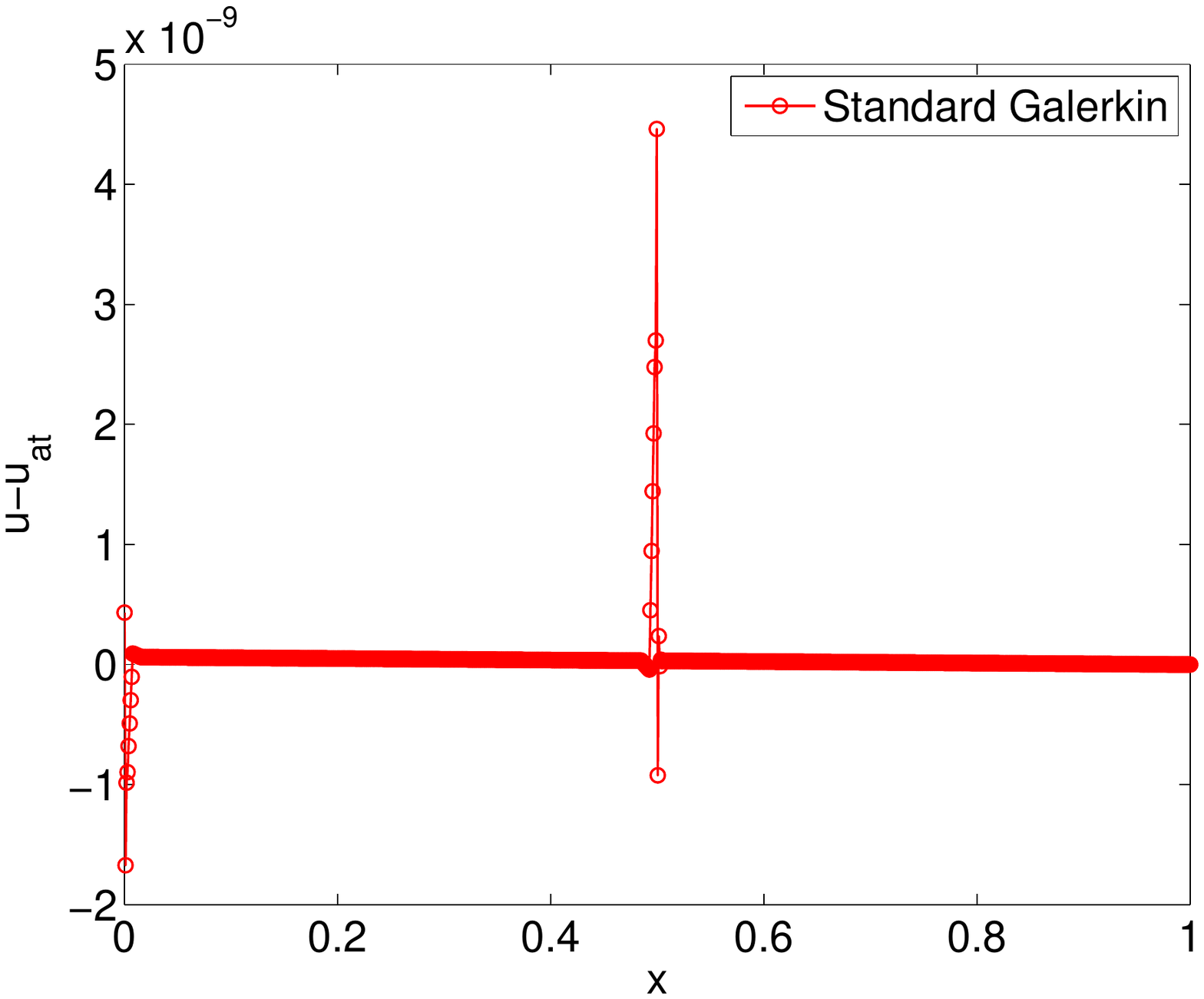}}%
\subfigure[With extrapolation]{\label{fig:DisplacementErrorStandard2}
\includegraphics[width=1.6in]{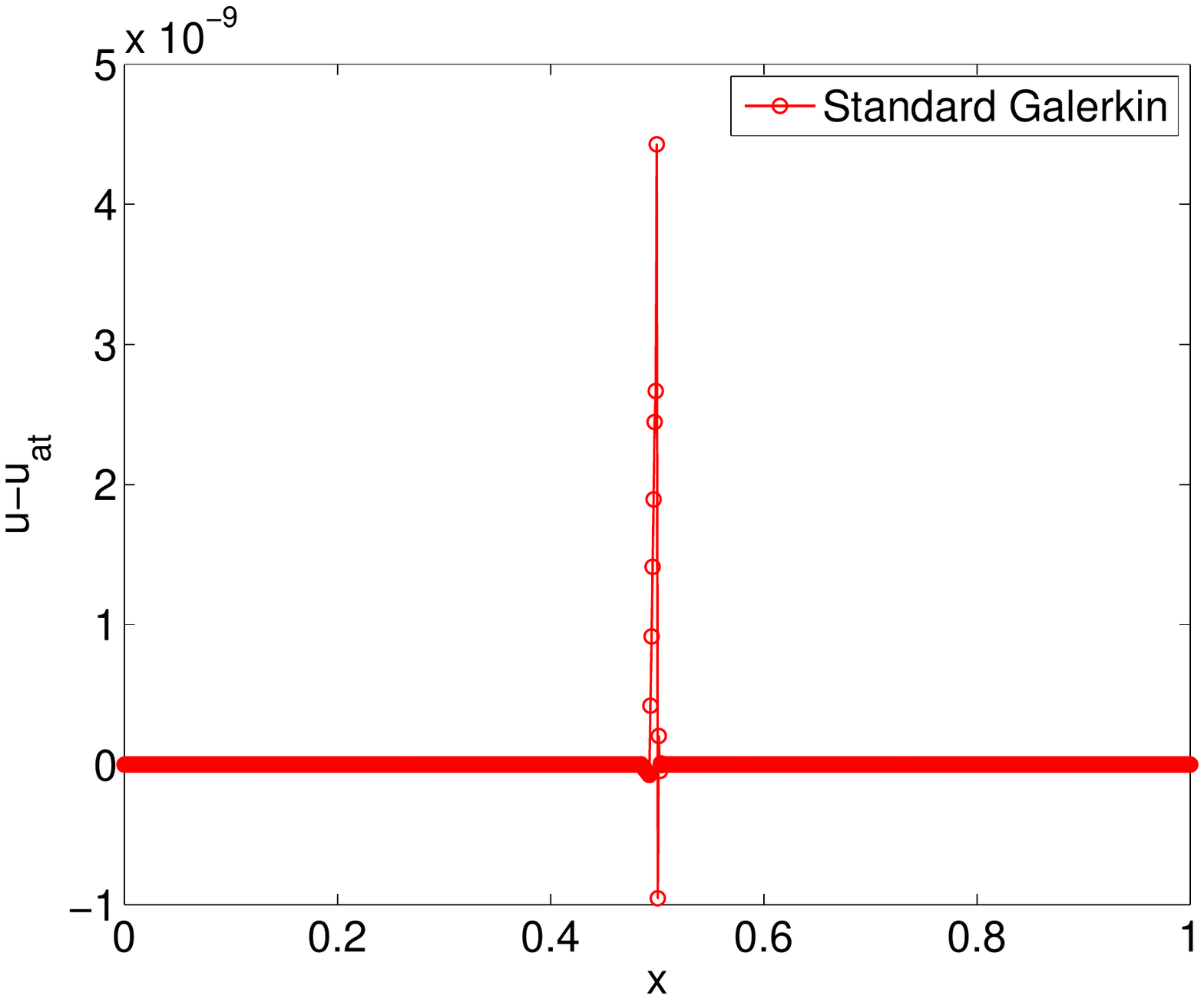}}%
\subfigure[$l^2$ norm of each row of $A_1$]{\label{fig:DisplacementErrorStandard3}
\includegraphics[width=1.6in]{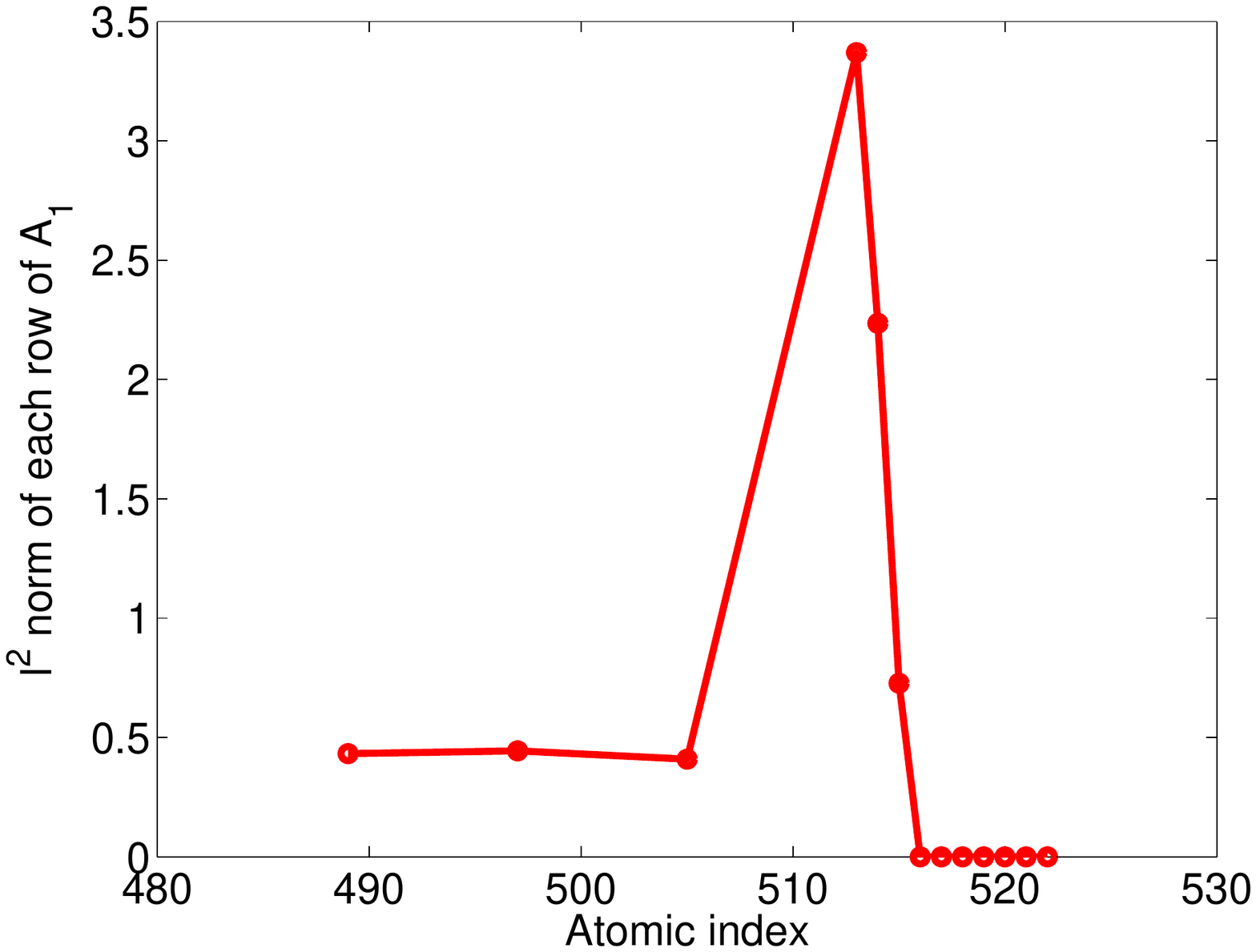}}
\vspace{-3em}
\caption{\small Displacement errors of the standard Galerkin method.
(a): Without extrapolation for atoms near the left boundary;
(b): With extrapolation for atoms near the left boundary;
(c): $l^2$ norm of each row of $A_1$ as a function of atomic index.
The error is localized around the interface in the standard Galerkin method, where
$A_1$ cannot be neglected.}\label{fig:DisplacementErrorStandard}
\end{figure}

For the enriched bases method, six basis functions
around the interface are selected to construct enriched bases from the Krylov subspace.
Recall that $\ell$ denotes the index of Krylov subspaces \eqref{eqn:Krylov2}.
Figs.
\ref{fig:DisplacementErrorExtended1} and
\ref{fig:DisplacementErrorExtended2}  demonstrate the displacement error for the solutions obtained
from the  enriched bases method, with and without the extrapolation at the boundaries. In particular,
from Fig.  \ref{fig:DisplacementErrorExtended2}, we find that with the extrapolation, the boundary effect
is eliminated. Furthermore, the error, which still concentrates at the interface; But it is of much smaller magnitude when compared to the result of the standard Galerkin method.

\begin{figure}[htbp]

\vspace{-6em}
\centering
\subfigcapskip -5em
\subfigure[Without extrapolation]{\label {fig:DisplacementErrorExtended1}
\includegraphics[width=2.5in]{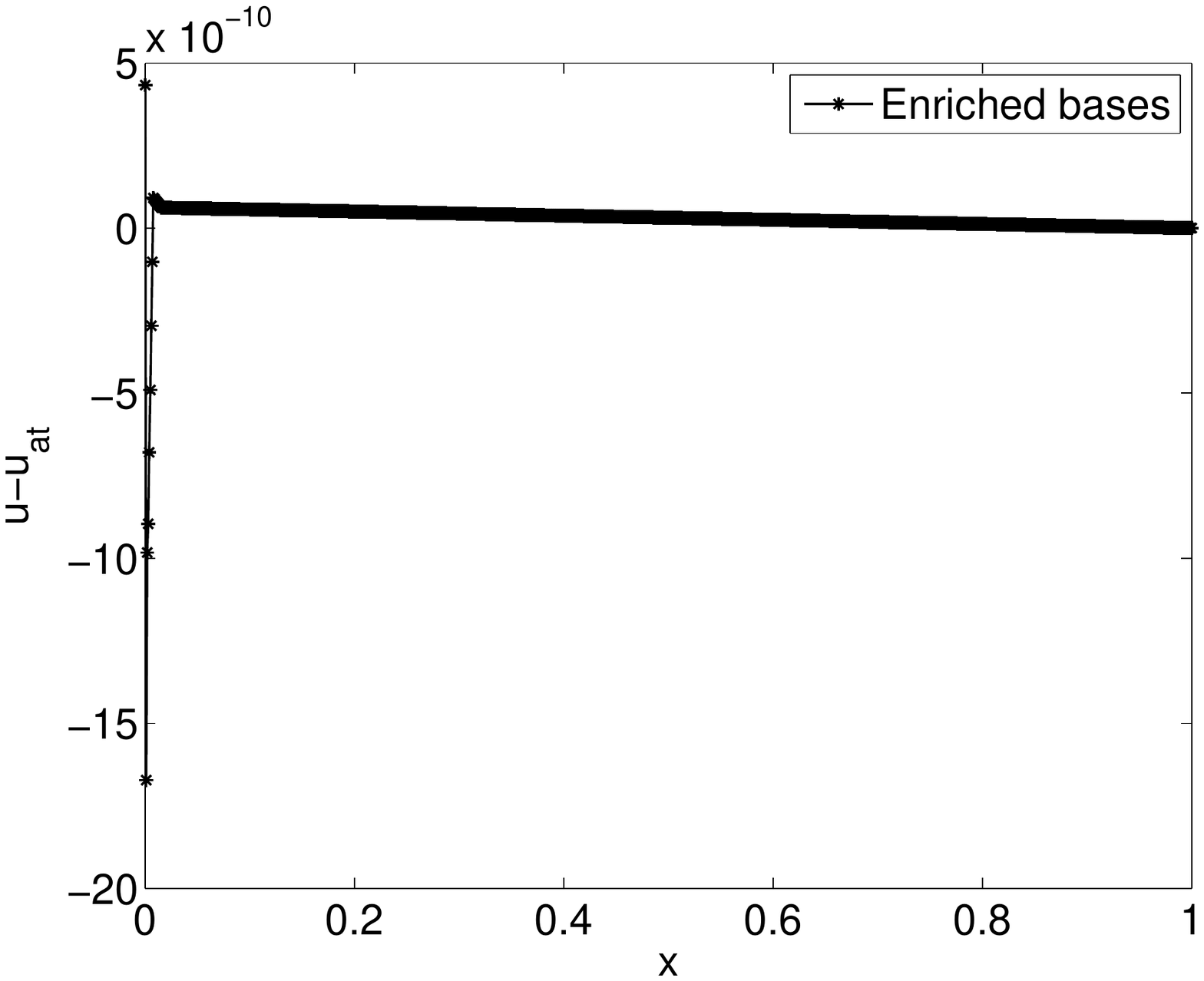}}%
\subfigure[With extrapolation]{\label {fig:DisplacementErrorExtended2}
\includegraphics[width=2.5in]{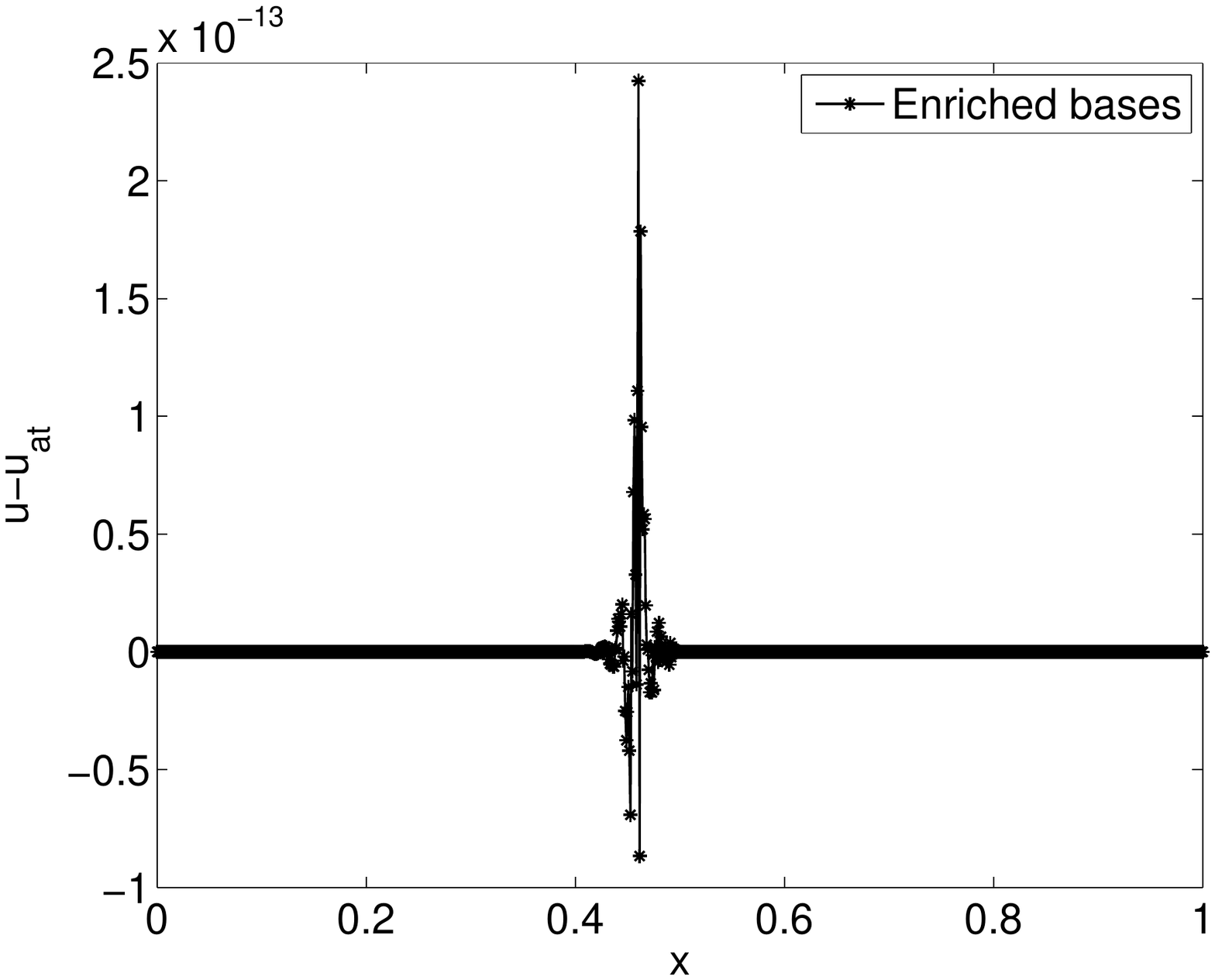}}%
\vspace{-4em}
\caption{\small Displacement errors of the enriched bases method
 with $\ell=5$.
(a): Without extrapolation for atoms near the left boundary;
(b): With extrapolation for atoms near the left boundary.}\label {fig:DisplacementErrorExtended}
\end{figure}

\smallskip

We further look at the accuracy of the enriched bases method with the order of the Krylov subspaces changed from $1$ to $5$.
The results are collected in Fig. \ref{fig:DisplacementErrorExtendedVaryL}. Dramatic reduction of the error is observed when $\ell$ changes from $1$ to $3$, and the error begins to saturate gradually when $\ell$ further increases. In order to reduce the error further, one has to include more basis functions in $\wt{\Phi}$ near the interface.

\begin{figure}[htbp]

\vspace{-6em}
\centering
\subfigcapskip -5em
\subfigure[$\ell=1$]{\label {fig:DisplacementErrorExtendedVaryL1}
\includegraphics[width=2.5in]{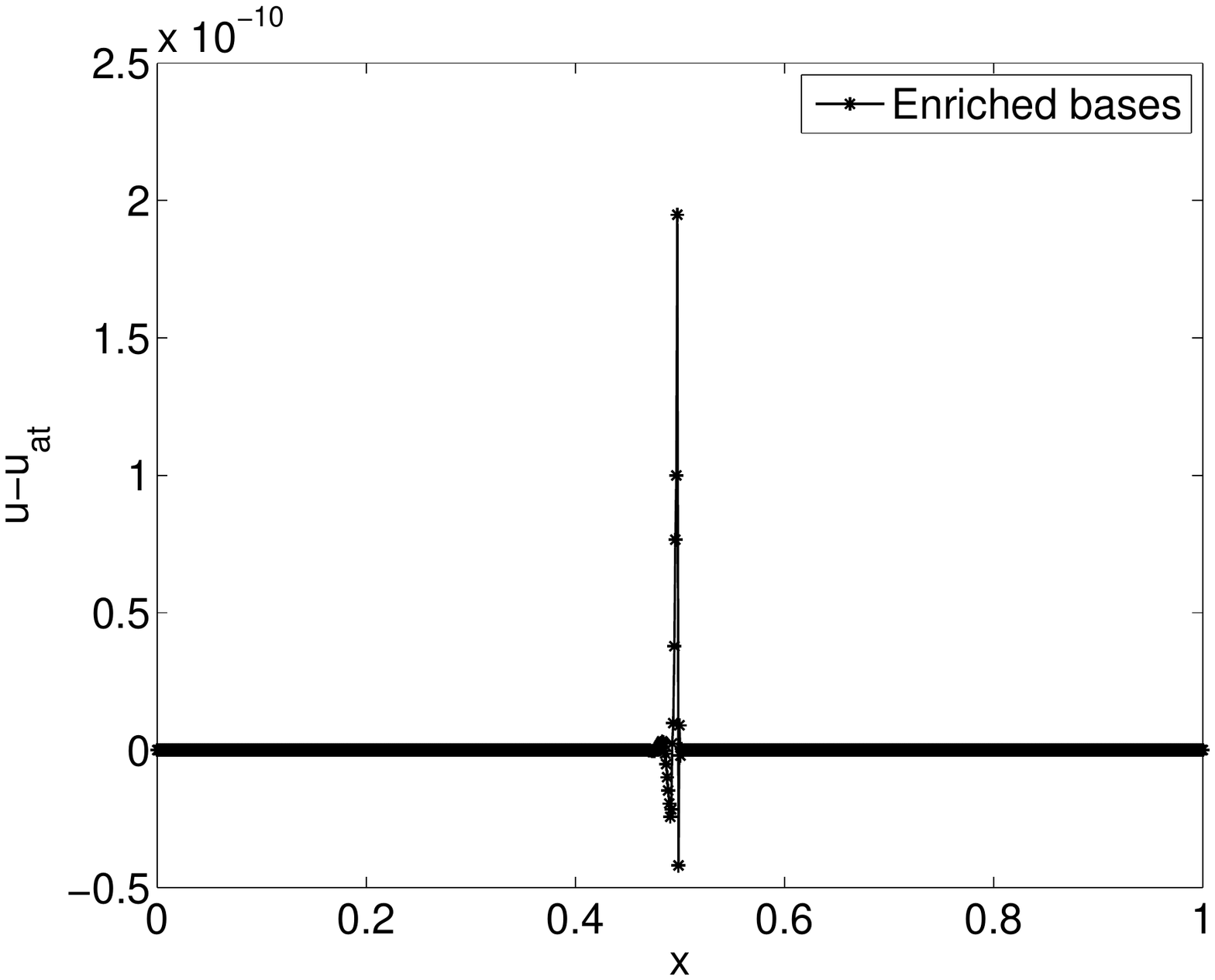}}%
\subfigure[$\ell=2$]{\label {fig:DisplacementErrorExtendedVaryL2}
\includegraphics[width=2.5in]{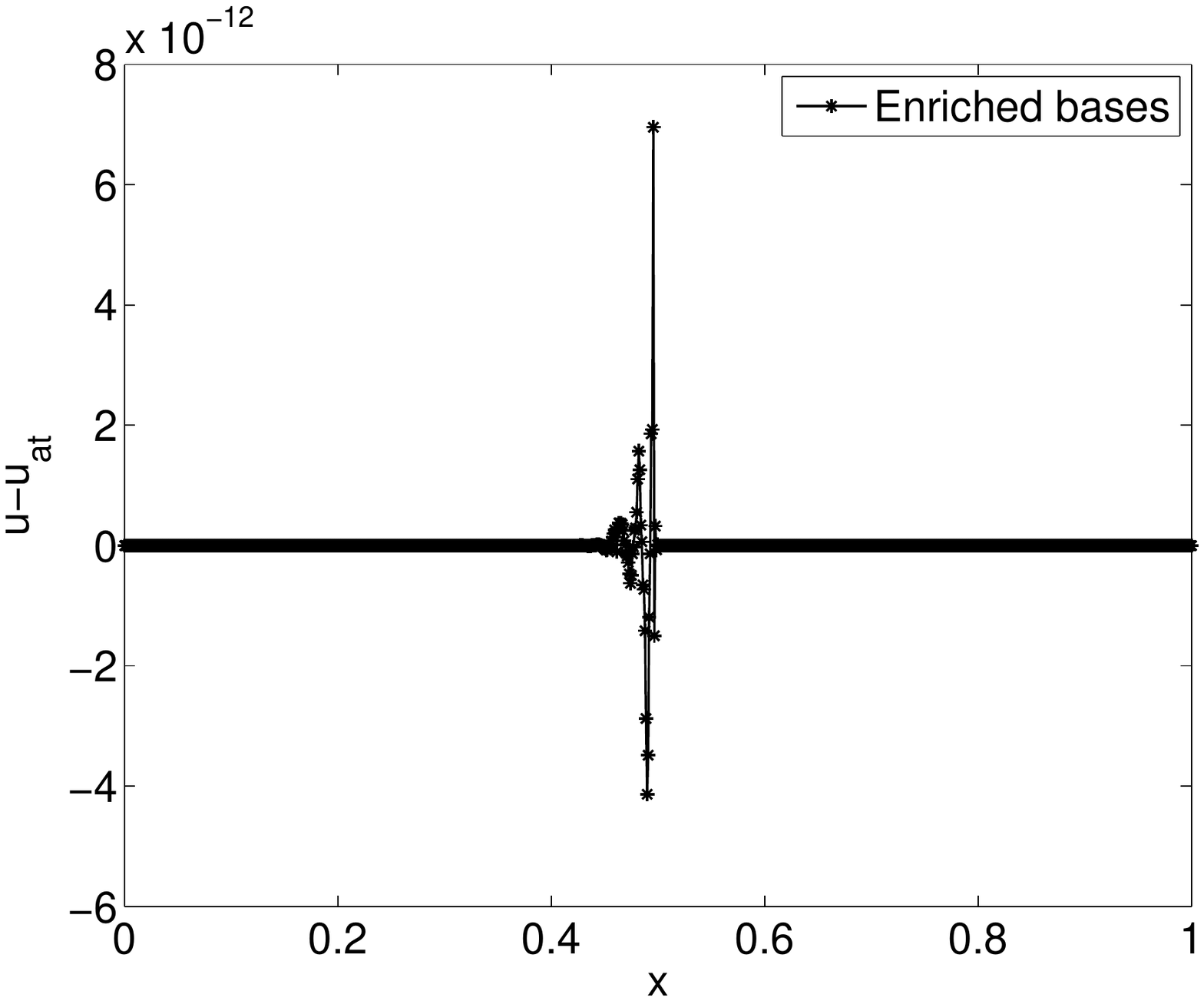}}%
\vspace{-9em}
\subfigure[$\ell=3$]{\label {fig:DisplacementErrorExtendedVaryL3}
\includegraphics[width=2.5in]{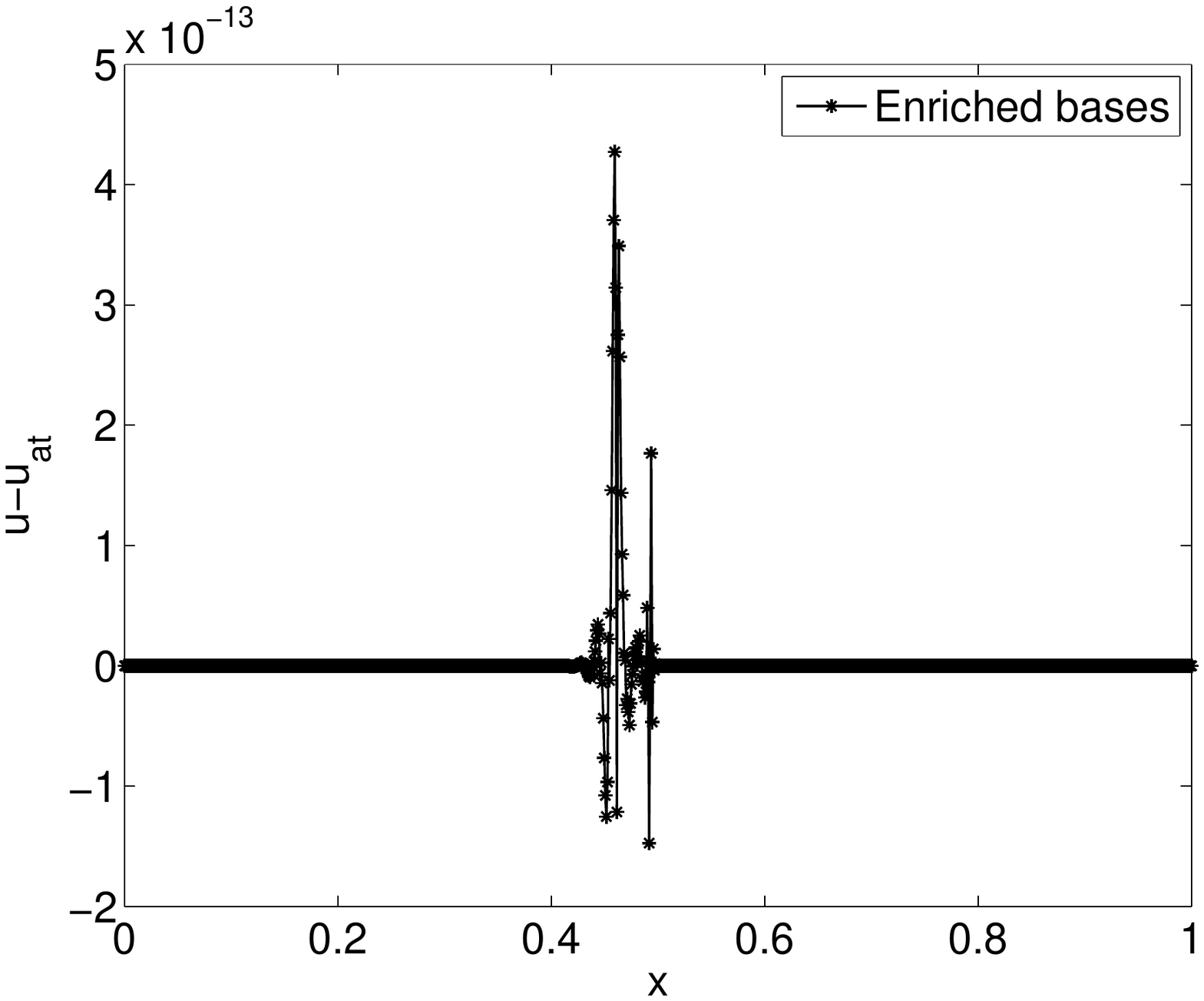}}%
\subfigure[$\ell=4$]{\label {fig:DisplacementErrorExtendedVaryL4}
\includegraphics[width=2.5in]{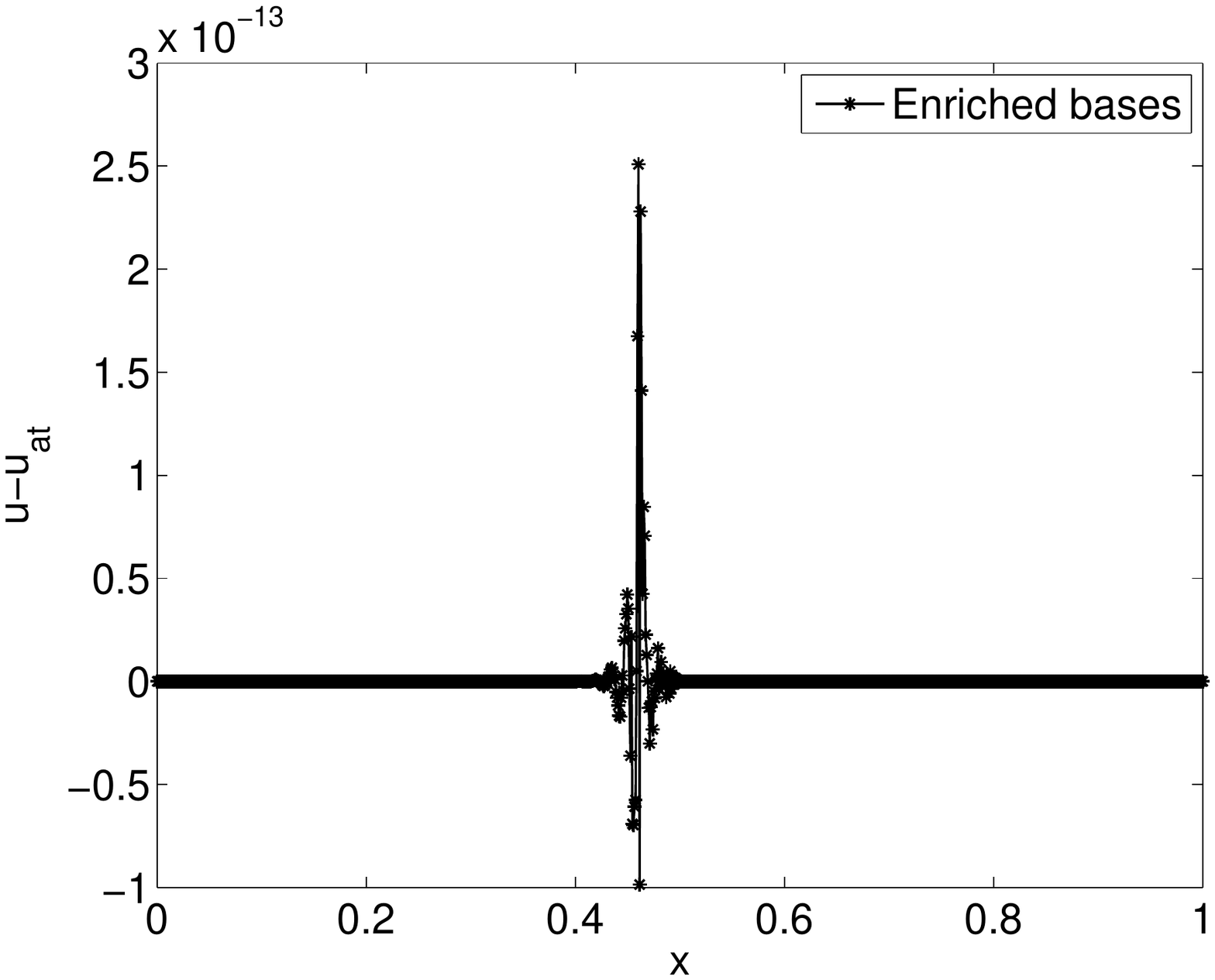}}%
\vspace{-4em}
\caption{\small Displacement errors of the enriched bases method
 with extrapolation for the left boundary. (a): $\ell=1$;
(b): $\ell=2$; (c): $\ell=3$; (d): $\ell=4$. }\label {fig:DisplacementErrorExtendedVaryL}
\end{figure}

Thus far, all the approximate solutions are obtained with the exact quadrature, in which case only the Galerkin projection is responsible for the error. To check the additional error introduced by the quadrature approximation, we implemented the aforementioned Galerkin methods with
the quadrature approximation, and the result, along with the results without quadrature approximation, can be found in
Fig. \ref{fig:DisplacementErrorQuad}. It is clear that the error introduced in the quadrature approximation is much smaller than the approximation error introduced in Galerkin methods. This observation holds true for all examples we studied here.
\begin{figure}[htbp]
\vspace{-6em}
\centering
\subfigcapskip -5em
\subfigure[Standard Galerkin]{\label {fig:DisplacementErrorQuad1}
\includegraphics[width=2.5in]{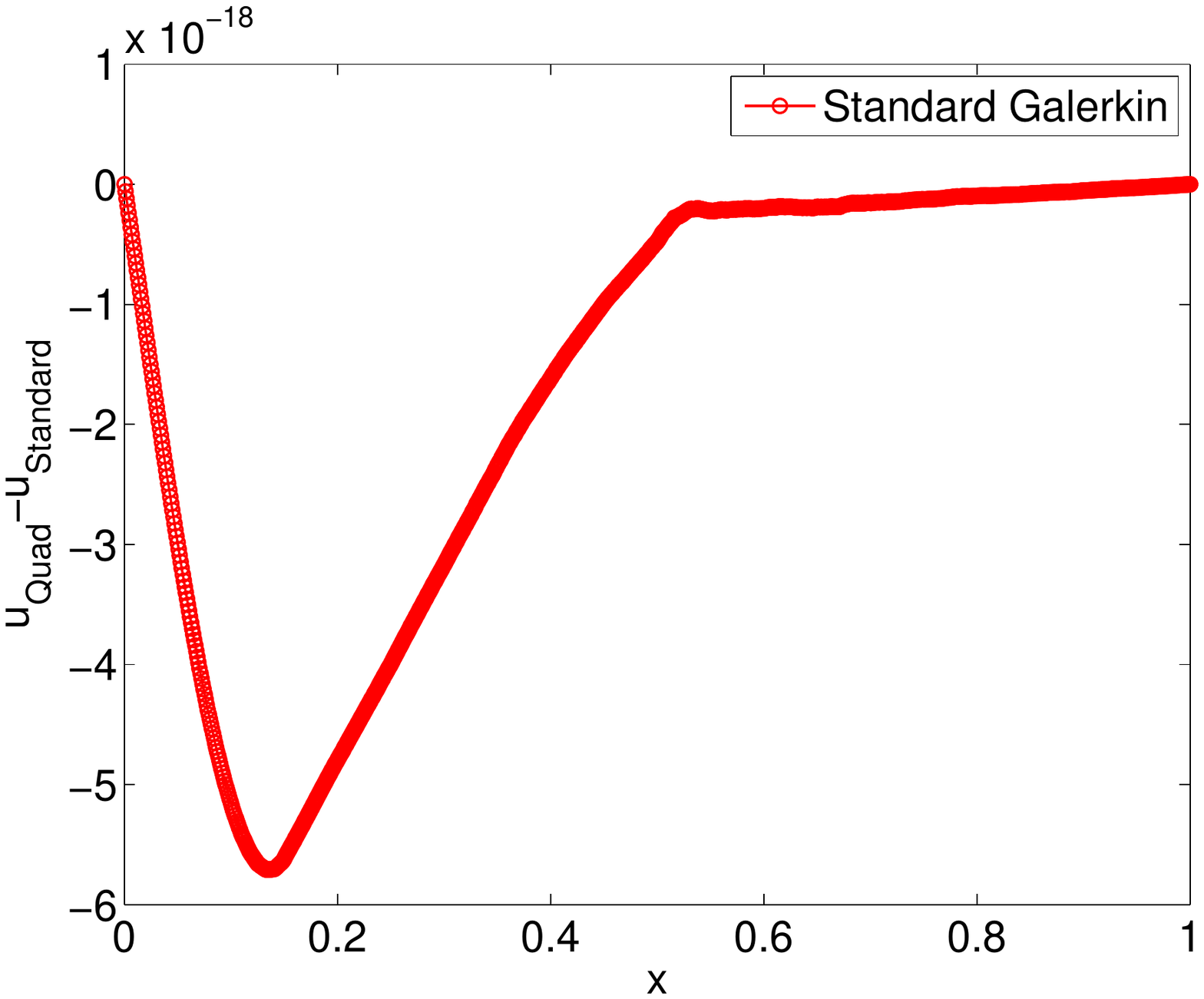}}%
\subfigure[Enriched bases]{\label {fig:DisplacementErrorQuad2}
\includegraphics[width=2.5in]{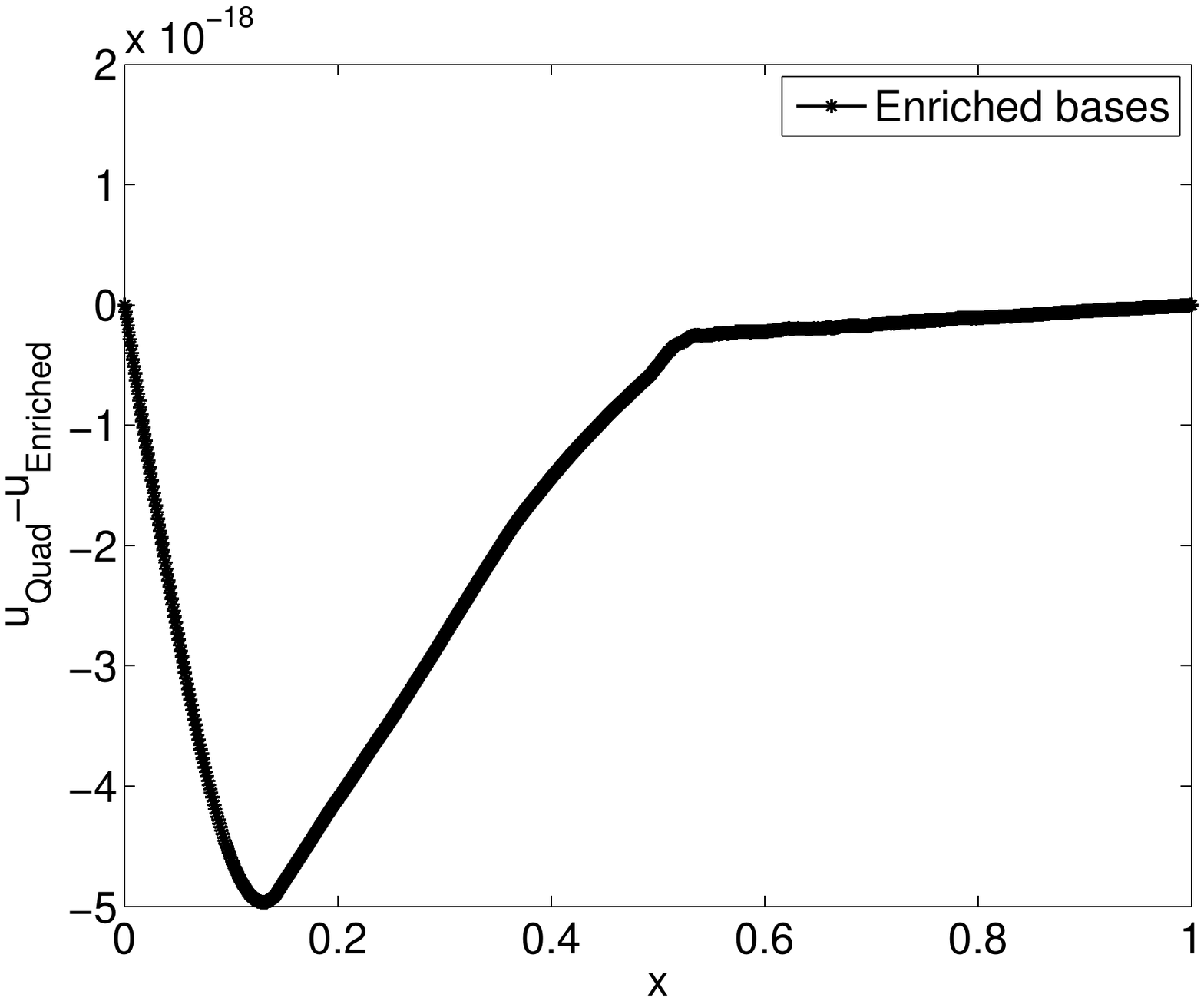}}%
\vspace{-4em}
\caption{\small Error due to the quadrature
approximation. (a): The standard Galerkin method;
(b): The enriched bases method with $\ell=1$. }\label {fig:DisplacementErrorQuad}
\end{figure}

\medskip

Finally, we tested Galerkin methods with a slowly varying external force given by
\begin{equation}\label{eq:fex}
 f^{\text{ex}}(x) =
 \left\{
 \begin{aligned}
   \sin(\pi x), &\quad \text{if}\; 1/2<x\leq 1,\\
   0, &\quad \text{if}\; 0\leq x\leq 1/2.
 \end{aligned}
 \right.
\end{equation}
This is smooth in the sense that $\bmf^{\text{ex}}\in\text{Range}(\Phi^T)$. The exact solution is shown in Fig. \ref{fig:DisplacementInhomogeneous} with displacement gradient in Fig. \ref{fig:DisplacementGradient1dInhomogeneous}.
In addition, errors of the two methods are shown in
 Figs. \ref{fig:DisplacementErrorStandardInhomogeneous} and \ref{fig:DisplacementErrorExtendedInhomogeneous}. Again better performance is observed for the enriched bases method.
\begin{figure}[htbp]

\vspace{-6em}
\centering
\subfigcapskip -5em
\subfigure[Displacement]{\label{fig:Displacement1dInhomogeneous}
\includegraphics[width=2.5in]{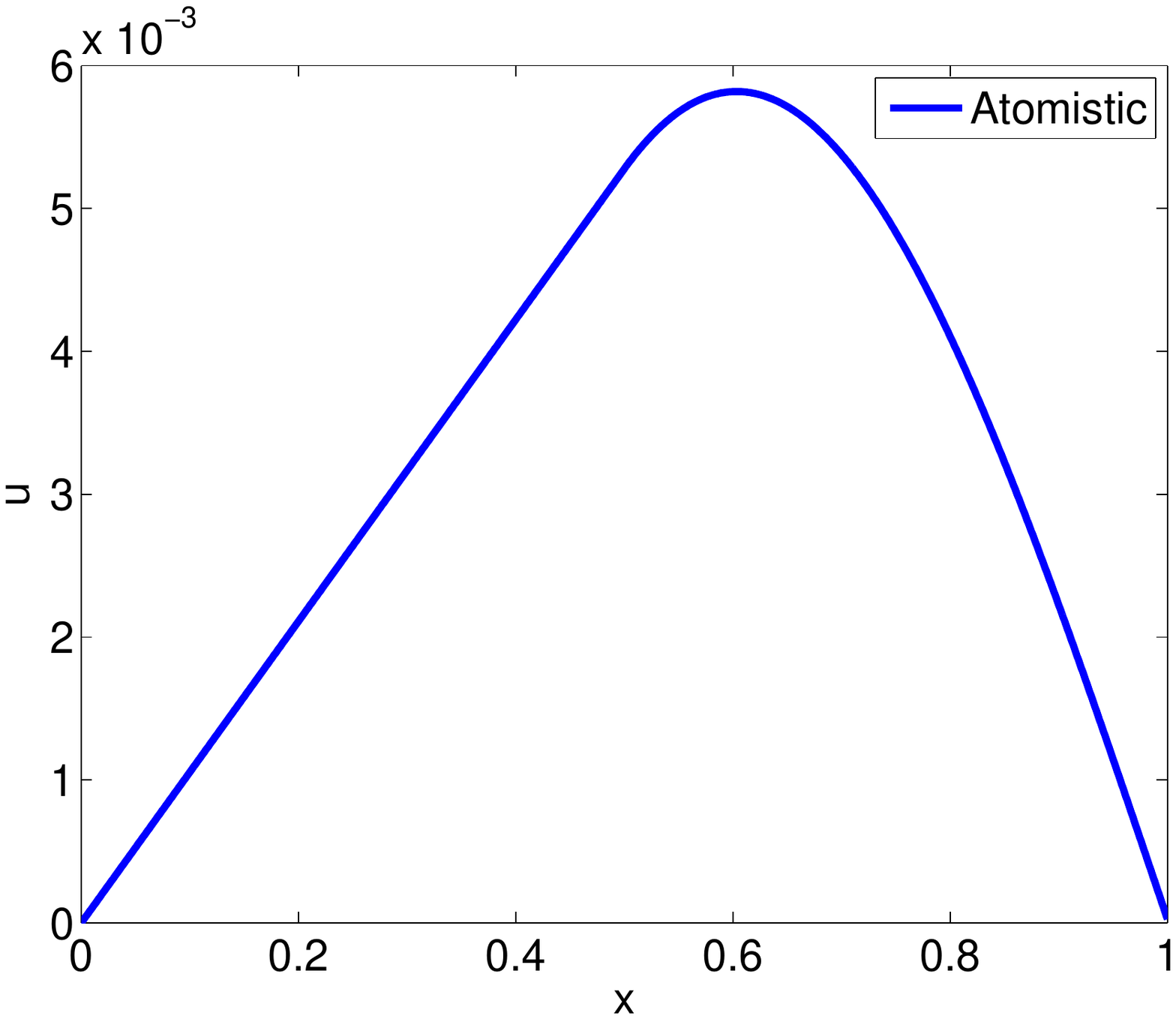}}%
\subfigure[Displacement gradient]{\label{fig:DisplacementGradient1dInhomogeneous}
\includegraphics[width=2.5in]{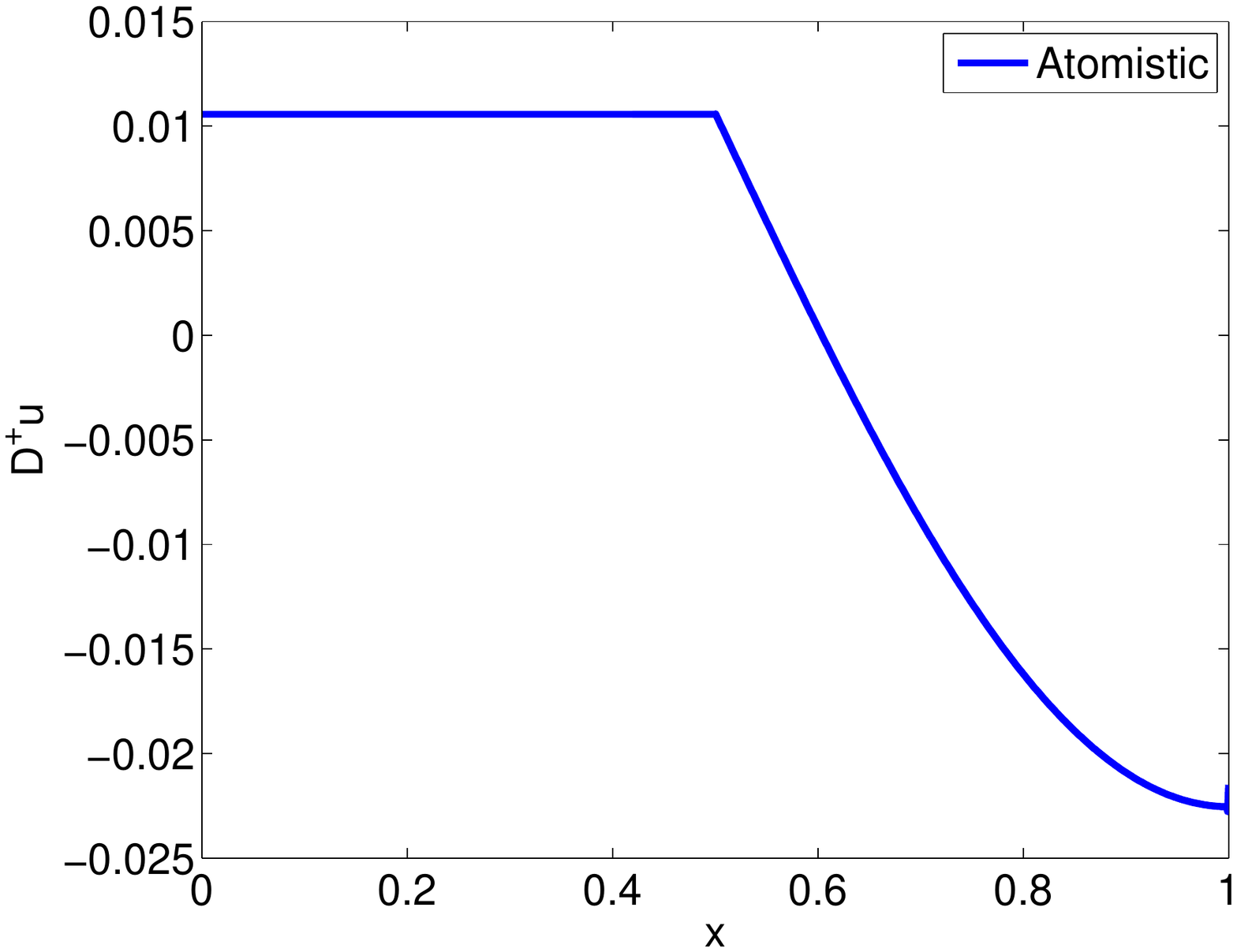}}%
\vspace{-9em}
\subfigure[Standard Galerkin]{\label {fig:DisplacementErrorStandardInhomogeneous}
\includegraphics[width=2.5in]{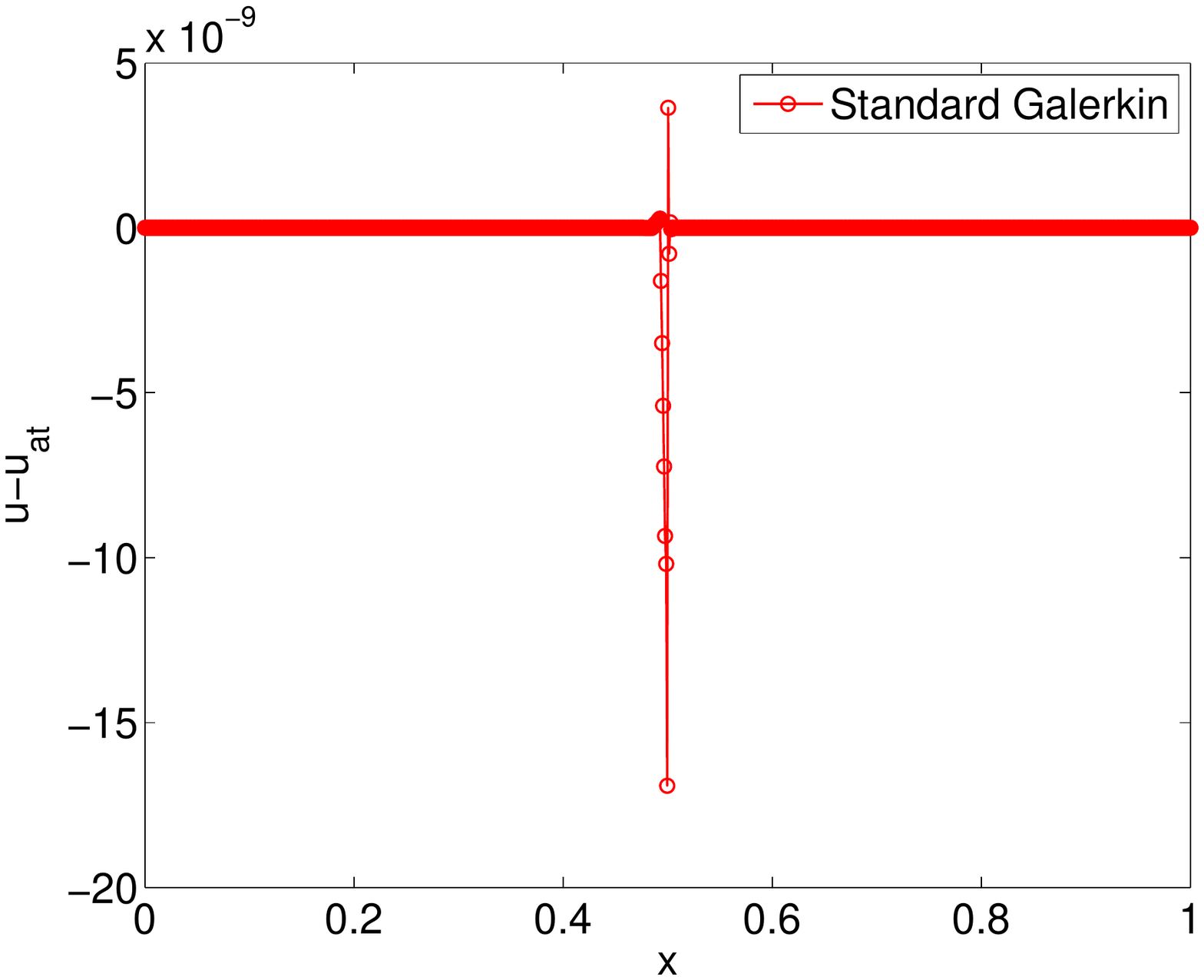}}%
\subfigure[Enriched bases]{\label {fig:DisplacementErrorExtendedInhomogeneous}
\includegraphics[width=2.5in]{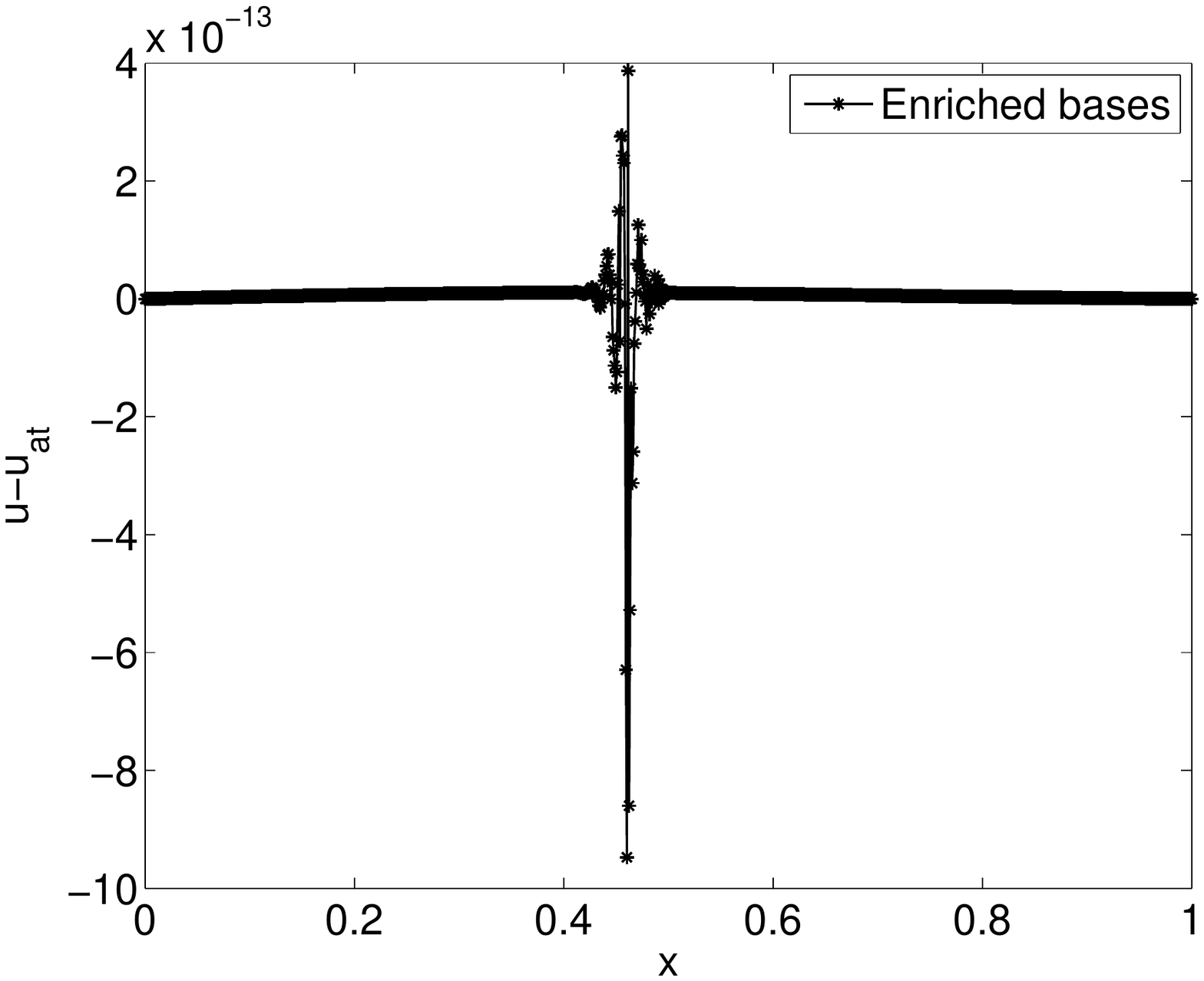}}%
\vspace{-4em}
\caption{\small Displacement and displacement gradient of a system with $1024$ atoms given by the atomistic model with external force given by \eqref{eq:fex},
and displacement errors of Galerkin methods.
 (a): Displacement; (b): Displacement gradient;
(c): The standard Galerkin method; (d): The enriched bases method with $\ell=4$.}
\label{fig:DisplacementInhomogeneous}
\end{figure}

\subsubsection{Convergence rate of Galerkin methods}
\label{sec:Convergence}
In this section, we investigate the behavior of the error as $\veps \to 0.$
Similar to the example in Fig. \ref{fig:DisplacementDelta}, we consider a case where a point external
force is applied at $x=0.5$, at the interface between the continuum region and the atomistic region.
It is clear that $\|f\|_{L^2([0,1])}\leq C_1$, $\|f\|_{\rm{BV}([0,1])} \le C_2$ and  $\| f\|_{H^1([0,1])}\leq C_3/{\sqrt{\veps}}$, where $C_1$, $C_2$  and $C_3$ are constants independent of $\veps$.

The following methods will be considered: The standard Galerkin method with a uniform mesh, the standard Galerkin method with a nonuniform
mesh, the quasi-nonlocal QC method \cite{Shimokawa:2004}, the force-based QC method \cite{Ming:2008, DobsonLuskinOrtner:2010a},
the enriched bases method with a uniform mesh,
and the enriched bases method with a nonuniform mesh. All the meshes considered here are fixed as the system
size is enlarged. All atoms are chosen as rep-atoms in the quasi-nonlocal QC method  and the force-based
QC method.

Fig. \ref{fig:ErrorConvergence} shows the rate of convergence of the first four methods in $W^{1,1}$ norm, $H^1$ norm, and $W^{1,\infty}$
norm, respectively. All the methods have a first-order convergence in $W^{1,\infty}$ norm, even
though the external force considered here is only in $L^2([0,1])\cap\rm{BV}([0,1])$. Galerkin methods and the quasi-nonlocal QC
method have a 3/2-order convergence
in $H^1$ norm, while the force-based QC method only has a first-order convergence in $H^1$ norm. Galerkin methods and the quasi-nonlocal QC method have a second-order convergence
in $W^{1,1}$ norm, while the force-based QC method has a first-order convergence in $W^{1,1}$ norm.

These results are somewhat unexpected since the external force is only in $L^2([0,1])\cap\rm{BV}([0,1])$, while
higher regularities are required for theoretical studies of the quasi-nonlocal QC method
\cite{MingYang:2009, DobsonLuskin:2009b} and the force-based QC method \cite{Ming:2008, DobsonLuskinOrtner:2010a}.
It has been proven in \cite{MingYang:2009} that, for the quasi-nonlocal QC method, if $f\in W^{m,p},\;p\ge 1,\; m\ge 4$ then a uniform convergence rate in $W^{1,\infty}$ is obtained; and in \cite{DobsonLuskin:2009b} that if the solution of the atomistic model is in $W^{4,2}$ then the convergence rate in the $W^{1,p}$ norm is $\veps^{1+1/p}$ for $1\leq p\leq \infty$. For the force-based QC method, it has also been proven that if $f\in W^{m,p},\;p\ge 1,\; m\ge 4$ one gets a uniform convergence rate in $W^{1,\infty}$ \cite{Ming:2008}, and if the solution of the atomistic model is in $W^{3,\infty}$, then the convergence rate in $W^{1,\infty}$ is quadratic. In this example $f\in L^2([0,1])\cap{\rm BV}([0,1])$, and it is not yet clear how the theoretical analysis can be extended to this case.

Quantitatively, the standard Galerkin method on a nonuniform mesh provides a better approximation to the atomistic model. Even though the same convergence rate is
obtained in the $W^{1,\infty}$ norm for these methods, the prefactor in the standard Galerkin method on a nonuniform mesh is smaller that those in the quasi-nonlocal QC method and the force-based QC method.
The prefactor in the Galerkin method with a nonuniform mesh can be one order of magnitude smaller than that in the quasi-nonlocal QC method and a couple of orders of magnitude smaller than that in the force-based QC method.
We can also see that the Galerkin method has a larger prefactor than those in the quasi-nonlocal QC method and the force-based QC method if a uniform mesh is used for $W^{1,\infty}$ norm. This suggests
that a smooth transition of the mesh size from the atomistic region to the continuum region is desirable to reduce the error. After the quadrature approximation, atoms near the interface in the continuum region work as quasiatoms,
which might explain why the corresponding Galerkin method outperforms other methods.
\begin{figure}[htbp]

\vspace{-4em}
\centering
\subfigcapskip -3em
\subfigure[$W^{1,1}$]{\label {fig:ErrorConvergenceW11}
\includegraphics[width=1.6in]{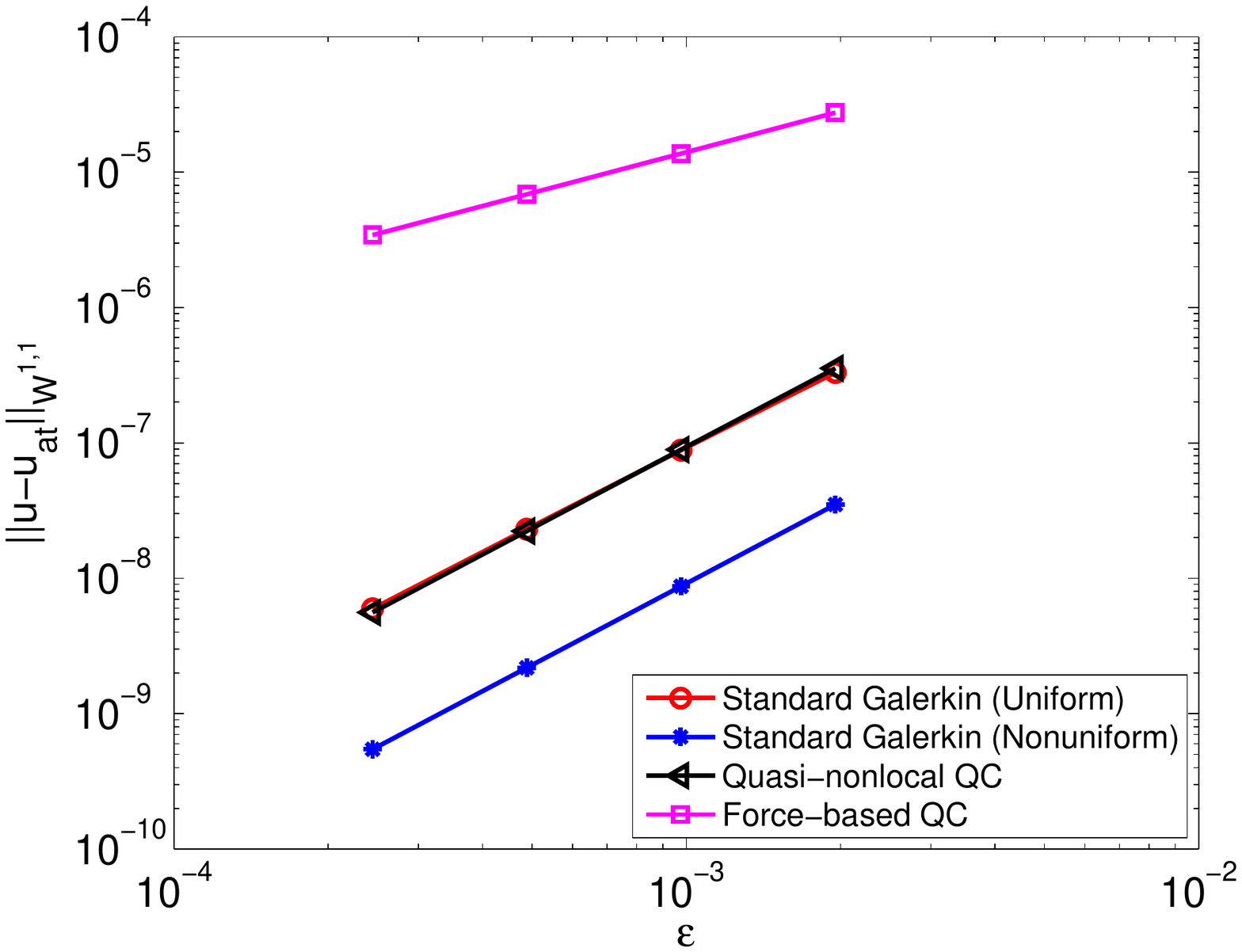}}%
\subfigure[$H^1$]{\label {fig:ErrorConvergenceH1}
\includegraphics[width=1.6in]{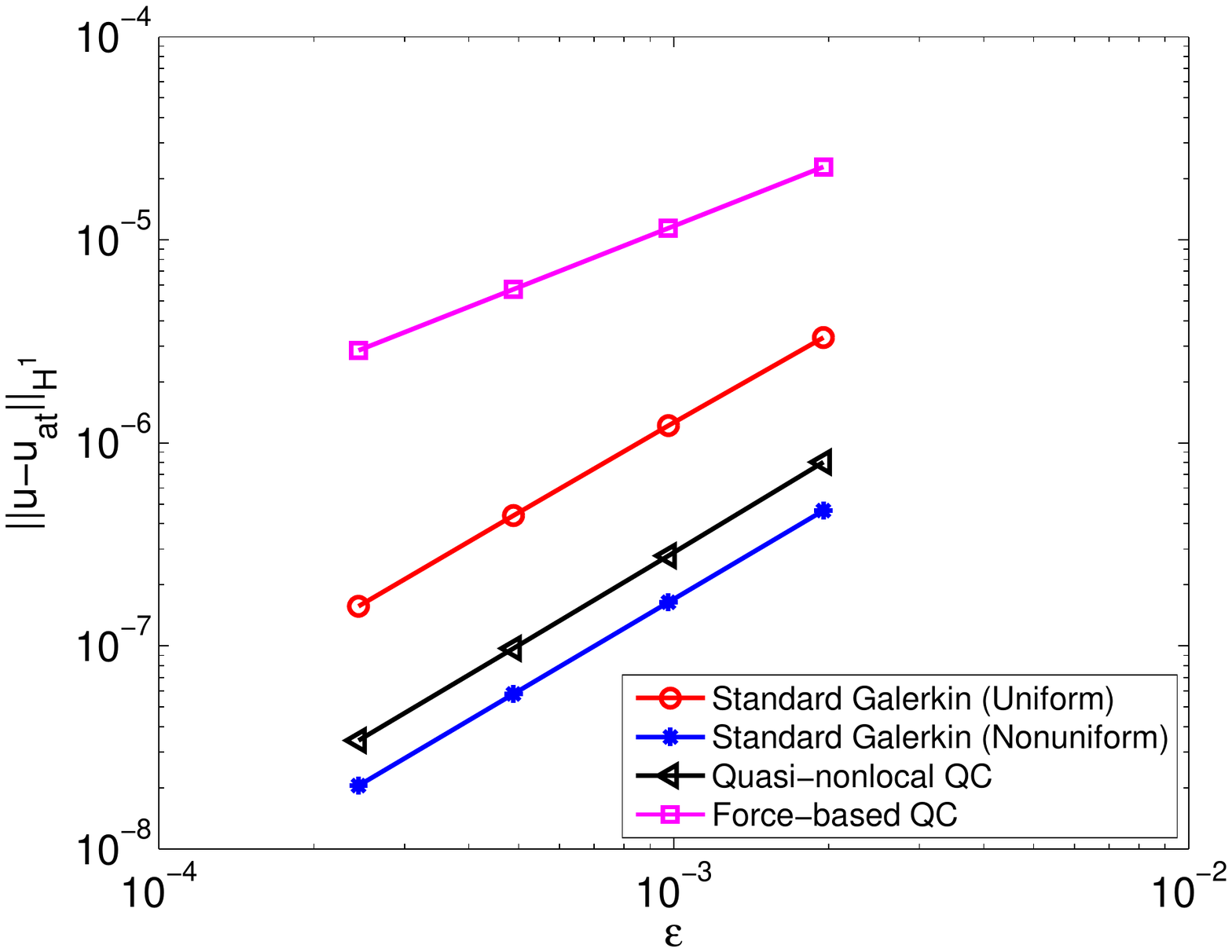}}%
\subfigure[$W^{1,\infty}$]{\label {fig:ErrorConvergenceW1infty}
\includegraphics[width=1.6in]{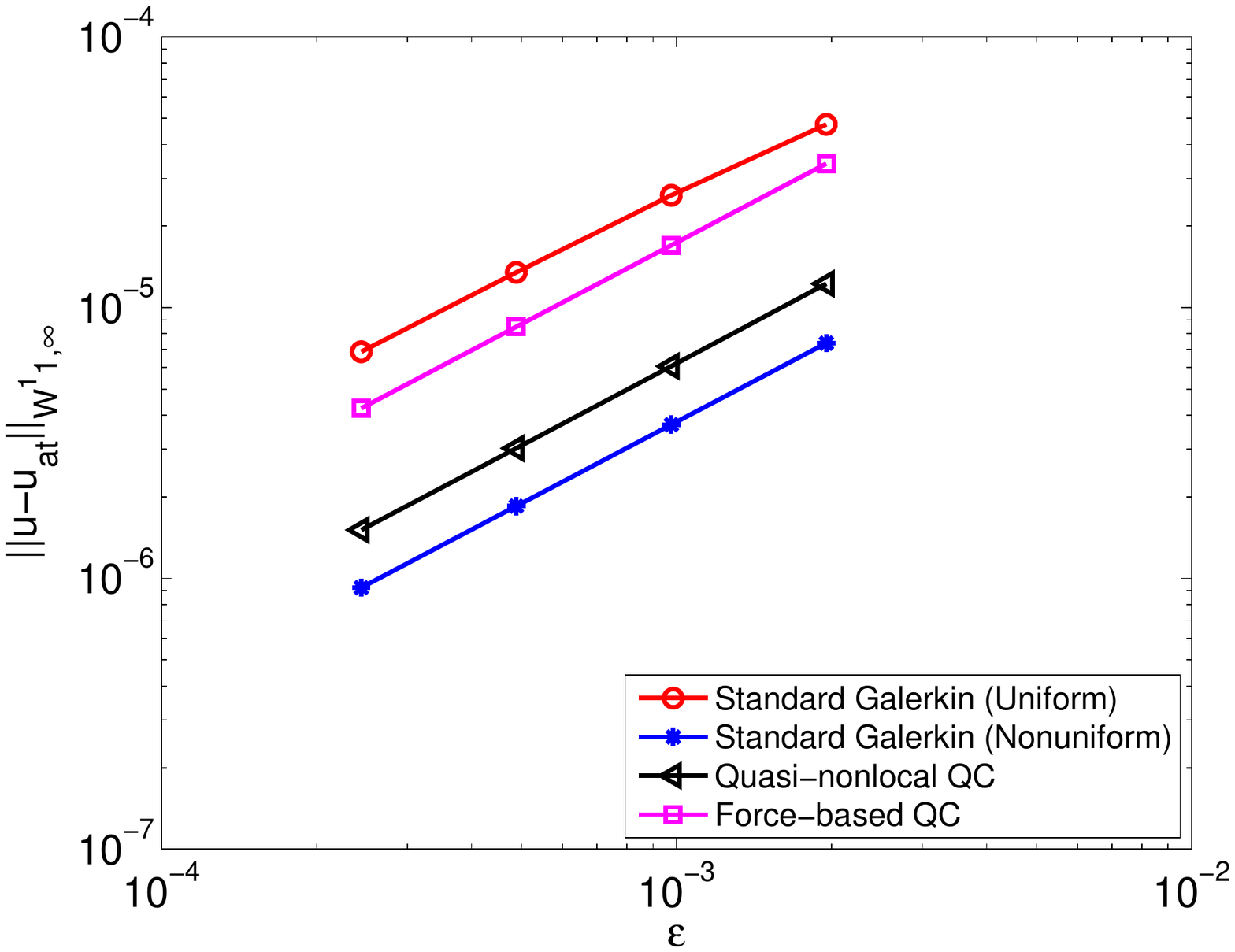}}%
\vspace{-2em}
\caption{\small Errors of the standard Galerkin method with a uniform method, the standard Galerkin method with a nonuniform mesh, the quasi-nonlocal QC method, and the force-based QC method. (a): $W^{1,1}$ norm. Convergence rates are $1.93$, $2.00$, $2.00$, and $1.00$, respectively; (b): $H^1$ norm. Convergence rates are $1.47$, $1.50$, $1.52$, and $1.00$, respectively; (c): $W^{1,\infty}$ norm. Convergence rates are $0.93$, $1.00$, $1.01$, and $1.00$, respectively. }\label {fig:ErrorConvergence}
\end{figure}

In Table \ref{tab:extendeduniform} we show the $W^{1,1}$, $H^1$, and $W^{1,\infty}$ norms of the error for the enriched bases method with a uniform mesh. With a fixed number of enriched bases $(m=20)$ around
the interface, the error drops quickly at the beginning and then follows corresponding convergence
rates of the method in three norms. This is mainly due the to deficiency of enriched bases for smaller $\veps$. Moreover, the prefactor in this example is
roughly one order of magnitude smaller than that of the standard Galerkin method with a uniform mesh in
the previous example, which shows the technique of enrichment is still an efficient way to reduce the error.
\begin{table}
\centering
\begin{tabular}{|c|c|c|c|c|}
\hline
$\veps$  & $1/512$ & $1/1024$ & $1/2048$ & $1/4096$ \\ \hline
$W^{1,1}$  &  $9.94e-17$  & $1.91e-08$ & $7.66e-09$ & $1.88e-09$  \\ \hline
$H^1$  &  $1.14e-16$  & $1.95e-07$ & $1.03e-07$ & $3.67e-08$  \\ \hline
$W^{1,\infty}$  &  $1.53e-15$  & $3.98e-06$ & $2.21e-06$ & $1.10e-06$  \\ \hline
$m$ &  $20$  & $20$ & $20$ & $20$\\ \hline
\end{tabular}
\caption{$W^{1,1}$, $H^1$, and $W^{1,\infty}$ norms of the error for the enriched bases method
with a uniform mesh. With a fixed number of enriched bases $20$ around
the interface, the error drops quickly at the beginning and then follows corresponding convergence
rates for the method in three norms. The reduction of the error slows down due to the deficiency
of enriched bases for smaller $\veps$.}
\label{tab:extendeduniform}
\end{table}

In Table \ref{tab:extendednonuniform} we list the $W^{1,1}$, $H^1$, and $W^{1,\infty}$ norms of the error for the enriched bases method with a nonuniform mesh. As $\veps$ becomes smaller, a larger number of enriched bases around the interface are required to keep the error almost independent of $\veps$. The small variation of the error in this example is mainly because the enrichment here is not done adaptively. As the number of bases
increases, the profile of the error spreads out around the interface. The number of enriched bases is roughly proportional to $1/|\veps\log\veps|$.
\begin{table}
\centering
\begin{tabular}{|c|c|c|c|c|}
\hline
$\veps$  & $1/512$ & $1/1024$ & $1/2048$ & $1/4096$ \\ \hline
$W^{1,1}$  &  $2.61e-16$  & $1.82e-16$ & $1.24e-14$ & $4.09e-15$  \\ \hline
$H^1$  &  $9.20e-16$  & $4.69e-16$ & $4.90e-14$ & $1.47e-14$  \\ \hline
$W^{1,\infty}$  &  $8.38e-15$  & $5.83e-15$ & $7.09e-13$ & $5.29e-13$  \\ \hline
$m$ &  $36$  & $81$ & $162$ & $364$\\ \hline
\end{tabular}
\caption{$W^{1,1}$, $H^1$, and $W^{1,\infty}$ norms of the error for the enriched bases method
with a nonuniform mesh. As $\veps$ becomes smaller, a larger number of enriched bases
around the interface is required to keep the error almost
independent of $\veps$. The number of enriched bases is roughly proportional to $1/|\veps\log\veps|$.}
\label{tab:extendednonuniform}
\end{table}

\subsubsection{A slowly varying external force}
We consider a slowly varying external force given by
\begin{equation}\label{eq:fnonlocal}
 f^{\text{ex}}(x) =
   \sin(\pi x), \quad 0\leq x\leq 1,
 \end{equation}
which does not satisfy the condition $\bmf^{\text{ex}}\in\text{Range}(\Phi^T)$.
This force vanishes at both end points.
The exact solution is shown in Fig. \ref{fig:Displacement1dnonlocal} with displacement gradient
in Fig. \ref{fig:DisplacementGradient1dnonlocal}. Since piecewise linear functions are used to
construct $\Phi$, the exact solution cannot be well approximated in the continuum region by
the standard Galerkin method; see Fig. \ref{fig:DisplacementErrorStandardnonlocal}.
Consequently, the difference decays algebraically. Nevertheless, as enriched bases are used, the error
is still reduced; see Figs. \ref{fig:DisplacementErrorExtendednonlocal1}
and \ref{fig:DisplacementErrorExtendednonlocal2}. By comparing Figs \ref{fig:DisplacementErrorExtendednonlocal1}
and \ref{fig:DisplacementErrorExtendednonlocal2}, we observe enriched bases in a larger support
are required to reduce the error significantly.
Overall, better performance is observed for the enriched bases method again.
\begin{figure}[htbp]

\vspace{-6em}
\centering
\subfigcapskip -5em
\subfigure[Displacement]{\label{fig:Displacement1dnonlocal}
\includegraphics[width=2.5in]{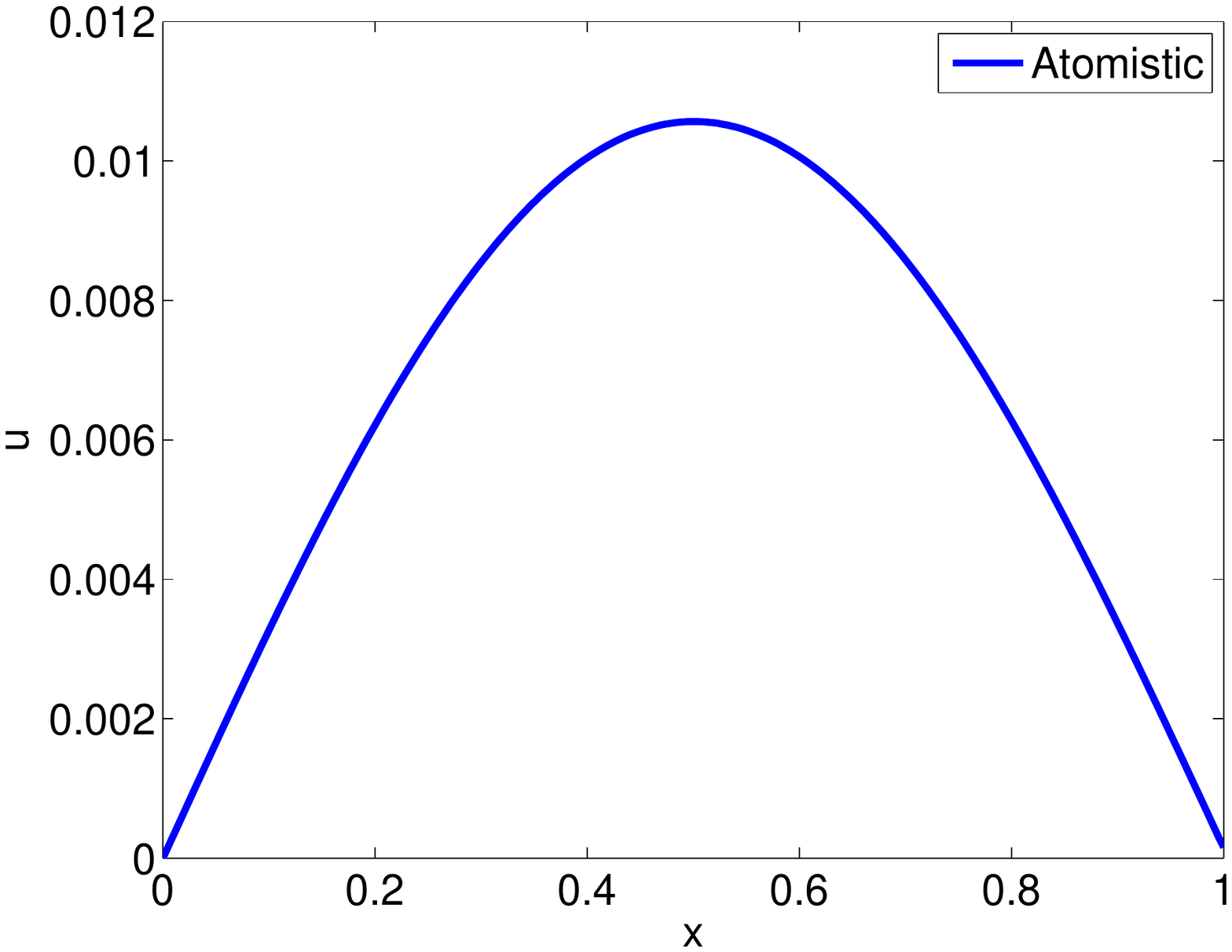}}%
\subfigure[Displacement gradient]{\label{fig:DisplacementGradient1dnonlocal}
\includegraphics[width=2.5in]{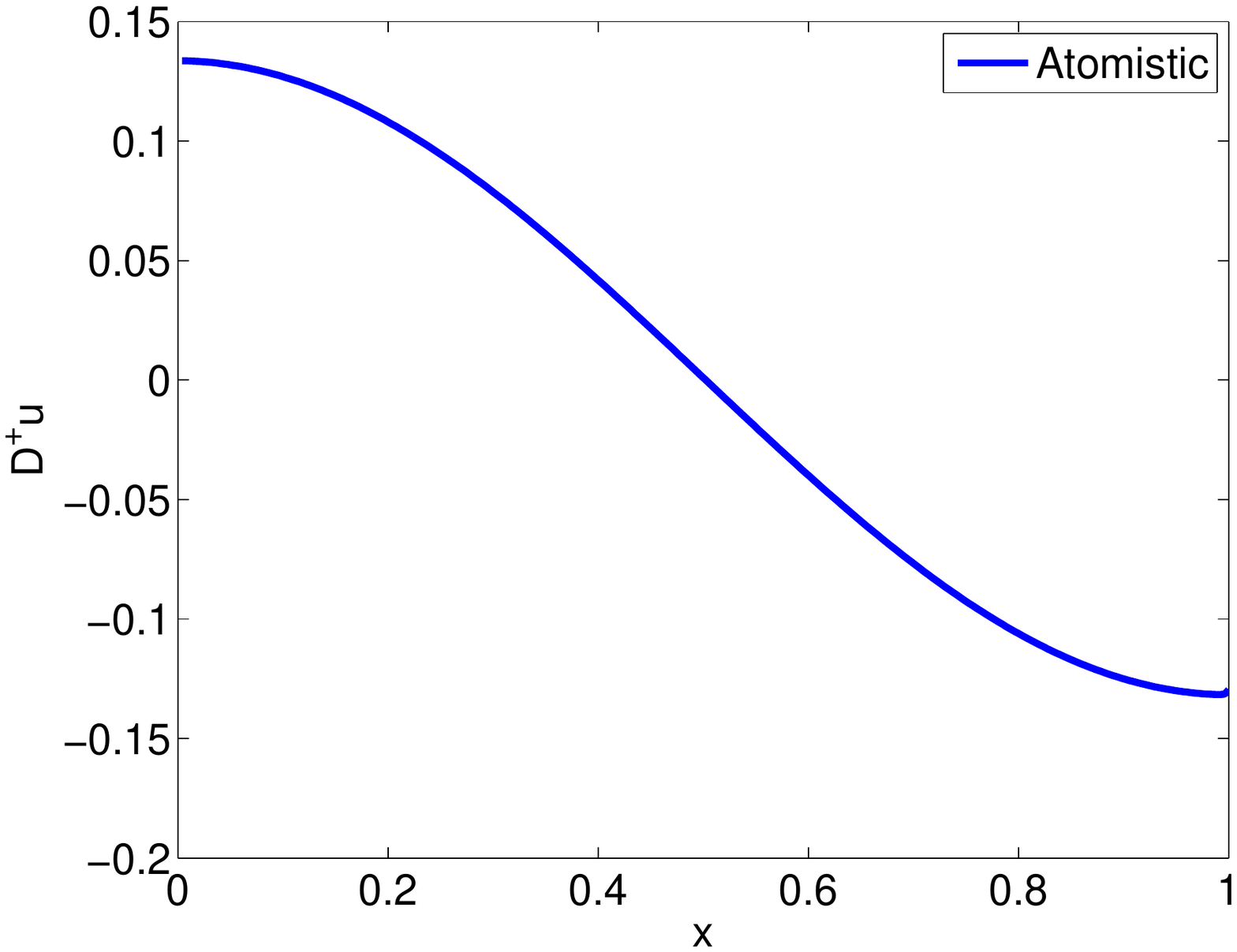}}%
\vspace{-9em}
\subfigure[Standard Galerkin]{\label {fig:DisplacementErrorStandardnonlocal}
\includegraphics[width=2.5in]{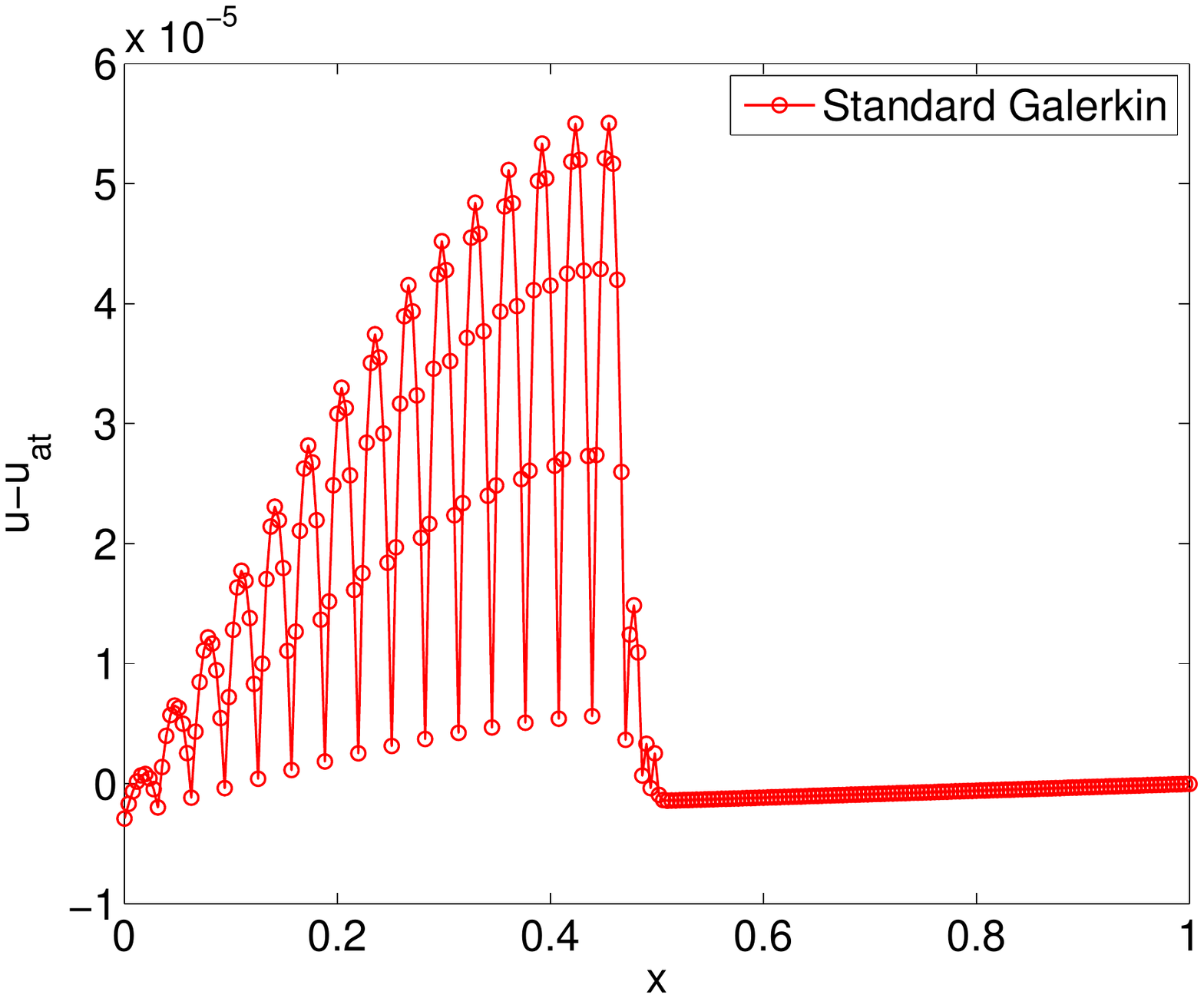}}%
\subfigure[Enriched bases]{\label {fig:DisplacementErrorExtendednonlocal1}
\includegraphics[width=2.5in]{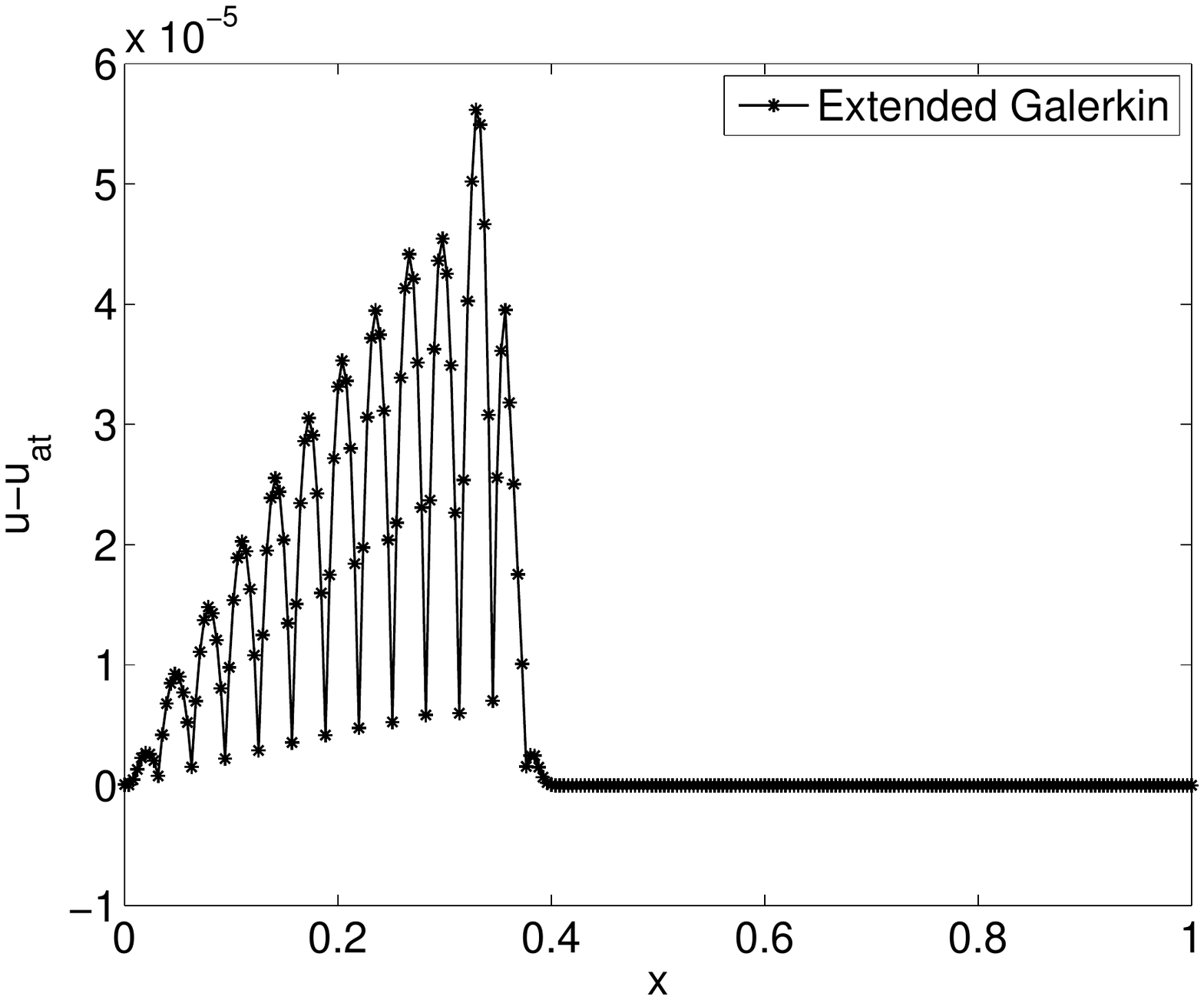}}%
\vspace{-8em}
\subfigure[Enriched bases]{\label {fig:DisplacementErrorExtendednonlocal2}
\includegraphics[width=2.5in]{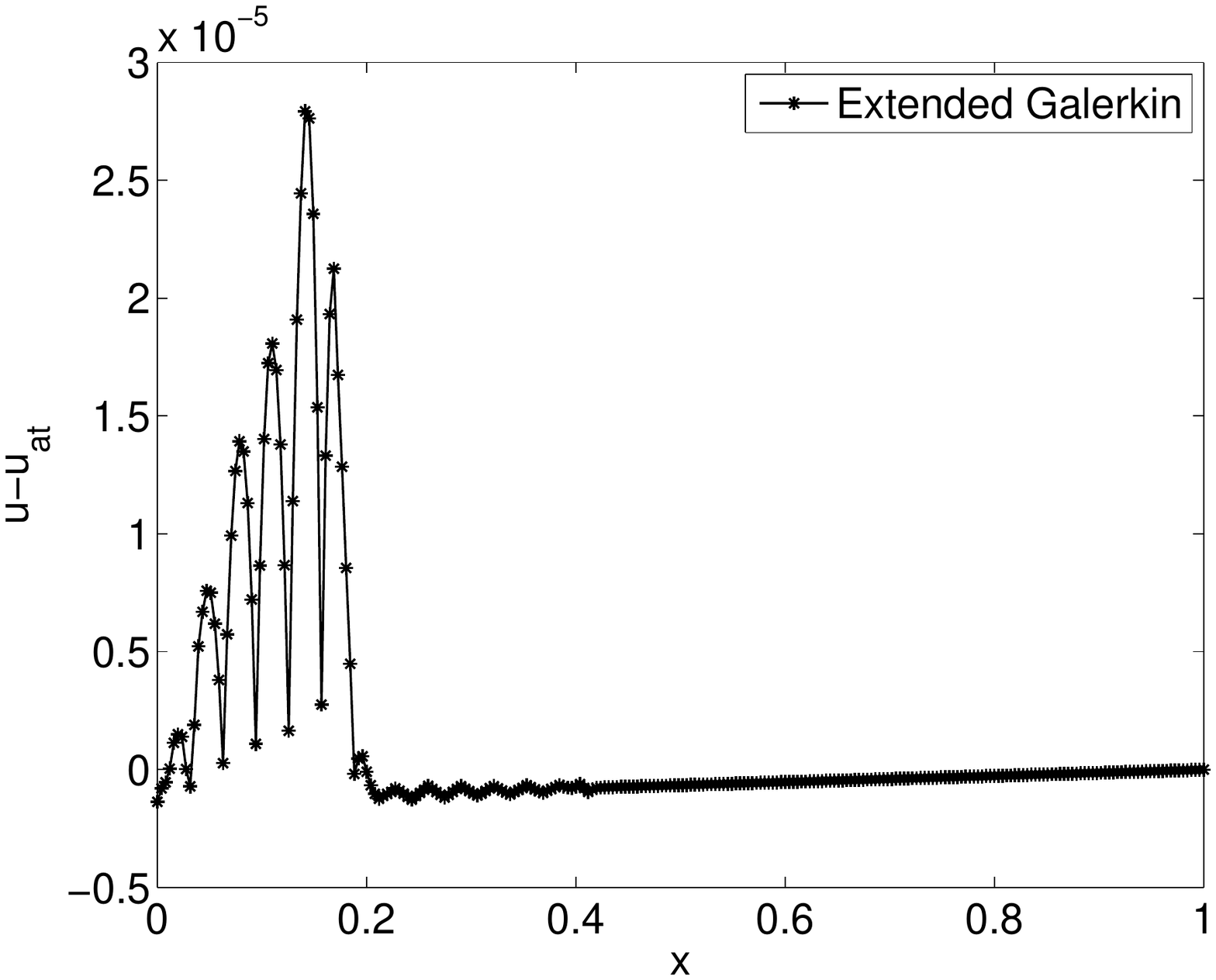}}%
\vspace{-4em}
\caption{\small Displacement and displacement gradient of a system with $256$ atoms given by the atomistic model with external force given by \eqref{eq:fnonlocal},
and displacement errors of Galerkin methods.
 (a): Displacement; (b): Displacement gradient;
(c): The standard Galerkin method; (d): The enriched bases method with 6 nodes selected and
$\ell=4$; (e): The enriched bases method with 12 nodes selected and
$\ell=4$.}
\label{fig:Displacementnonlocal}
\end{figure}

Similar to the convergence study in Section \ref{sec:Convergence}, we also
conducted the convergence study for the slowly varying force \eqref{eq:fnonlocal}.
The standard Galerkin method, the force-based QC method, and the quasi-nonlocal QC
method have first-order convergence in $W^{1,1}$, $H^1$ and
$W^{1,\infty}$ norms.
This suggests convergence rates of these methods are problem-dependent.
A nonlocal force may reduce the convergence rate of a method.
To make the error in the enriched bases method independent of $\veps$,
the number of enriched bases scales $\mathcal{O}(1/|\veps\log\veps|)$ and
the support of these enriched bases is much larger than that in
Section \ref{sec:Convergence}.
}

\subsection{An one-dimensional crack model}
In this section, we consider a one-dimensional crack model studied in \cite{LiMing:2014}, where a detailed description of the
model can be found. The model is adopted from the lattice models \cite{ThHsRa71}, for which the issue of lattice trapping is first addressed. It contains two symmetric chains of atoms above and below an open crack. Due to the symmetry, only the upper chain is considered. Each atom is interacting with two neighboring atoms from the left and two atoms from the right. In addition, each atom in the upper chain is bonded to the corresponding atom in the lower chain with a nonlinear spring, except for the first $n$
bonds, which are considered ``cracked". In spite of its simplicity, the model predicts the correct fracture initiation mechanism \cite{XLi2014,XLi2013A}. The total potential energy for the upper chain reads as
\begin{equation}\label{eq:enfrac}
\begin{split}
V=& -Pu_1 + \sum_{j\ge 1}\left(U \left (\frac{u_{j+1}-u_j}{\veps}\right )+U\left (\frac{u_{j+2}-u_j}{\veps}\right )\right)\\
& +(n-1)\gamma_0 + \gamma(u_n) + K_2\sum_{j>n}\left (\frac{u_j}{\veps} \right)^2
\end{split}
\end{equation}
with
\begin{eqnarray*}
\gamma(u) & = & \int_0^u \frac{K_2}{u_{\text{cut}}^2} v\left (\frac{v-u_{\text{cut}}}{\veps}\right )^2\chi_{[0,u_{\text{cut}}]}(v)\text{d}v,\\
\gamma_0 & = & \gamma(u_{\text{cut}}),
\end{eqnarray*}
and $\chi_{[0,u_{\text{cut}}]}$ the characteristic function.

Galerkin methods are tested for scenarios when $U$ is a harmonic potential with
$K_0=4$ and $K_1=0.4$, and the Lennard-Jones potential. The latter gives similar results.
The minimization problem and the nonlinear systems of equations are solved using Newton's method.
For simplicity, we keep the enriched bases fixed at each step in Newton's method.
We choose $K_2=0.5$, and $u_{\text{cut}}=0.5$.

\begin{figure}[htbp]

\centering
\vspace{-5em}
\includegraphics[width=3.5in]{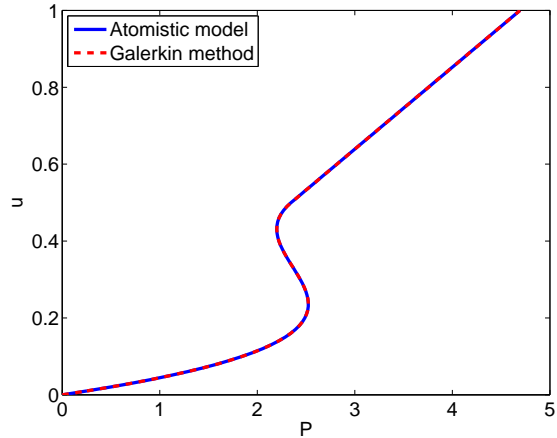}
\vspace{-5em}
\caption{\small Bifurcation diagrams for the atomistic model and the Galerkin method.
The middle branch contains unstable equilibria while the other two branches are stable. }\label {fig:Bifdiagram}
\end{figure}

In the context of fracture mechanics, of particular interest is the crack initiation
criteria. This is manifested in lattice models in the form of a bifurcation. In Fig. \ref{fig:Bifdiagram} we show the corresponding bifurcation
diagrams for the the atomistic model and the Galerkin method, and they agree very well. This implies that
the reduced model obtained from the Galerkin method inherits the correct property of stability transition.
In what follows, we set $P=1$ so that the system only has one solution, which makes the comparison between different methods unambiguous.  For the traction boundary condition, we set forces of the first two atoms to be
\[
f_1 = \frac{K_0+2K_1}{K_0+4K_1} P,\quad f_2 = \frac{2K_1}{K_0+4K_1} P
\]
so that a uniform displacement gradient will be generated near the boundary \cite{LiLu14}. The same idea is applied to the nonlinear case.

\begin{figure}[htbp]

\vspace{-6em}
\centering
\subfigcapskip -5em
\subfigure[Displacement]{\label {fig:DisplacementCrack1}
\includegraphics[width=2.5in]{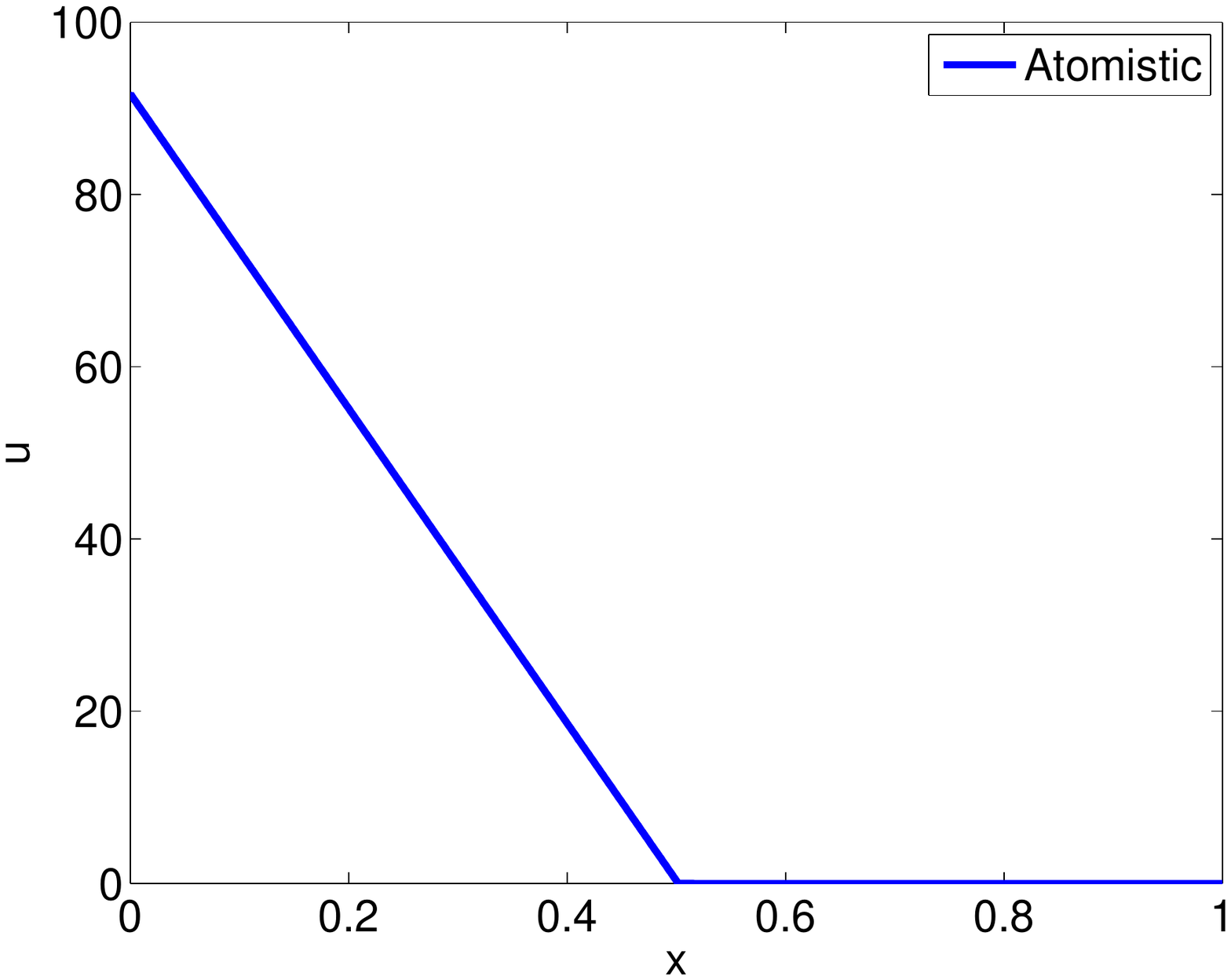}}%
\subfigure[Displacement gradient]{\label {fig:DisplacementCrack2}
\includegraphics[width=2.5in]{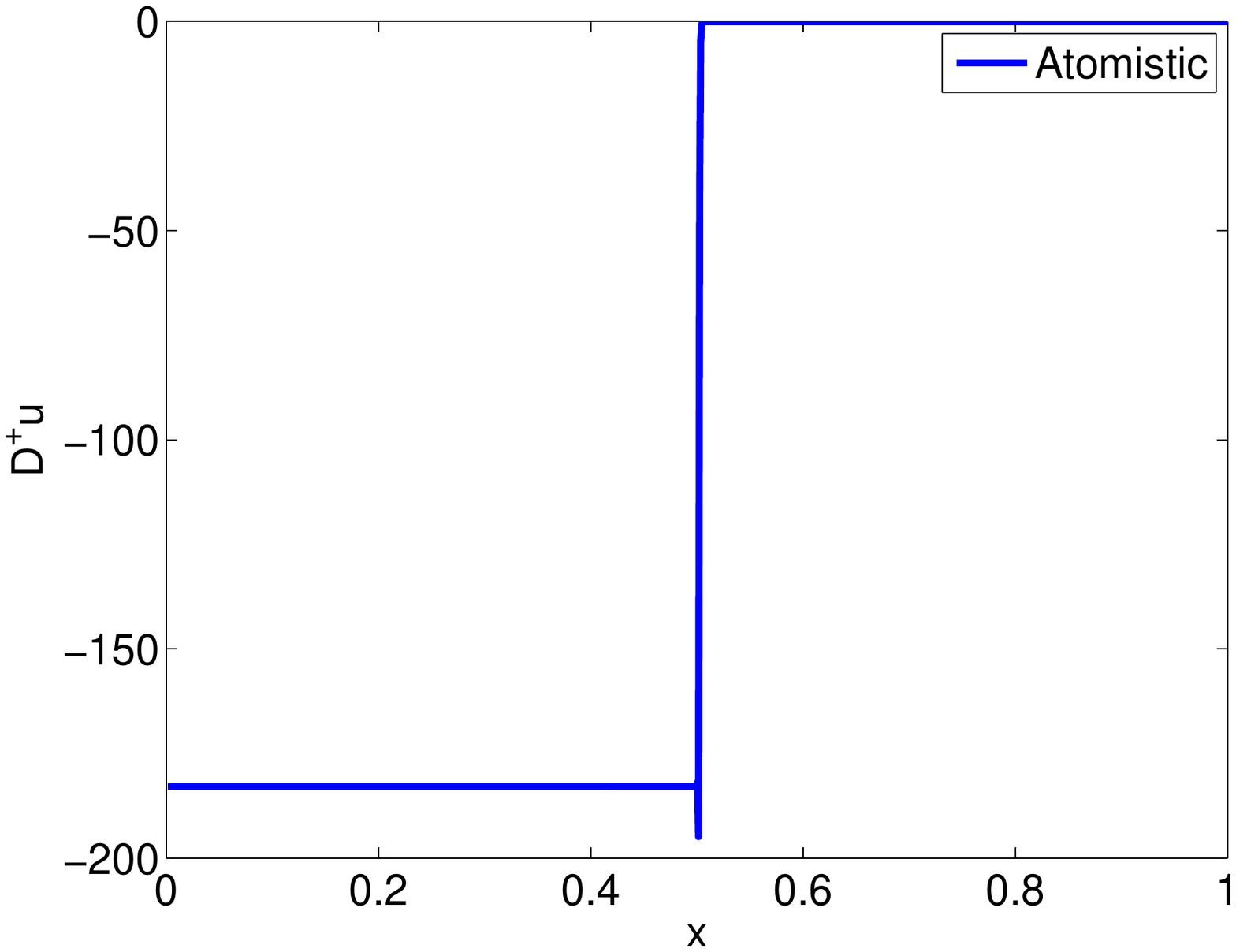}}%
\vspace{-4em}
\caption{\small Atomic displacement, and displacement gradient of the crack model.
(a): Displacement; (b): Displacement gradient.}\label {fig:DisplacementCrack}
\end{figure}

\begin{figure}[htbp]

\vspace{-6em}
\centering
\subfigcapskip -5em
\subfigure[Standard Galerkin]{\label {fig:DisplacementErrorCrack1}
\includegraphics[width=2.5in]{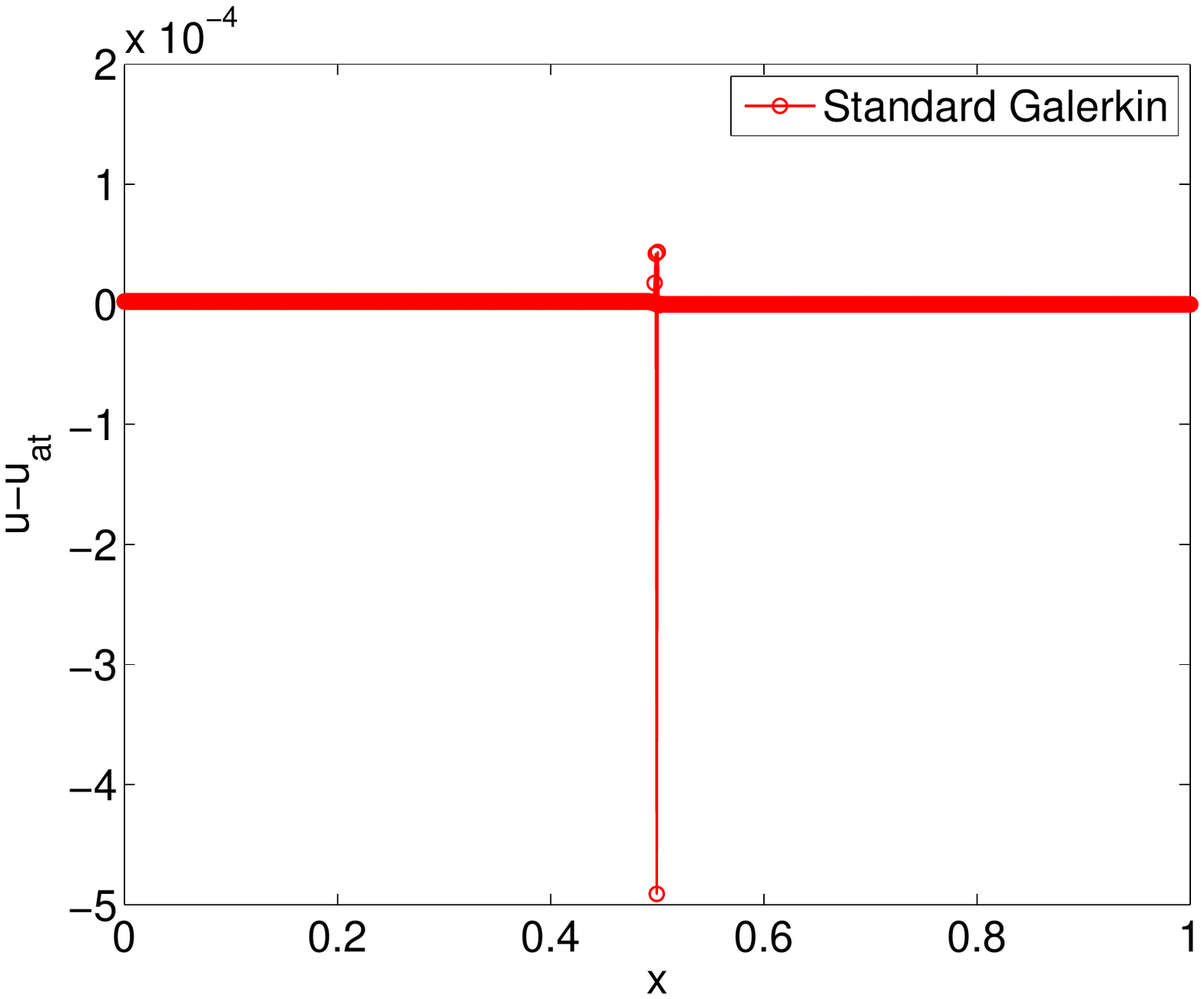}}%
\subfigure[$\ell=1$]{\label {fig:DisplacementErrorCrack2}
\includegraphics[width=2.5in]{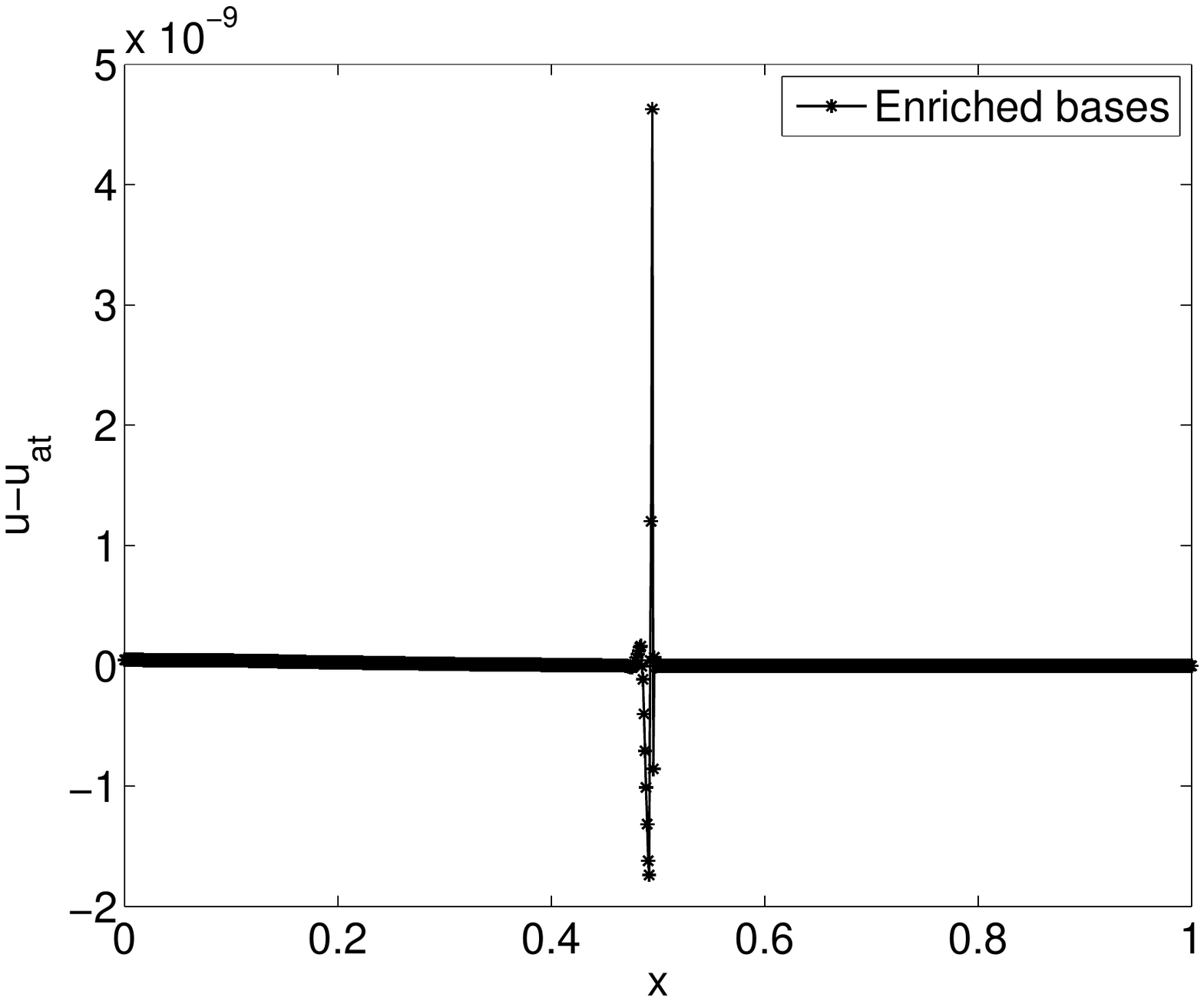}}%
\vspace{-9em}
\subfigure[$\ell=2$]{\label {fig:DisplacementErrorCrack3}
\includegraphics[width=2.5in]{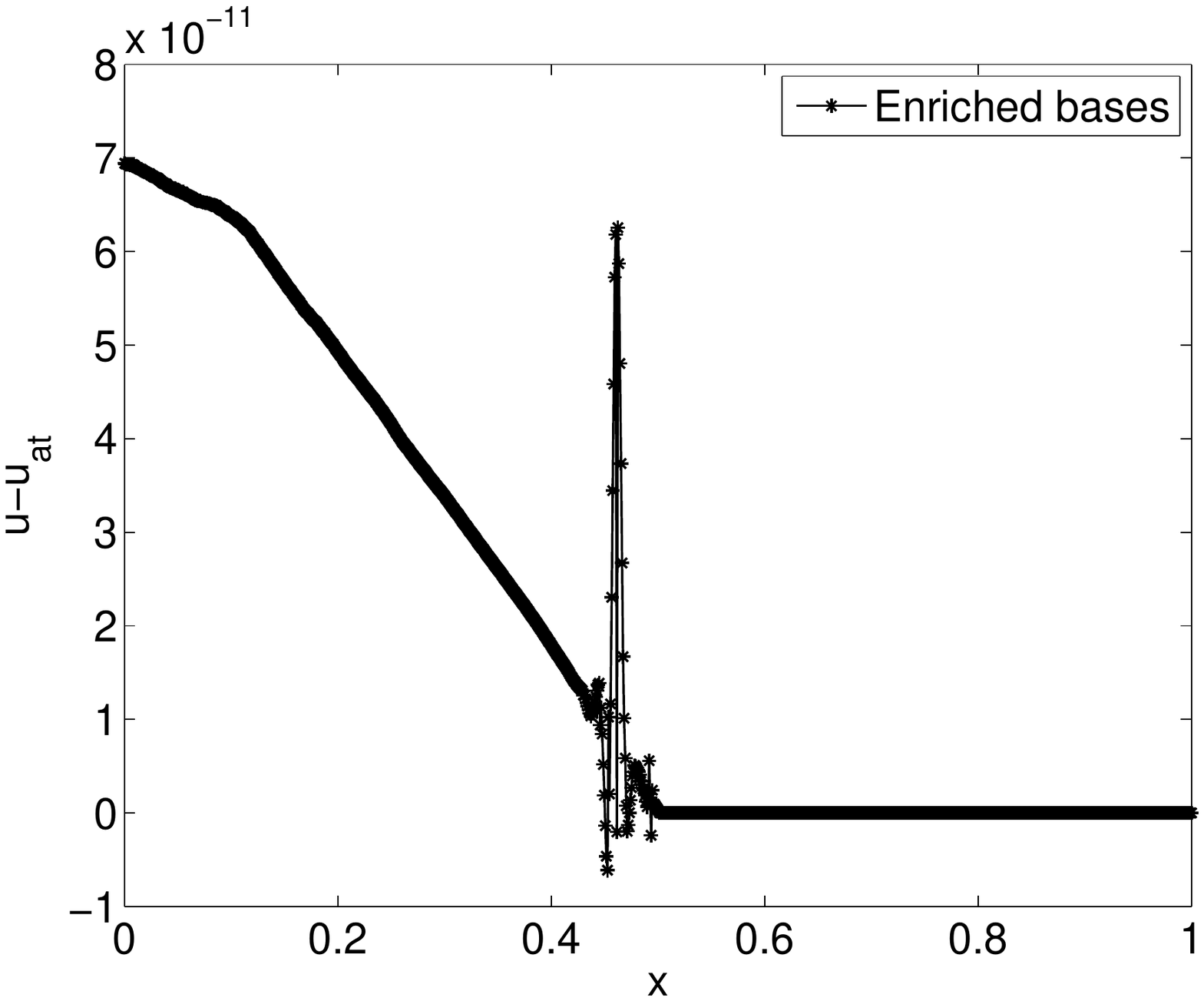}}%
\subfigure[$\ell=3$]{\label {fig:DisplacementErrorCrack4}
\includegraphics[width=2.5in]{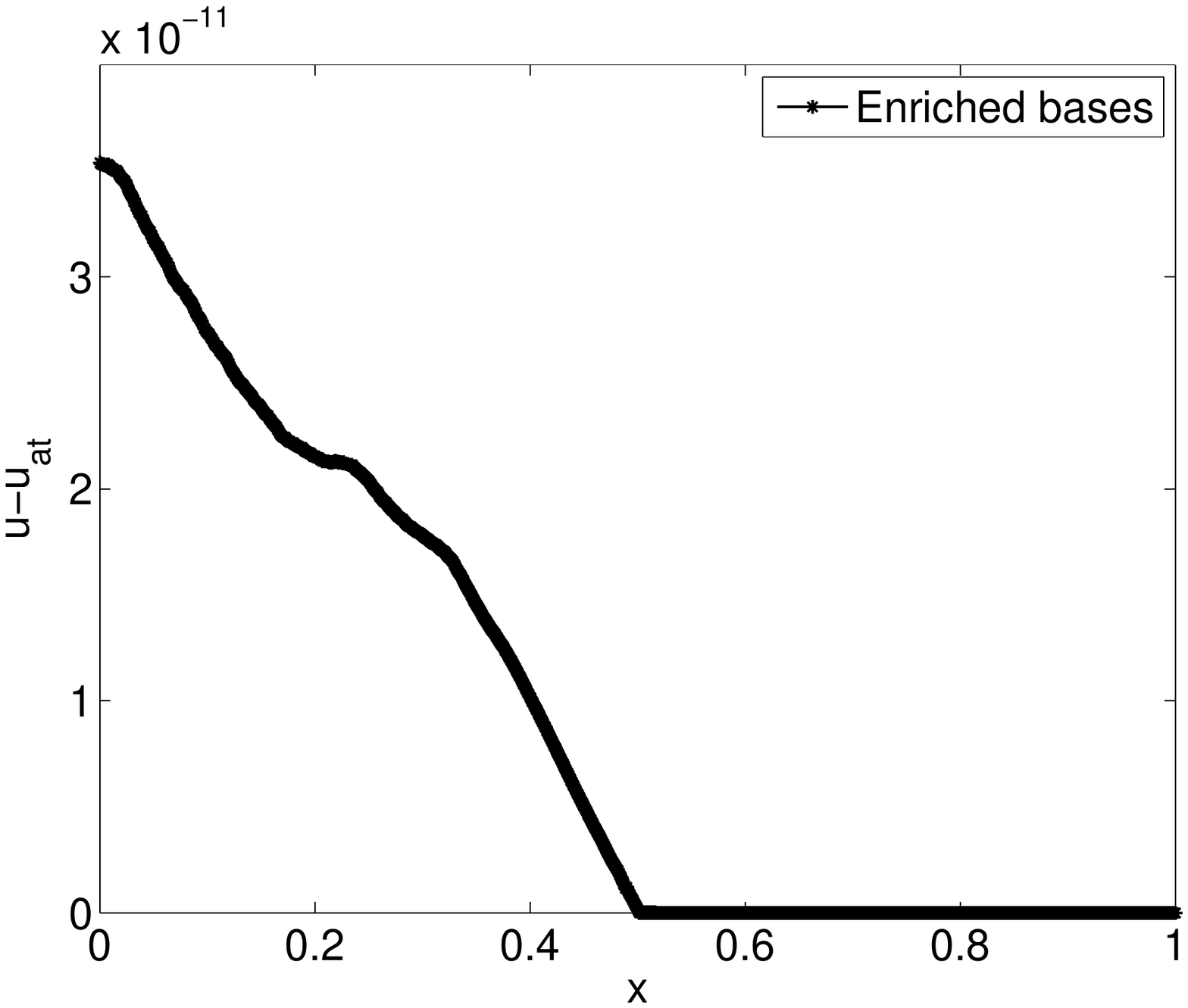}}%
\vspace{-4em}
\caption{\small Displacement errors of Galerkin methods.
(a): The standard Galerkin method; (b): The enriched bases method with $\ell=1$; (c): The
enriched bases method with $\ell=2$; (d): The enriched bases method with $\ell=3$. }\label {fig:DisplacementErrorCrack}
\end{figure}

We first study the numerical error introduced in Galerkin methods.
Consider the system with $1024$ atoms and the crack tip at $n=514$.
The left $512$ atoms are coarse-grained and the remaining atoms are kept to define the atomistic region.
Atomic displacement and displacement gradient are shown in Fig. \ref{fig:DisplacementCrack}. This will
be regarded as the exact solution.
Meanwhile, Fig. \ref{fig:DisplacementErrorCrack} shows displacement errors of
the standard Galerkin method on a nonuniform mesh, and the enriched bases method as a function of $\ell$
for $\ell=1, 2, 3$, respectively. We see that the error is reduced dramatically in the enriched bases method.


\medskip

We now turn to the study of convergence rates. In particular,
we compare the standard Galerkin method with a nonuniform
mesh, the quasi-nonlocal QC method, the force-based QC method,
and the enriched bases method with a nonuniform mesh. The limiting process is defined as follows: All the meshes considered here are fixed as the system
size is enlarged. All atoms are chosen as rep-atoms in the quasi-nonlocal QC method and the force-based QC method.

\begin{figure}[htbp]

\vspace{-4em}
\centering
\subfigcapskip -3em
\subfigure[$W^{1,1}$]{\label {fig:ErrorConvergenceCrackW11}
\includegraphics[width=1.6in]{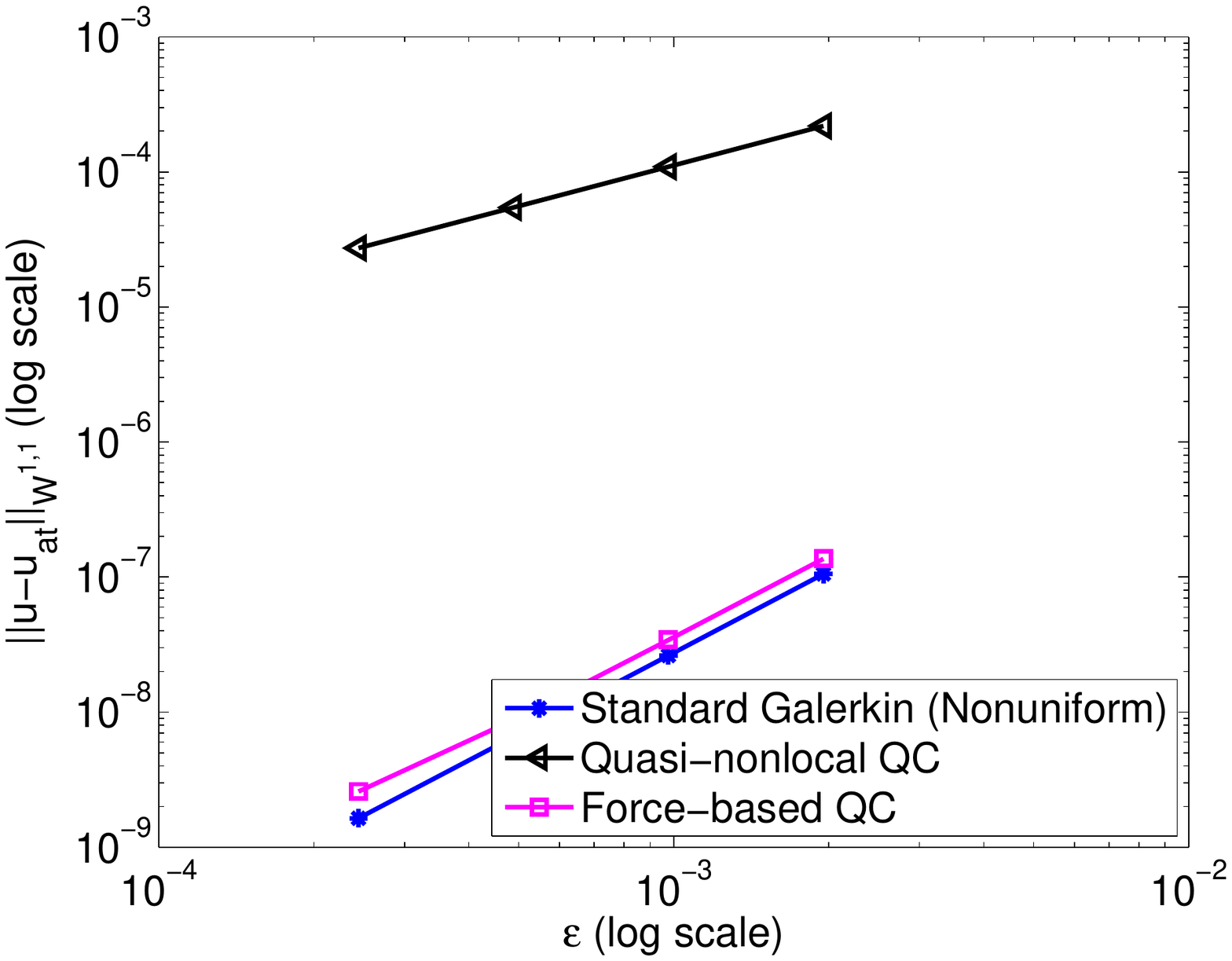}}%
\subfigure[$H^{1}$]{\label {fig:ErrorConvergenceCrackH1}
\includegraphics[width=1.6in]{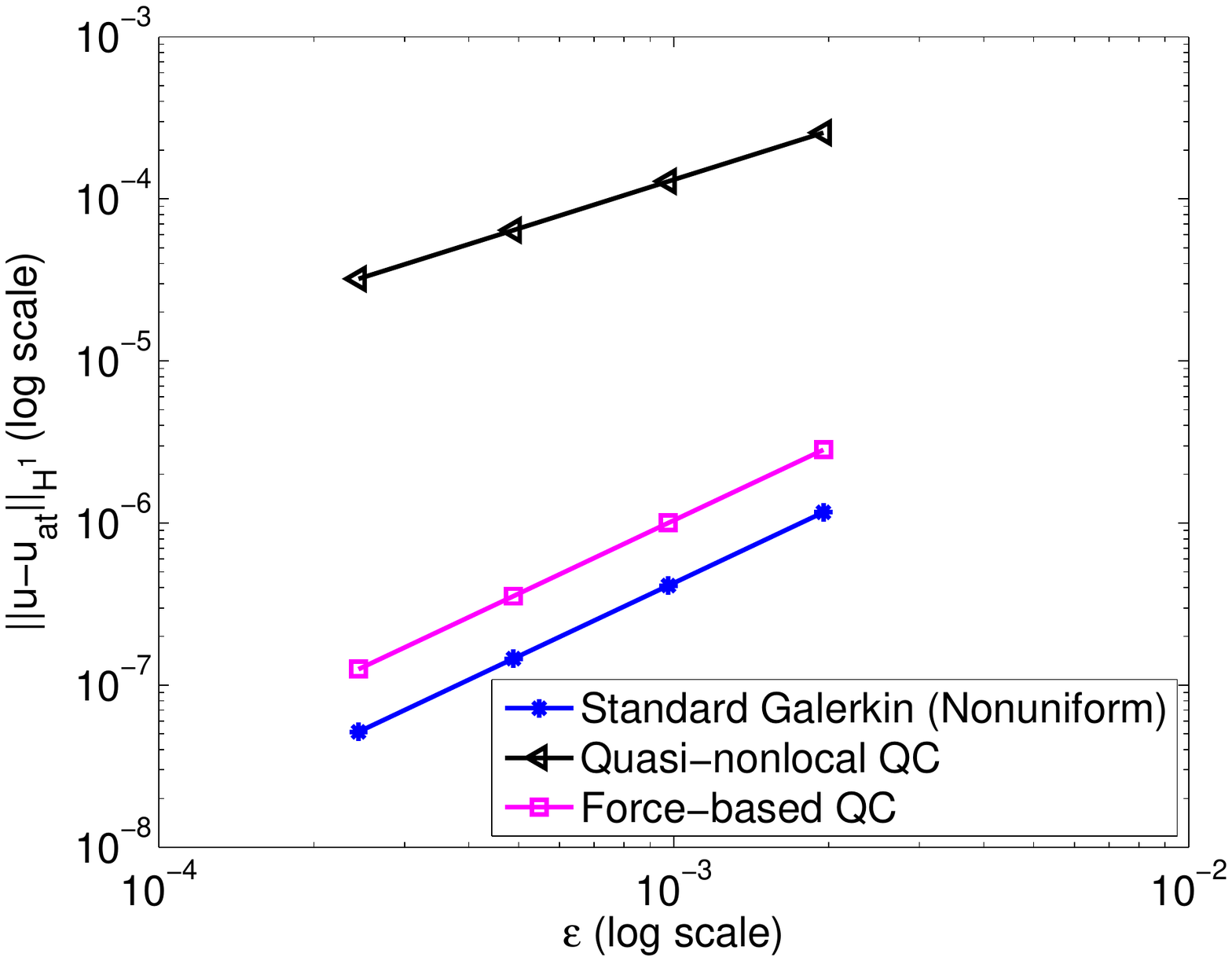}}%
\subfigure[$W^{1,\infty}$]{\label {fig:ErrorConvergenceCrackW1infty}
\includegraphics[width=1.6in]{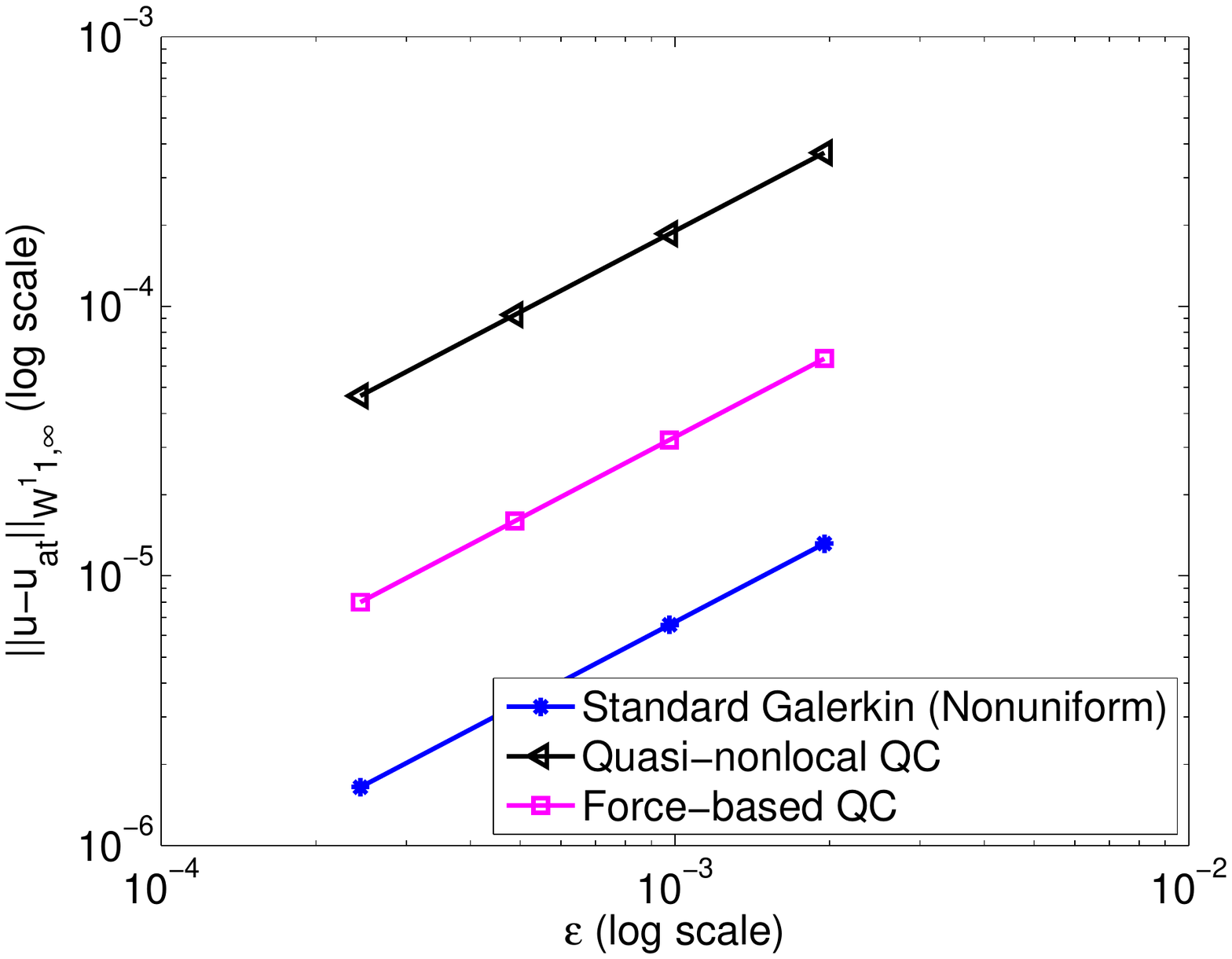}}%
\vspace{-2em}
\caption{\small Errors of the standard Galerkin method with a nonuniform mesh, the quasi-nonlocal QC method,
and the force-based QC method. (a): $W^{1,1}$ norm. Convergence rates are $2.00$, $1.00$ and $1.91$, respectively; (b): $H^1$ norm. Convergence rates are $1.50$, $1.00$ and $1.50$, respectively; (c): $W^{1,\infty}$ norm. Convergence rates are $1.00$, $1.00$ and $1.00$, respectively.}\label {fig:ErrorConvergenceCrack}
\end{figure}

Fig. \ref{fig:ErrorConvergenceCrack} shows convergence rates of the first two methods in $W^{1,1}$ norm, $H^1$ norm, and $W^{1,\infty}$
norm, respectively. The standard Galerkin method on a nonuniform mesh
and the force-based QC method converge
in $W^{1,1}$, $H^1$, and $W^{1,\infty}$ norms with rates of convergence
$2$, $1.5$, and $1$, respectively.
The quasi-nonlocal QC method converges linearly in all three norms.
Moreover, the standard Galerkin method on a nonuniform mesh converges with a prefactor which is more than
two orders of magnitude smaller than that of the quasi-nonlocal QC method.
Compared with the previous example, the external force in this example is also in $L^2$, but
the convergence rates are different. Therefore, how the theoretical analysis can be extended to this case
will be of great interest.

Table \ref{tab:extendeduniformcrack} shows $W^{1,1}$, $H^1$, and $W^{1,\infty}$ norms of the error for the enriched bases method with a nonuniform mesh. With a fixed number of enriched bases around the interface, the error remains small, but slightly grows as $\veps$ becomes smaller, indicating that more enriched bases are needed.

\begin{table}

\centering
\begin{tabular}{|c|c|c|c|c|}
\hline
$\veps$  & $1/512$ & $1/1024$ & $1/2048$ & $1/4096$ \\ \hline
$W^{1,1}$  &  $7.44e-17$  & $1.09e-17$ & $3.00e-17$ & $9.29e-17$  \\ \hline
$H^1$  &  $8.99e-17$  & $1.99e-17$ & $4.19e-17$ & $1.19e-16$  \\ \hline
$W^{1,\infty}$  &  $1.94e-16$  & $1.23e-16$ & $1.4e-16$ & $3.33e-16$  \\ \hline
$m$ &  $18$  & $18$ & $18$ & $18$\\ \hline
\end{tabular}
\caption{$W^{1,1}$, $H^1$, and $W^{1,\infty}$ norms of the error for the enriched bases method
with a uniform mesh. With a fixed number of enriched bases $18$ around
the interface, the error remains small, but slightly
grows as $\veps$ becomes smaller, indicating that more enriched bases are needed.}
\label{tab:extendeduniformcrack}
\end{table}

\begin{figure}[htbp]

\vspace{-6em}
\centering
\subfigcapskip -5em
\subfigure[Standard Galerkin]{\label {fig:DisplacementCrackCAC1}
\includegraphics[width=2.5in]{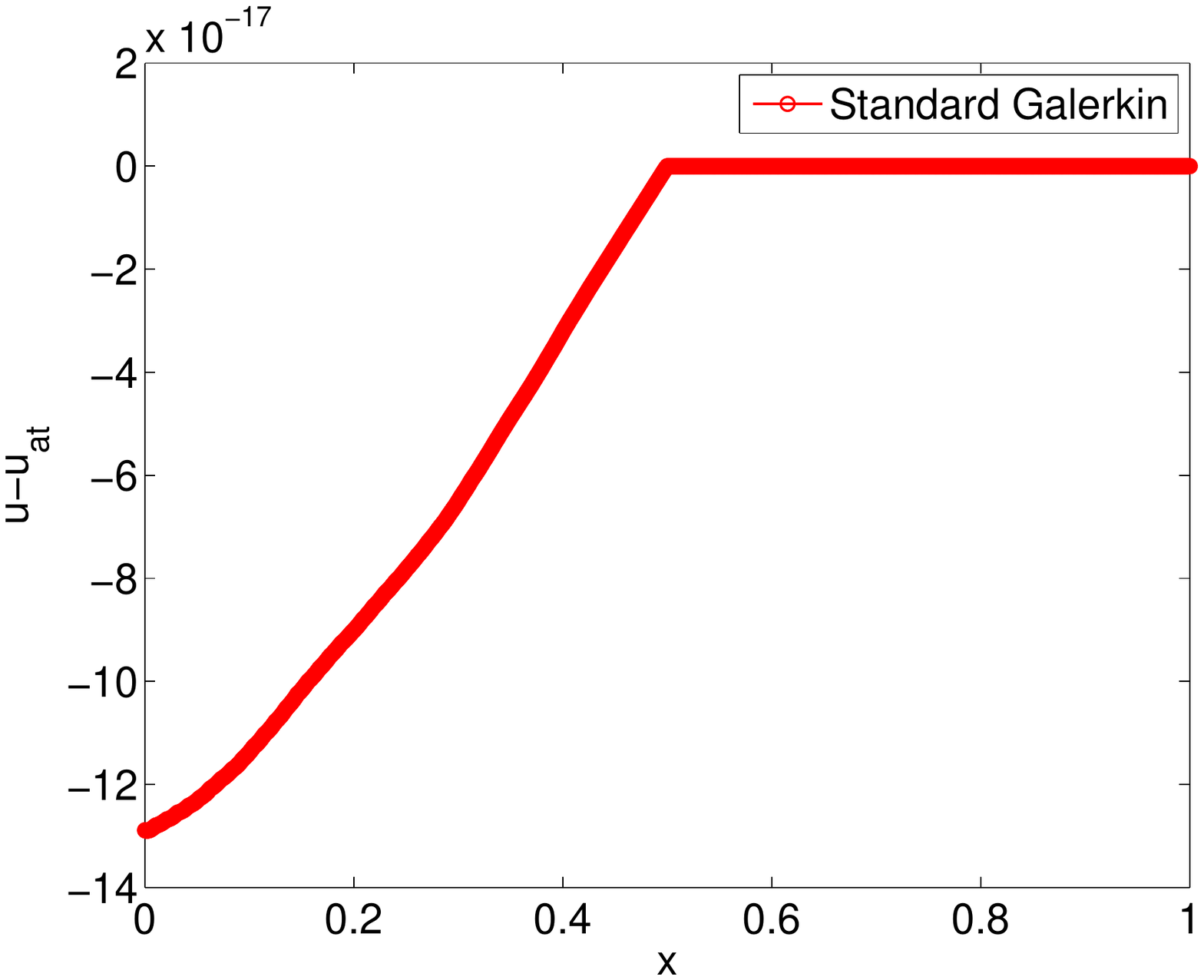}}%
\subfigure[Enriched bases]{\label {fig:DisplacementCrackCAC2}
\includegraphics[width=2.5in]{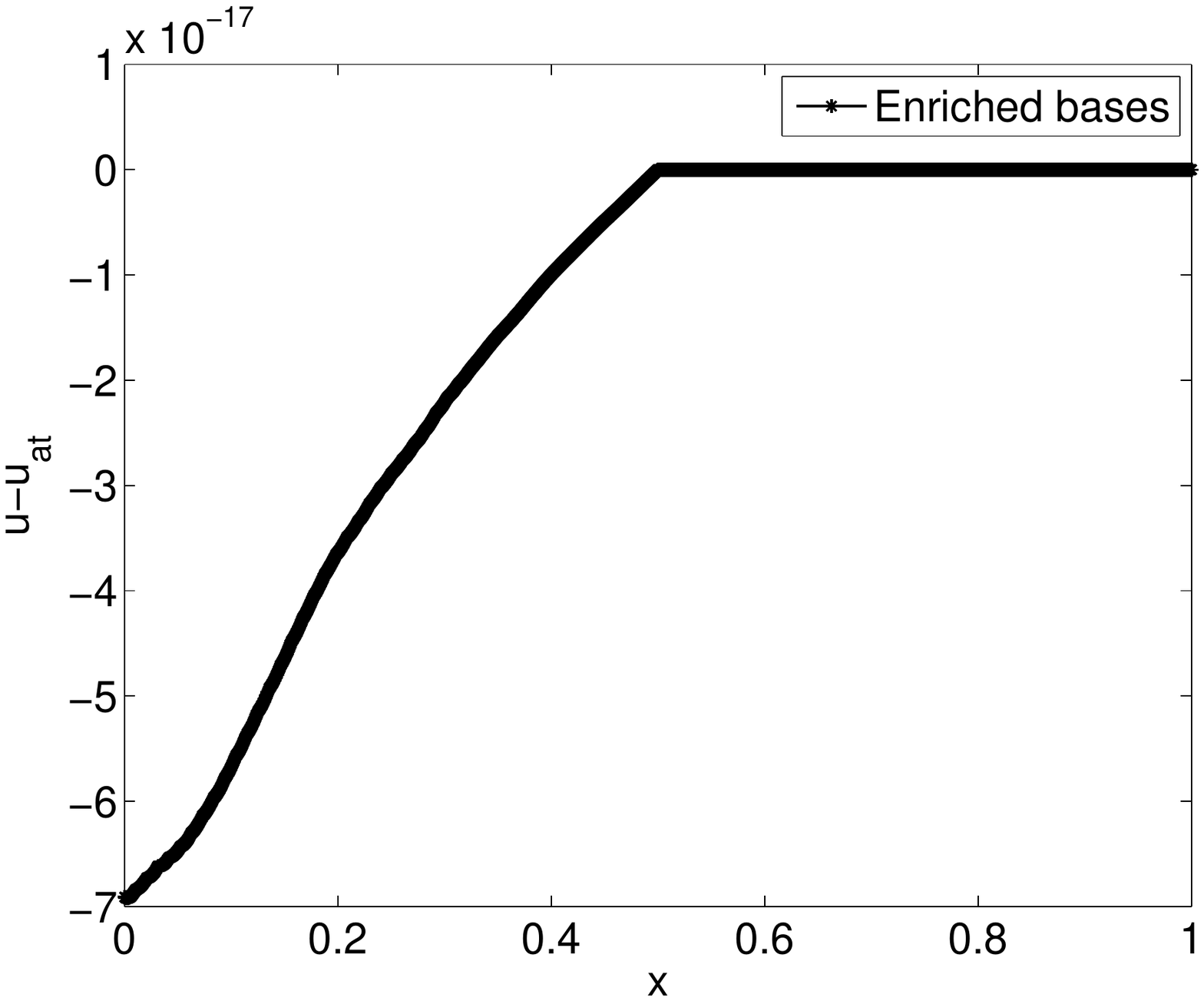}}%
\vspace{-4em}
\caption{\small Displacement errors of Galerkin methods using a continuum/interbedded/atomistic/interbedded/continuum coupling for $384$ atoms. (a): Standard Galerkin method; (b): Enriched bases method.}\label {fig:DisplacementCrackCAC}
\end{figure}

From a practical viewpoint, only a small fraction of atoms around localized defects should be kept,
and all remaining atoms should be coarse-grained. So in the last test, we consider the case where only a few atoms around
the crack tip are kept and all remaining atoms are coarse grained. This
gives a   \protect{continuum/interbedded/atomistic/interbedded/continuum} partition of the system. Fig. \ref{fig:DisplacementCrackCAC} shows displacement errors of Galerkin methods.
It is clear that this coupling, which involves much fewer degrees of freedom, offers a good approximation.

\section{Conclusion}
\label{sec:conclusion}

We derived an effective coarse-grained model starting with molecular mechanics
as an exact model. It is worthwhile to emphasize that we do not introduce the continuum model
to begin with. Rather, we work with a pre-selected CG variables.
The exact model can be projected to a subspace, which naturally leads to the standard
Galerkin method. Due to the lack of direct physical significance of the basis functions, the approximation may
not have sufficient accuracy. We proposed to improve the modeling accuracy by using the enriched bases method, which is based on introducing more basis functions at the interface. We have used the Krylov subspaces and the block Lanczos algorithm to automate this procedure. Both methods have been
systematically tested and compared with some other coupling methods.
Although currently there is no theoretical analysis, the numerical experiments
have indicated excellent performances both quantitatively and qualitatively.
The introduction of the quadrature approximation, the interbedded region, and quasiatoms has made the computational
complexity of our methods comparable to other methods but with better performances.
Application of our methods for high dimensional problems will be presented in
future works. Another interesting issue is the performance of the enriched bases method
if enriched bases are used in the energy-based framework of the Galerkin method since
ghost forces appear as a consequence of the quadrature approximation.

\section*{Acknowledgments}
We are grateful to Prof. Pingbing Ming for helpful discussions.
The work of Chen and Garc\'{\i}a-Cervera was supported by National Natural Science
Foundation via grant DMS1217315.


\begin{thebibliography}{10}

\bibitem{AbdulleLinShapeev2012}
A.~Abdulle, P.~Lin, and A.~V. Shapeev, \emph{Numerical methods for
  multilattices}, Multiscale Model. Simul. \textbf{10} (2012), 696--726.

\bibitem{AbdulleLinShapeevII}
A.~Abdulle, P.~Lin, and A.~V. Shapeev, \emph{A priori and a posteriori
  ${W}^{1,\infty}$ error analysis of a {Q}{C} method for complex lattices},
  SIAM J. Numer. Anal. \textbf{51} (2013), 2357--2379.

\bibitem{Blancetal:2002}
X.~Blanc, C.~Le~Bris, and P.-L. Lions, \emph{From molecular models to continuum
  mechanics}, Arch. Rational Mech. Anal. \textbf{164} (2002), no.~4, 341--381.

\bibitem{BornHuang:1954}
M.~Born and K.~Huang, \emph{Dynamical {T}heory of {C}rystal {L}attices}, Oxford
  University Press, 1954.

\bibitem{B77}
A.~Brandt, \emph{Multi-level adaptive solutions to boundary-value problems},
  Math. Comp. \textbf{31} (1977), no.~138, 333--390.

\bibitem{ChenMing:2012}
J.~Chen and P.B. Ming, \emph{Ghost force influence of a quasicontinuum method
  in two dimension}, J. Comput. Math. \textbf{30} (2012), 657--683.

\bibitem{CuiMing:2013}
L.~Cui and P.B. Ming, \emph{The effect of ghost forces for a quasicontinuum
  method in three dimension}, Sci. in China Ser. A: Math. \textbf{56} (2013),
  2571--2589.

\bibitem{DobsonLuskin:2009b}
M.~Dobson and M.~Luskin, \emph{An optimal order error analysis of the
  one-dimensional quasicontinuum approximation}, SIAM J. Numer. Anal.
  \textbf{47} (2009), 2455--2475.

\bibitem{DobsonLuskinOrtner:2010a}
M.~Dobson, M.~Luskin, and C.~Ortner, \emph{Stability, instability, and error of
  the force-based quasicontinuum approximation}, Arch. Ration. Mech. Anal.
  \textbf{197} (2010), 179--202.

\bibitem{E:book}
W.~E, \emph{Principles of multiscale modeling}, Cambridge University Press,
  2011.

\bibitem{EEngquist02}
W.~E and B.~Engquist, \emph{The heterogeneous multi-scale methods}, Commun.
  Math. Sci. \textbf{1} (2002), 87 -- 132.

\bibitem{ELuYang:2006}
W.~E, J.~Lu, and J.Z. Yang, \emph{Uniform accuracy of the quasicontinuum
  method}, Phys. Rev. B \textbf{74} (2006), 214115.

\bibitem{EMing:2005}
W.~E and P.B. Ming, \emph{Analysis of the local quasicontinuum methods},
  Frontiers and {P}rospects of {C}ontemporary {A}pplied {M}athematics, {Li,
  Tatsien and Zhang, P.W.} ({E}ditor), Higher Education Press, World
  Scientific, Singapore, 2005, pp.~18--32.

\bibitem{EMing:2007}
\bysame, \emph{{C}auchy-{B}orn rule and the stability of crystalline solids:
  {S}tatic problems}, Arch. Rational Mech. Anal. \textbf{183} (2007), no.~2,
  241--297.

\bibitem{el2004block}
A.~El~Guennouni, K.~Jbilou, and H.~Sadok, \emph{The block {L}anczos method for
  linear systems with multiple right-hand sides}, Appl. Numer. Math.
  \textbf{51} (2004), no.~2, 243--256.

\bibitem{Ericksen:1984}
J.L. Ericksen, \emph{The {C}auchy and {B}orn hypothesis for crystals}, Phase
  Transformations and Material Instabilities in Solids (M.E. Gurtin, ed.),
  Academic Press, 1984, pp.~61--77.

\bibitem{Gunzburger:2010}
M.~Gunzburger and Y.~Zhang, \emph{A quadrature-rule type approximation to the
  quasi-continuum method}, Multiscale Model. Simul. \textbf{8} (2010),
  571--590.

\bibitem{KnapOrtiz:2001}
J.~Knap and M.~Ortiz, \emph{An analysis of the quasicontinuum method}, J. Mech.
  Phys. Solids \textbf{49} (2001), 1899--1923.

\bibitem{LaTh03}
S.~Larsson and V.~Thomee, \emph{Partial differential equations with numerical
  methods}, Springer, 2003.

\bibitem{Leach01}
A.R. Leach, \emph{Molecular modelling: {P}rinciples and applications}, Prentice
  Hall, 2001.

\bibitem{Li:2010}
X.~Li, \emph{A coarse-grained molecular dynamics model for crystalline solids},
  Int. J. Numer. Meth. Engng \textbf{83} (2010), no.~8-9, 986--997.

\bibitem{LI12a}
X.~Li, \emph{An atomistic-based boundary element method for the reduction of
  the molecular statics models}, Comp. Meth. Appl. Mech. Engreg \textbf{225}
  (2012), 1--13.

\bibitem{XLi2013A}
\bysame, \emph{A bifurcation study of crack initiation and kinking}, Eur. Phys.
  J. B \textbf{86} (2013), 258.

\bibitem{Li:2014}
X.~Li, \emph{Coarse-graining molecular dynamics models using an extended
  galerkin projection method}, Int. J. Numer. Meth. Engng \textbf{99} (2014),
  no.~3, 157--182.

\bibitem{XLi2014}
X.~Li, \emph{Dynamic crack initiation through bifurcation}, Preprint. (2014).

\bibitem{LiLu14}
X.~Li and J.~Lu, \emph{Traction boundary conditions for atomistic models}, In
  preparation. (2014).

\bibitem{Lietal:2013}
X.~Li, M.~Luskin, C.~Ortner, and A.V. Shapeev, \emph{Theory-based benchmarking
  of the blended force-based quasicontinuum method}, Comput. methods Appl.
  Mech. Engrg. \textbf{268} (2014), 763--781.

\bibitem{LiMing:2014}
X.~Li and P.B. Ming, \emph{A study on the quasicontinuum approximations of a
  one-dimensional fractural model}, Multiscale Model. Simul. \textbf{12}
  (2014), 1379--1400.

\bibitem{Lin07}
P.~Lin, \emph{Convergence analysis of a quasi-continuum approximation for a
  two-dimensional material}, SIAM J. Numer. Anal. \textbf{45} (2007), 313--332.

\bibitem{LuMing:2013}
J.~Lu and P.B. Ming, \emph{Convergence of a force-based hybrid method in three
  dimensions}, Comm. Pure Appl. Math. \textbf{66} (2013), 83--108.

\bibitem{LuMing:2014}
\bysame, \emph{Convergence of a force-based hybrid method with planar sharp
  interface}, SIAM J. Numer. Anal. (2014), in press.

\bibitem{LuskinOrtner:2009}
M.~Luskin and C.~Ortner, \emph{An analysis of node-based cluster summation
  rules in the quasicontinuum method}, SIAM J. Numer. Anal. \textbf{47} (2009),
  3070--3086.

\bibitem{LuskinOrtner:2013}
\bysame, \emph{Atomistic-to-continuum coupling}, Acta Numerica \textbf{22}
  (2013), 397--508.

\bibitem{McCormick87}
J.~Mandel, S.~McCormick, and R.~Bank, \emph{Variational multigrid theory},
  ch.~5, pp.~131--188, SIAM, 1987.

\bibitem{MillerTadmor:2002}
R.E. Miller and E.B. Tadmor, \emph{The quasicontinuum method: Overview,
  applications and current directions}, J. Comput. Aided Mater. Des. \textbf{9}
  (2002), no.~3, 203--239.

\bibitem{MillerTadmor:2009}
\bysame, \emph{A unified framework and performance benchmark of fourteen
  multiscale atomistic/continuum coupling methods}, Modelling Simul. Mater.
  Sci. Eng. \textbf{17} (2009), no.~5, 053001--053051.

\bibitem{Ming:2008}
P.B. Ming, \emph{Error estimate of force-based quasicontinuum method}, Commun.
  Math. Sci. \textbf{6} (2008), 1087--1095.

\bibitem{MingYang:2009}
P.B. Ming and J.Z. Yang, \emph{Analysis of a one-dimensional nonlocal
  quasi-continuum method}, Multiscale Model. Simul. \textbf{7} (2009),
  1838--1875.

\bibitem{OdPrRoBa06}
J.T. Oden, S.~Prudhomme, A.~Romkes, and P.~Bauman, \emph{Multi-scale modeling
  of physical phenomena: Adaptive control of models}, SIAM J. Sci. Comput.
  \textbf{28(6)} (2006), 2359--2389.

\bibitem{OrtnerShapeev}
C.~Ortner and A.V. Shapeev, \emph{Analysis of an energy-based
  atomistic/continuum coupling approximation of a vacancy in the 2d triangular
  lattice}, Math. Comp. \textbf{82} (2013), 2191--2236.

\bibitem{OrtnerZhang:2012}
C.~{Ortner} and L.~{Zhang}, \emph{{Construction and sharp consistency estimates
  for atomistic/continuum coupling methods with general interfaces: a 2D model
  problem}}, SIAM J. Numer. Anal. \textbf{50} (2012), 2940--2965.

\bibitem{OrtnerZhang}
C.~Ortner and L.~Zhang, \emph{Energy-based atomistic-to-continuum coupling
  without ghost forces}, Comput. methods Appl. Mech. Engrg. \textbf{279}
  (2014), 29--45.

\bibitem{citeulike:2898408}
M.~Parks, P.~Bochev, and R.~Lehoucq, \emph{Connecting atomistic-to-continuum
  coupling and domain decomposition}, Multiscale Modeling \& Simulation
  \textbf{7} (2008), no.~1, 362--380.

\bibitem{prudhomme2006error}
S.~Prudhomme, P.T. Bauman, and J.T. Oden, \emph{Error control for molecular
  statics problems}, Int. J. Multiscale Com. \textbf{4} (2006), no.~5-6,
  648--662.

\bibitem{prudhomme2009adaptive}
S.~Prudhomme, L.~Chamoin, H.B. Dhia, and P.T. Bauman, \emph{An adaptive
  strategy for the control of modeling error in two-dimensional
  atomic-to-continuum coupling simulations}, Comput. Methods Appl. Mech. Eng.
  \textbf{198} (2009), no.~21, 1887--1901.

\bibitem{saad1980rates}
Y.~Saad, \emph{On the rates of convergence of the {L}anczos and the
  block-{L}anczos methods}, SIAM J. Numer. Anal. \textbf{17} (1980), no.~5,
  687--706.

\bibitem{Saad:2011}
\bysame, \emph{Numerical methods for large eigenvalue problems}, SIAM, 2011.

\bibitem{schaffner2013analytical}
M.~Sch{\"a}ffner and A.~Schl{\"o}merkemper, \emph{About an analytical
  verification of quasi-continuum methods with {$\Gamma$}-convergence
  techniques}, MRS Proceedings, vol. 1535, Cambridge Univ Press, 2013,
  pp.~mmm12--a.

\bibitem{Shapeev:2011}
A.V. Shapeev, \emph{Consistent energy-based atomistic/continuum coupling for
  two-body potentials in one and two dimensions}, Multiscale Model. Simul.
  \textbf{9} (2011), 905--932.

\bibitem{Shapeev:2012}
A.V. Shapeev, \emph{Consistent energy-based atomistic/continuum coupling for
  two-body potentials in three dimensions}, SIAM J. Sci. Comput. \textbf{34}
  (2012), no.~3, B335--B360.

\bibitem{Shenoyetal:1999}
V.~B. Shenoy, R.~Miller, E.~B. Tadmor, D.~Rodney, R.~Phillips, and M.~Ortiz,
  \emph{An adaptive finite element approach to atomic-scale mechanics: {T}he
  quasicontinuum method}, J. Mech. Phys. Solids \textbf{47} (1999), no.~3,
  611--642.

\bibitem{Shimokawa:2004}
T.~Shimokawa, J.J. Mortensen, J.~Schioz, and K.W. Jacobsen, \emph{Matching
  conditions in the quasicontinuum method: {R}emoval of the error introduced at
  the interface between the coarse-grained and fully atomistic region}, Phys.
  Rev. B \textbf{69} (2004), 214104.

\bibitem{TadmorOrtizPhillips:1996}
E.B. Tadmor, M.~Ortiz, and R.~Phillips, \emph{Quasicontinuum analysis of
  defects in solids}, Philos. Mag. A \textbf{73} (1996), 1529--1563.

\bibitem{ThHsRa71}
R.~Thompson, C.~Hsieh, and V.~Rana, \emph{Lattice trapping of fracture cracks},
  J. Appl. Phys. \textbf{42} (1971), 3154--3160.

\bibitem{VanLuskin:2011}
B.~Van~Koten and M.~Luskin, \emph{Analysis of energy-based blended
  quasi-continuum approximations}, SIAM J. Numer. Anal. \textbf{49} (2011),
  no.~5, 2182--2209.

\bibitem{WaLi03}
G.J. Wagner and W.K. Liu, \emph{Coupling of atomistic and continuum simulations
  using a bridging scale decomposition}, J. Comput. Phys. \textbf{190} (2003),
  249 -- 274.

\end{thebibliography}
\end{document}